\documentclass[paper,twocolumn]{geophysics}

\usepackage{amsmath,amsfonts,amssymb,mathrsfs}
\usepackage{xcolor}

\usepackage{graphicx}
\usepackage{ftnright}%footnode to the right column to fit abstract

\graphicspath{{./FIG/}}

\usepackage[colorlinks = true,
            linkcolor = blue,
            urlcolor  = blue,
            citecolor = blue,
            anchorcolor = blue]{hyperref}

\newtheorem{thm}{Theorem}
\newtheorem{algorithm}[thm]{Algorithm}

\DeclareMathAlphabet{\itbf}{OML}{cmm}{b}{it}
\DeclareMathAlphabet\mathbfcal{OMS}{cmsy}{b}{n}

\renewcommand{\tilde}{\widetilde}
\def\RR{\mathbb{R}}
\def\bx{{{\itbf x}}}

\def\bu{{{\itbf u}}}

\def\bg{{\itbf g}}
\def\be{{\itbf e}}
\def\br{{\itbf r}}

\def\bd{{\itbf d}}

\def\bP{{\itbf P}}
\def\bU{{\itbf U}}
\def\bJ{{\itbf J}}

\def\bv{{\bf v}}
\def\bet{{\boldsymbol{\eta}}}
\def\ss{{(s)}}

\def\bV{{\itbf V}}
\def\bR{{\itbf R}}
\def\bZ{{\itbf Z}}
\def\bI{{\itbf I}}
\def\bS{{\itbf S}}
\def\bQ{{\itbf Q}}
\def\bM{{\itbf M}}
\def\bD{{\itbf D}}
\def\bA{\boldsymbol{\cal A}}
\def\caP{\mathcal{P}}
\def\cP{\boldsymbol{\caP}}

\def\cA{{\mathcal A}}

\def\cf{{\scalebox{0.7}[0.6]{f}}}
\def\cfe{{\scalebox{0.7}[0.6]{fe}}}
\def\cN{{\scalebox{0.5}[0.4]{N}}}

\def\bLa{{\boldsymbol{\Lambda}}}
\def\bGa{{\boldsymbol{\Gamma}}}

\def\FWI{{\scalebox{0.5}[0.4]{FWI}}}

\def\RM{{\scalebox{0.5}[0.4]{ROM}}}

\def\cF{\mathcal{F}}

\def\om{\omega}
\def\la{\lambda}
\def\12{{\frac{1}{2}}}

\begin{document}

\renewcommand{\figdir}{FIGS} % figure directory

\title{Waveform inversion via reduced order modeling}

\renewcommand{\thefootnote}{\fnsymbol{footnote}} 

\address{
\footnotemark[1]University of Michigan, Department of Mathematics, Ann Arbor, MI  48109-1043.
\footnotemark[2]Centre  de Math\'ematiques Appliqu\'ees, Ecole Polytechnique, Institut Polytechnique de Paris, 91128
  Palaiseau Cedex, France.
\footnotemark[3]University of Houston, Department of Mathematics, Houston, TX 77204-3008.
\footnotemark[4]Uppsala Universitet, Department of Information Technology, Division of Scientific Computing, 75105 Uppsala, Sweden.}
\author{Liliana Borcea\footnotemark[1], Josselin Garnier\footnotemark[2],
Alexander V. Mamonov\footnotemark[3] and J\"{o}rn Zimmerling\footnotemark[1]\footnotemark[4]}

\footer{ROM waveform inversion}
\lefthead{L. Borcea, et al.}%\lefthead{L. Borcea, J. Garnier, A.V. Mamonov \& J. Zimmerling}
\righthead{ROM waveform inversion}

%\maketitle is unnecessary

%\renewcommand{\eqref}[1]{equation~\ref{#1}}
\renewcommand{\eqref}[1]{\ref{#1}}

\newcommand\m{{\color{black} \ensuremath{N_{\rm s}} }} %\newcommand\m{{\scalebox{0.7}[0.6]{\color{blue} \ensuremath{N_s} }}}
\newcommand\n{{\color{black} \ensuremath{N_{\rm t}} }} %\newcommand\n{{\scalebox{0.7}[0.6]{\color{blue} \ensuremath{N_t} }}}
\newcommand\q{{{\color{black} \ensuremath{n_{\rm iter}} }}}
\renewcommand\k{{{\color{black} \ensuremath{k} }}}%k is clearly not, not sure what he wantsNt_{\rm max}

\begin{abstract}
We introduce a novel approach to waveform inversion, based on a data driven reduced order model (ROM) of the wave operator. The presentation is for the acoustic wave equation, but the approach can be extended to elastic or electromagnetic waves. The data are time resolved measurements of the pressure wave gathered by an acquisition system which probes the unknown medium with pulses and measures the generated waves.  {We propose to solve the inverse problem of velocity estimation by minimizing the square misfit between the ROM  computed from the recorded data and the ROM computed from the modeled data, at the current guess of the velocity.  We give the step by step computation of the  ROM, which depends nonlinearly on the data and yet can be obtained from them in a non-iterative fashion, using efficient methods from linear algebra.  We also explain how to make the ROM robust to data inaccuracy. 
The ROM computation requires the full array response matrix gathered with  colocated sources and receivers. However, we show that the computation can deal with an approximation of this matrix,  obtained from towed-streamer data using interpolation and reciprocity on-the-fly.

While the full-waveform inversion approach of nonlinear least-squares data fitting is challenging without low frequency information, due to multiple minima of the data fit objective function, we show that the ROM misfit objective function has a better behavior, even for a poor initial guess. We also show by an explicit computation of the objective functions in a simple setting that the ROM misfit objective function {has convexity properties}, whereas the least squares 
data fit objective function displays multiple local minima. } 
\end{abstract}

%===========================================================
\section{Introduction}
\label{sec:intro}

We study the inverse problem of velocity estimation from reflection data gathered by an 
array of $\m$ {colocated  sources and receivers}. The  methodology applies to any linear wave 
equation, for sound or vectorial (electromagnetic or elastic)  waves, but for simplicity 
we work with the acoustic wave equation in a medium 
with constant density and unknown wave speed $c(\bx)$.

Let $p^\ss(t,\bx)$ model the pressure wave generated by the $s^{\rm th}$ 
source, for $s = 1, \ldots, \m$. It satisfies the wave equation
\begin{align}
\big[ \partial_t^2 - c^2(\bx) \Delta \big]  p^\ss(t,\bx) &= f'(t) \delta_{\bx_s}(\bx), \quad 
t \in \RR, 
\label{eq:I1} \\
p^\ss(t,\bx) & = 0, \quad t <-t_{\rm f} ,\label{eq:I2}
\end{align}
for $\bx \in \Omega$, a simply connected  domain, with boundary $\partial \Omega$. This domain can arise from the mathematical truncation of the space, since  over the 
finite duration $T$ of the measurements, the waves are not affected by the medium at distances 
exceeding $T \max_{\bx} c(\bx)$.  Thus, we can impose any homogeneous boundary conditions at $\partial \Omega$, for example Dirichlet.

{The right-hand side in equation \eqref{eq:I1} models the point-like source at location $\bx_s$, where $\delta_{\bx_s}(\bx)$ denotes the 
Dirac $\delta(\bx - \bx_s)$, $f(t)$ is the probing pulse and the prime stands for the time derivative.  It is convenient for the analysis to assume that $f(t)$ is an even function, with support in the interval $(-t_{\rm f},t_{\rm f})$. This may not be the case in practice, but  we explain later that data gathered with an arbitrary pulse that is known or can be estimated, can be transformed by simple processing to data for an even pulse $f(t)$. Prior to the excitation the medium is quiescent, as stated in equation \eqref{eq:I2}.}

The inverse problem is {to} find the velocity $c(\bx)$ from the measured array response matrix 
$\boldsymbol{\mathcal{M}}(t)$, with entries 
\begin{equation}
\mathcal{M}^{(r,s)}(t) = p^\ss(t,\bx_r) ,~ 1 \le r,s \le \m, ~ {t \in (-t_{\rm f},T]}. \label{eq:I3}
\end{equation}
{Note that knowing $\boldsymbol{\mathcal{M}}(t) $ requires colocated sources and receivers. }
This is typically not the case in geophysics applications, but the formulation extends, for example, 
to the {towed-streamer} data acquisition. The missing off-diagonal entries of $\boldsymbol{\mathcal{M}}(t)$ 
are obtained from {towed-streamer} data using source-receiver reciprocity on-the-fly, and the diagonal entries 
can be approximated by interpolation.

Common velocity estimation approaches are travel time tomography \cite[]{dines1979computerized} 
and its more general version studied in the mathematics community \cite[]{stefanov2019travel},
linearized, a.k.a. Born inversion \cite[]{clayton1981born},
migration velocity analysis \cite[]{symes1991velocity,sava2004wave} and 
full-waveform inversion  \cite[]{tarantola1984inversion,virieux2009overview}. 
The first three are based on assumptions {such as} the velocity changes slowly on the scale of the wavelength 
(for travel time tomography), or the velocity variations are small (for Born inversion) or there is 
separation of scales between the smooth components of the velocity and the rough part that gives 
the reflectivity of the medium (for migration). 
Full-waveform inversion (FWI) circumvents such assumptions. It is a partial differential equation 
constrained optimization that fits the data with its model prediction, {typically} in the $L^2$ (least-squares) sense. 
The increase in computing power has lead to growing interest in FWI, but  there is a 
fundamental impediment, which manifests especially for high-frequency data. The objective function is 
nonconvex even in the absence of noise \cite[]{gauthier1986two,santosa1989analysis} and displays 
numerous local minima. This issue, which is due to nonlinear (multiple scattering) effects and cycle-skipping, 
makes any gradient based,  local optimization algorithm, unlikely to succeed in the absence of an accurate 
starting guess \cite[]{virieux2009overview}.
	
There are several approaches to mitigate cycle skipping. {For instance,} multiscale methods
pursue a good starting guess by inverting first very low frequency data  \cite[]{bunks1995multiscale}.
However, such data may not  be available and there is no guarantee that what seems a reasonable 
starting guess  will not create cycle skipping issues for high-frequency data. 
Extended modeling approaches \cite[]{symes2008migration} like the differential semblance method 
\cite[]{symes1991velocity,symes1994inversion} and the source-receiver extension method \cite[]{huang2017full},
introduce in a systematic way additional degrees of freedom in the optimization and then use 
some objective function to drive the extended model toward a velocity estimate.
There are also approaches that use a better alternative than the $L^2$ norm for measuring the data misfit
\cite[]{brossier2010data,bozdaug2011misfit,guitton2003robust}. A prominent alternative is the  optimal transport (Wasserstein) 
metric proposed and analyzed  for seismic inversion in \cite[]{EngquistFroese,yang2018application}.

We introduce a different approach to velocity estimation, based on a data driven reduced order model (ROM) 
of the wave operator. The mapping between the measurements defined in equation \eqref{eq:I3} and the ROM is nonlinear and yet, 
it can be calculated efficiently with methods from numerical linear algebra. %\cite[]{golubVanLoan}. 
The main point of the paper is that the objective function given by the ROM misfit has better behavior 
than the FWI objective function, so optimization methods can converge for a poor initial guess.

There is an ever-growing list of data driven ROM approaches to operator inference and dynamical system 
identification \cite[]{brunton2016discovering,peherstorfer2016data}. 
However, they require data that are not available in our inverse problem. They assume 
knowledge of the state of the system, the wave $p^\ss(t,\bx)$ in our case, at a finite set of time 
{instances} and for all $\bx \in \Omega$.  In contrast, {seismic surveys only provide} the measurements  $\boldsymbol{\mathcal{M}}(t)$ of the wave {at the receiver positions}.

The first array data driven ROM for wave propagation was introduced and used in 
\cite[]{druskin2016direct} in one dimension and in 
\cite[]{borcea2018untangling,borcea2019robust,borcea2020reduced} in higher dimensions. 
The ROM in these studies is not for the wave operator, but for the ``propagator" operator which 
maps the wavefield from one {instance} to the next one {and on} a uniform time grid. 
The ROM propagator has proved useful for imaging the reflectivity of a medium
\cite[]{druskin2018nonlinear,borcea2020reduced,borcea2021reduced}.
In this paper we introduce another ROM, for the wave operator, which is better suited for velocity 
estimation. In fact, we demonstrate with  explicit computations, carried out for a low-dimensional 
velocity model, that the wave operator ROM misfit objective function %is convex.
{has convexity properties.}
 This is not the case 
for the FWI  misfit objective function, computed for the same velocity model. For high-dimensional 
models, where it is not possible to display the objective function, we show via numerical simulations 
that the ROM-based inversion converges to a good estimate of $c(\bx)$, even for a poor initial guess, 
whereas FWI does not. 

\section{Theory}
\label{sect:Meth}

We begin with a general description that motivates our ROM based approach to velocity estimation and 
gives the key ideas behind the ROM construction. Then, {we discuss the mathematical details that establish 
the relationship between the ROM and the data and we summarize the ROM construction in the form of an algorithm}.
 The methodology introduced  in this section {assumes} noiseless data and 
 full knowledge of the array response matrix $\boldsymbol{\cal M}(t)$. This allows us to {describe the objective 
 function for velocity estimation without using a penalty regularization term.} However, regularization is important 
and must be done carefully, as explained later in  the paper, {when dealing} with noisy data and the 
approximation of $\boldsymbol{\cal M}(t)$ from towed-streamer type of measurements.

\subsection{Outline and motivation of the method}

The FWI approach to velocity estimation seeks  {an approximate inverse of the nonlinear forward map
$c(\bx) \stackrel{\mathcal{F}}{\mapsto} \boldsymbol{\cal M}(t)$} using the data misfit minimization
\begin{equation}
\min_{v \in \mathcal{C}} \int_{{-t_{\rm f}}}^{T} dt \, \|\boldsymbol{\cal M}(t) - \mathcal{F}[v](t) \|_F^2 
+ \mbox{regularization},
\label{eq:FWI_obj}
\end{equation}
where $v$ denotes the search velocity in the search space $\mathcal{C}$ and 
$\| \cdot \|_F$ is the matrix Frobenius norm. Our approach introduces an additional mapping,
from $\boldsymbol{\cal M}(t)$ to an approximation of the symmetrized wave operator
$\partial_t^2 + \cA$. The symmetrization is carried out with a similarity transformation of the 
usual wave operator $\partial_t^2 - c^2(\bx) \Delta$. It amounts to scaling  $p^\ss(t,\bx)$ by $c^{-1}(\bx)$ 
and gives
\begin{equation}
\cA = c^{-1}(\bx) \left[ -c^2(\bx) \Delta \right] c(\bx) = - c(\bx) \Delta \big[ c(\bx) \cdot \big] 
\label{eq:defA}.
\end{equation} 

The approximation that we seek is the ROM wave operator $\partial_t^2 + \bA^\RM$, 
where $\bA^\RM$ is a symmetric and positive definite matrix, a Galerkin approximation of 
the self-adjoint and positive definite operator $\cA$. Roughly speaking, we can think of the 
data to ROM mapping $\mathcal{R}$ as a preconditioner of the forward mapping $\mathcal{F}$
\begin{equation}
c(\bx) \stackrel{\mathcal{F}}{\mapsto} \boldsymbol{\cal M}(t) \stackrel{\mathcal{R}}{\mapsto} \bA^\RM,
\end{equation}
because the composition $\mathcal{R} \circ \mathcal{F}$, which gives $\bA^\RM = \mathcal{R} \big(\mathcal{F}[c]\big)$, 
is easier to ``invert".  

The Galerkin method is a standard way of approximating an operator, like $\cA$, by a matrix. 
Typically, the approximation is in spaces of piecewise polynomial functions {with support over a few grid cells}
\cite[]{brenner2008mathematical}. {If we  gather these functions in a row vector field $\boldsymbol{\Psi}(\bx)$, 
the {matrix approximation of $\cA$ is} 
\begin{equation}
\label{eq:APsi}
\bA^\Psi = \int_\Omega d \bx \, \boldsymbol{\Psi}^T(\bx) \cA \boldsymbol{\Psi}(\bx).
\end{equation}
{This matrix $\bA^\Psi$ has a much simpler dependence on $c(\bx)$ than $\boldsymbol{\cal M}(t) = \cF[c](t)$, because its 
 entries depend quadratically on the coefficient $c(\bx)$ integrated locally, in a few grid cells}. 
{It would be easy to find $c(\bx)$ from  $\bA^\Psi$, but this matrix cannot be computed
from the measurements $\boldsymbol{\cal M}(t)$.} 

Our ROM matrix $\bA^\RM$ is a Galerkin approximation of $\cA$ on the space spanned by the snapshots 
of the wavefield, at $\n$ discrete and equidistant time {instances}. Such approximation spaces are common 
in model order reduction \cite[]{brunton2019data,hesthaven2016certified}, where  the idea is to use the 
history of the wavefield to extrapolate or interpolate its behavior. Our projection of $\cA$ is carried out using an 
orthonormal basis of the space of snapshots, gathered in the row vector field $\bV(\bx)$,
\begin{equation}
\bA^\RM = \int_{\Omega} d \bx \, \bV^T(\bx) \cA \bV(\bx) \in \RR^{\n\m \times \n\m}.
\label{eq:ROMAProj}
\end{equation}\\

{Here are the important observations about $\bV(\bx)$:}
\begin{enumerate}
\itemsep -0.01in
\item The ROM matrix $\bA^\RM$ can be obtained directly from the measurements $\boldsymbol{\cal M}(t)$,
without knowing {the snapshots ${\itbf V}(\bx)$ nor the operator  ${\cal A}$.
This is one of the most striking results of this paper.
We summarize the transform ${\cal R}$ from $\boldsymbol{\cal M}(t)$ to $\bA^\RM$  in Algorithm~\ref{alg:arom}
% and in the flow chart in Figure~\ref{Fig.FC}
and we explain the relationship between the ROM and the data that leads to Algorithm~\ref{alg:arom}  in the next subsection.}
\item {$\bV(\bx)$ cannot be computed from the measurements. However, the analysis in   \cite[~Appendix A]{borcea2021reduced} and numerical studies in \cite[~Section 6.3]{borcea2021reduced}  suggest  that $\bV(\bx)$ is almost 
independent of the rough part of $c(\bx)$ i.e., the reflectivity.} 
\item The basis functions {in $\bV(\bx)$} associated with the $j^{\rm th}$ time 
{instance} are peaked near the maximum depth reached by the wavefield up to this {instance}. 
\item  $\bV(\bx)$ is causal. With the first $\k <\n$ snapshots, {the definition in equation~\eqref{eq:ROMAProj}} gives the principal $\k\m \times \k \m$ submatrix 
of $\bA^\RM$, 
obtained by removing its last $(\n-\k)\m$ rows and columns. 
\end{enumerate}

Since $\bV(\bx)$  depends on $c(\bx)$ in a complicated way,  
we cannot prove the convexity of the ROM misfit objective function 
$v \mapsto \|\bA^\RM - \mathcal{R} \big(\mathcal{F}[v]\big)\|_F^2$ for a general medium. 
It is only in layered media that the result follows from the proof in 
\cite[~Appendix A]{borcea2021reduced}. {Explicitly, it is proved there that in a layered medium with variable wave speed and density,
containing multiple reflectors of arbitrary strength,  the orthonormal basis written in travel time coordinates is almost the same as the one in a homogeneous medium. This means that at least in the vicinity of the right kinematics, the dependence of $\bA^\RM$ on $c(\bx)$ is mainly through $\cA$, 
and the objective function is locally convex.}  

{
In general media we expect that, for a rich enough space of snapshots,
which allows a good approximation of $ \boldsymbol{\Psi}(\bx)$ in equation \eqref{eq:APsi} in terms of $\bV(\bx)$, 
the ROM matrix $\bA^\RM$ contains roughly the same information as $\bA^\Psi$.  The numerical study in 
\cite[~Section 6.3]{borcea2021reduced} shows that  ``rich enough" means 
for sources/receivers separated by roughly 
half a wavelength and for time sampling satisfying the Nyquist criterium. }
The {third} attribute of $\bV(\bx)$ listed above and equation \eqref{eq:ROMAProj} also show that the 
entries of $\bA^\RM$ depend mostly on the locally integrated $c(\bx)$, similar to $\bA^\Psi$. 
Thus, we expect that the velocity estimation from the computable $\bA^\RM$ behaves similarly to that 
from the uncomputable $\bA^\Psi$, which is why we propose using the minimization
\begin{equation}
 \min_{v \in \mathcal{C}} \|\bA^\RM - \mathcal{R} \big(\mathcal{F}[v]\big)\|_F^2 + \mbox{regularization}.
\label{eq:ROM_obj}
\end{equation}

{
 The minimization problem (\ref{eq:ROM_obj}) can be solved with a Gauss-Newton iterative method that is summarized in Algorithm~\ref{alg:prowi}. But first, we explain  the relationship between the ROM and the data.
}
 
%{
%One of the most striking results of the paper is that we can obtain $\bA^\RM$ defined by (\ref{eq:ROMAProj}) directly from $\boldsymbol{\cal M}(t)$, without knowing the snapshots ${\itbf V}(\bx)$ nor the operator  ${\cal A}$. 
%We summarize the transform ${\cal R}$ from $\boldsymbol{\cal M}(t)$ to $\bA^\RM$  in Algorithm~\ref{alg:arom} and in the flow chart in Figure~\ref{Fig.FC}.
% The minimization problem (\ref{eq:ROM_obj}) can then be solved by a Gauss-Newton method that is summarized in Algorithm~\ref{alg:prowi}.
%}
%
%{ 
%We explain the relationship between the ROM and the data that leads to Algorithm~\ref{alg:arom}  in the next subsection. We then study the ROM based velocity estimation method that gives Algorithm~\ref{alg:prowi}.
% }

%The question is how can we obtain $\bA^\RM$ directly from $\boldsymbol{\cal M}(t)$, without knowing 
%the snapshots and therefore the Galerkin space? To achieve this, {we write the wave field 
%in a convenient expression involving $\bA$. Then, we  manipulate this expression using operator calculus,
% to get $\bA^\RM$.}

\subsection{{Relationship between the ROM and the data}}
{We begin by transforming} equation \eqref{eq:I1} 
to a homogeneous wave equation for a new wave $u^\ss(t,\bx)$, with {an} initial state determined by the source. 
This new wave is defined in the next section  and the transformation involves {working with the  
even in time wave 
\begin{equation}
\label{eq:evenp}
p_e^\ss(t,\bx) = [p^\ss(t,\bx) + p^\ss(-t,\bx)],
\end{equation} 
where $p^\ss(t,\bx)$ solves equations \ref{eq:I1} and \ref{eq:I2}.}
We can think of the transformation as a 
Duhamel principle, although it is not in the usual form \cite[]{FritzJohn}, because at $t = 0$ we get  
\begin{equation}
u^\ss(0,\bx) = u_0^\ss(\bx), \quad \partial_t u^\ss(0,\bx) = 0, \quad \bx \in \Omega,
\label{eq:INIC}
\end{equation}
with $u_0^\ss(\bx)$ determined by the source location $\bx_s$ and the probing pulse $f(t)$.

{Note that we do not lose any information by working with the even wave in equation \ref{eq:evenp} and therefore the simple initial conditions in equation \eqref{eq:INIC}, 
as long as we know the medium near the colocated sources/receivers. Near means within the distance of travel over the small time interval  
$(-t_{\rm f},t_{\rm f})$ of support of  $f(t)$. We assume henceforth that the medium near the colocated sources/receivers is known and homogeneous, 
with velocity $\bar{c}$.  Due to the initial condition in equation \eqref{eq:I2},  we observe that 
\begin{equation}
p_e^\ss(t,\bx_r) = p^\ss(t,\bx_r), \qquad \mbox{for} ~ t \geq t_{\rm f}, ~~ s, r = 1, \ldots, N_s.
\end{equation}
The waves differ at $t \in {[0,t_{\rm f})}$,  but  since for such time the measurements are 
insensitive to the unknown part of the medium, no information is lost. }

{Note also that the measurements $p_e^\ss(t,\bx_r)$ are obtained easily from those of $p^\ss(t,\bx_r)$, if the latter are gathered at $t \ge -t_{\rm f}$,
for $s,r = 1, \ldots, N_s$, as assumed in equation \ref{eq:I3}. But even if the measurements are made at $t \geq t_{\rm f}$ only, we can still compute $p_e^\ss(t,\bx_r)$ at $t \in { [0,t_{\rm f})}$ by solving the wave equation with velocity $\bar{c}$. Thus, in either case,}  we can map the measured $\boldsymbol{\mathcal{M}}(t)$ to a new 
data matrix $\bD(t)$, with entries {at $t \ge 0$ given by}
\begin{align}
D^{(r,s)} (t)&=p^\ss(t,\bx_r)+p^\ss(-t,\bx_r)\nonumber \\
&=\mathcal{M}^{(r,s)}(t)+\mathcal{M}^{(r,s)}(-t), ~~1 \le r,s \le \m.
\label{eq:defD}
\end{align}

To define our Galerkin approximation space, let us consider a time discretization $t_j = j \tau,$  
with uniform stepping $\tau$, for $ j \ge 0$. We gather the waves $u^\ss(t,\bx)$ evaluated at $t_j$, for all the $\m$ sources, 
in the $j^{\rm th}$ snapshot vector field 
\begin{equation}
\bu_j(\bx) = \left(u^{(1)}(t_j,\bx), \ldots, u^{(\m)}(t_j,\bx) \right), \quad \bx \in \Omega.
\label{eq:snapshot}
\end{equation}
We are interested only in the first $\n$  snapshots, and organize them in  the $\n\m $ dimensional row vector field 
\begin{equation}
\bU(\bx) = \left(\bu_0(\bx), \ldots, \bu_{\n-1}(\bx) \right), \quad \bx \in \Omega.
\label{eq:defmatU}
\end{equation}
The space spanned by the components of $\bU(\bx)$, denoted  $\mbox{range} \big(\bU(\bx) \big)$, 
is our approximation space and the Galerkin approximation of the wavefield is  
\begin{equation}
\hspace{-0.02in} \bu_{\rm G}(t, \bx) =  \Big(u^{(1)}_{\rm G}(t,\bx), \ldots, u^{(\m)}_{\rm G}(t,\bx) \Big) = \bU(\bx) \bg(t)
\label{eq:GalAp}
\end{equation}
with time dependent coefficients gathered in the matrices $\bg(t) \in \RR^{\n\m \times \m}$. 
These coefficients are such that when substituting equation \eqref{eq:GalAp}
into the homogeneous wave equation, the residual is orthogonal to the approximation space. 
This gives the following system of second order ordinary differential equations
\begin{equation}
\underbrace{\int_{\Omega} d \bx \, \bU^T(\bx) \bU(\bx)}_{\bM} \bg''(t) + 
\underbrace{ \int_{\Omega} d \bx \, \bU^T(\bx) \cA \bU(\bx)}_{\bS} \bg(t) = 0,
\label{eq:Gal1}
\end{equation} 
for $t > 0$, with initial condition
\begin{equation}
\bg(0) = \be_0, \quad \bg'(0) = {\bf 0}.
\label{eq:Gal2}
\end{equation}
Here  $\be_0$ is the first $\n\m \times \m$  block of the $\n\m \times \n\m$ identity matrix $\bI_{\n\m}$. 
Equation \eqref{eq:Gal2} ensures that the Galerkin approximation \eqref{eq:GalAp} satisfies 
the initial conditions
\begin{equation}
\bu_{\rm G}(0, \bx) = \bU(\bx) \be_0 = \bu_0(\bx), \quad \partial_t \bu_{\rm G}(0, \bx) = {\bf 0}.
\end{equation}

The Galerkin approximation described above would be {straightforward} if we knew $\bU(\bx)$, 
{but we do not know it}. 
Our key observation is that the $\n\m \times \n\m$ Gramian matrix
\begin{equation}
\bM =  \int_{\Omega} d \bx \, \bU^T(\bx) \bU(\bx) \in \RR^{\n\m \times \n\m},
\label{eq:defM}
\end{equation}
called the ``mass matrix" in Galerkin jargon, and the ``stiffness matrix" 
\begin{equation}
\bS = \int_{\Omega} d \bx \, \bU^T(\bx) \cA \bU(\bx)\in \RR^{\n\m \times \n\m},
\label{eq:defS}
\end{equation}
can be calculated directly from  $\bD(t)$ and the second derivative  $\bD''(t)$, 
evaluated at instances $\{t_j = j \tau\}_{j=0}^{2\n-2}$, as explained in the next section (see Appendix \ref{app:numdata}
for the estimation of $\bD''(t)$,  using a filtered Fourier transform).
Thus, even though we do not know the operator $\cA$ and the vector field $\bU(\bx)$, we can compute the Galerkin 
coefficients $\bg(t)$ for all $t \ge 0$, by solving the system of equations \eqref{eq:Gal1} 
with the data driven $\bM$ and $\bS$, and the initial conditions given in equation \eqref{eq:Gal2}.

The final step of the ROM construction is to put equation \eqref{eq:Gal1} in an algebraic form
that describes the evolution of a causal wave $\bu^\RM(t) \in \RR^{\n\m \times \m}$. 
Each column of this wave corresponds to a source index $s$, with $1 \le s  \le \m$. 
Initially, the true wave is  supported near the sources, which is reflected in the algebraic 
structure of $\bu^\RM(0)$, whose only nonzero entries are in the first $\m \times \m$ block. 
At later times there is  block row fill-in in $\bu^\RM(t)$, which models   wave propagation 
further away from the sources. 

The desired transformation of equation \eqref{eq:Gal1} is achieved using the block Cholesky 
square root \cite[]{golubVanLoan} of the data driven mass matrix
\begin{equation}
\bM = \bR^T \bR,
\label{eq:Cholesky}
\end{equation}
where $\bR$ is block upper triangular (with blocks of size $\m\times \m$). 
The wave in the ROM space is defined by 
\begin{equation}
\bu^\RM(t) = \bR \bg(t),
\label{eq:ROMwave}
\end{equation}
and we note from equation \eqref{eq:Gal2} that at $t = 0$ it satisfies 
\begin{equation}
\bu^\RM(0) =  \bR \be_0 = \begin{pmatrix} \bR_{0,0} \\ {\bf 0} \\ \vdots \\{\bf 0} 
\end{pmatrix}, \quad \frac{d\bu^\RM}{dt}(0) = {\bf 0},
\end{equation}
where $\bR_{0,0} \in \RR^{\m \times \m}$. The wave equation in the ROM space is 
obtained after multiplying equation \eqref{eq:Gal1} on the left by $\bR^{-T} = (\bR^{-1})^T$, 
\begin{equation}
\frac{d^2\bu^\RM }{dt^2}(t) + \bA^\RM \bu^\RM(t) = 0, \quad t > 0,
\label{eq:ROMeq}
\end{equation}
and the ROM approximation of $\cA$ is the $\n\m \times \n\m$ matrix
\begin{equation}
\bA^\RM = \bR^{-T} \bS \bR^{-1}.
\label{eq:ROMA}
\end{equation}

Note that the same block upper triangular matrix $\bR$ arises in the Gram-Schmidt orthogonalization 
of the components of $\bU(\bx)$ given by
\begin{equation}
\bU(\bx) = \bV(\bx) \bR,
\label{eq:GS}
\end{equation}
where $\bV(\bx)$ is an $\n\m$ dimensional row vector field, with orthonormal components, i.e., it satisfies 
\begin{equation}
\int_{\Omega} d \bx \, \bV^T(\bx) \bV(\bx) = \bI_{\n\m}.
\label{eq:orthV}
\end{equation}
This $\bV(\bx)$ stores the orthonormal basis mentioned earlier in the section. 
Its causality, in the sense that the $j^{\rm th}$ (m-dimensional) component of $\bV(\bx)$ 
is determined by $\bu_0(\bx), \ldots, \bu_j(\bx)$, is built into the Gram-Schmidt 
orthogonalization procedure, and therefore in the block upper triangular structure of  $\bR$. 
Substituting equation \eqref{eq:GS} into equation \eqref{eq:defM}, and using {equation}~\eqref{eq:orthV}, 
we observe that $\bR$ in equation \eqref{eq:GS} is the same as in equation \eqref{eq:Cholesky}, because
\begin{equation}
\bM = \bR^T \int_{\Omega} d \bx \, \bV^T(\bx) \bV(\bx) \bR = \bR^T \bR.
\label{eq:GS1}
\end{equation}

If we use the Gram-Schmidt equation \eqref{eq:GS} in {equation}~\eqref{eq:ROMA}, and recall {equation}~\eqref{eq:defS} {for} $\bS$, we get that $\bA^\RM$ satisfies equation~\eqref{eq:ROMAProj}. 
Therefore, the data driven $\bA^\RM$ defined in  {equation}~\eqref{eq:ROMA}, 
is in fact the orthogonal projection of the operator $\cA$ on the unknown space 
$\mbox{range} \big(\bU(\bx)\big)$, obtained with the unknown causal and orthonormal basis in $\bV(\bx)$. 

{We can now add a fifth observation about $\bV(\bx)$. It has been proved recently in \cite[~Proposition 3.2]{borcea2022waveform} that the snapshots gathered in  $\bU(\bx;v) = \bV(\bx;v) \bR$  satisfy exactly the data $\{\bD(j \tau)\}_{j=0}^{2n-2}$.
The difference between this  field and the true one in equation \eqref{eq:GS} is that the unknown $\bV(\bx)$ is replaced by $\bV(\bx;v)$,
whose components are  the orthonormal basis functions computed with the guess velocity 
$v(\bx)$. Any guess velocity works, even  $v(\bx) = \bar c$. That both $\bU(\bx)$ and $\bU(\bx;v)$ 
give an exact data fit, means that the data driven matrix $\bR$ contains all the information. This is why, as shown in \cite[]{borcea2022waveform},  $\bU(\bx;v)$  contains all the arrival events present in $\bU(\bx)$.  The purpose of $\bV(\bx)$ in equation \eqref{eq:GS} may be viewed as 
mapping the information in $\bR$,  from the algebraic (ROM) space to the physical space. When we have the wrong kinematics (smooth part of $v(\bx)$),  
 $\bV(\bx;v)$ maps the arrivals to incorrect depths.
But if the kinematics is only slightly wrong, the computable $\bV(\bx;v)$ is very close to the uncomputable $\bV(\bx)$. This is another way of explaining that at least close enough to the true velocity, $\bA^{\RM}$ defined by equation \eqref{eq:ROMAProj} depends on $c(\bx)$ mostly through 
$\cA$ and the objective function of the ROM misfit is locally convex.
}

\begin{figure}[h!]
\begin{center}
\includegraphics[width=0.48\textwidth]
{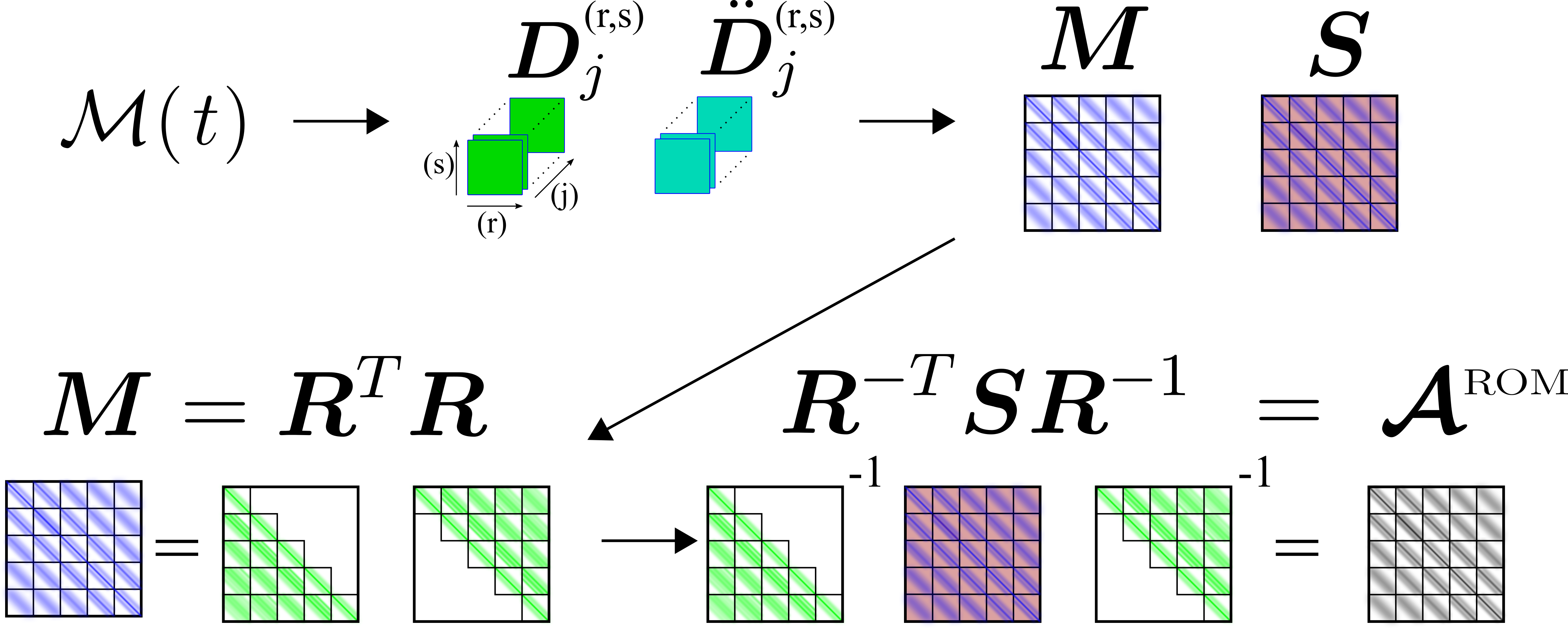}
\end{center}
\caption{Flow chart for the computation of the ROM from the measurements. 
There are four steps, each indicated with an arrow. {All the matrices are of size $\n \m \times \n \m$, with 
entries organized in $\m \times \m$ blocks.}}
\label{Fig.FC}
\end{figure}

\subsection{Technical details of ROM computation}

We show in Figure~\ref{Fig.FC} the flow chart of the computation of $\bA^\RM$ from 
the measurements $\boldsymbol{\cal M}(t)$. The first step computes the data matrices 
\begin{align}
\hspace{-0.1in}
\bD(t) = \left(\mathcal{M}^{(r,s)}(t)+\mathcal{M}^{(r,s)}(-t) \right)_{r,s=1}^{\m}
\label{eq:Dj}
\end{align}
and their second derivatives $\bD''(t)$ at { instances} $t = j \tau$, for $0 \le j \le 2\n-2$.
{Recall from the previous discussion that $\mathcal{M}^{(r,s)}(-t)$ contributes only at $t = j \tau \in [0,t_{\rm f})$
and it may either be measured or computed in the reference medium with velocity $\bar{c}$. }
The details on the computation of the second derivative $\bD''(t)$ are given in Appendix~\ref{app:numdata}.
Consistent with our previous notation convention, we denote henceforth
\begin{equation}
\bD_j = \bD(j \tau), ~~ \ddot \bD_j = \bD''( j \tau).
\label{eq:notDj}
\end{equation}

Before we explain the second step in the flow chart, let us give a few technical details of the definition 
of the new wave $u^\ss(t,\bx)$ and the derivation of the inner product expression in {equation}~\eqref{eq:innD} 
of the data matrices. These details are not needed to compute $\bA^\RM$, which is why they are 
not in the flow chart, but they allow us to derive the expression of the mass and stiffness matrices 
in terms of the data. 

It is proved in \cite[~Appendix A]{borcea2020reduced} that 
\begin{align}
\hspace{-0.04in}\frac{[p^\ss(t,\bx) + p^\ss(-t,\bx)]}{c(\bx)/\bar c} = 
\cos \big(t \sqrt{\cA}\big) \widehat f \big(\sqrt{\cA} \big) \delta_{\bx_s}(\bx)\nonumber \\
=\sum_{j=1}^\infty \cos\big(t \sqrt{\la_j}\big) \widehat f \big(\sqrt{\la_j} \big) y_j(\bx_s) y_j(\bx),
\label{eq:TD1}
\end{align}
where 
\begin{equation}
\widehat f(\omega) = \int_{\RR} f(t) e^{i \omega t} dt
\end{equation} 
is the Fourier transform of the probing pulse and we define functions of the self-adjoint 
and positive definite operator $\cA$ using its spectral decomposition. 
If $\cA$ has the eigenvalues $\{\la_j\}_{j \ge 1}$ and the 
eigenfunctions $\{y_j\}_{j \ge1}$, then $\cos \big(t \sqrt{\cA}\big) $ is the operator 
with eigenvalues $\{ \cos \big(t \sqrt{\la_j} \big)\}_{j \ge 1}$ and the same eigenfunctions.
The operator $\widehat f \big(\sqrt{\cA} \big)$ is defined similarly. {The derivation of equation \eqref{eq:TD1} 
 involves the expansion of the wavefield in the basis 
$\{y_j(\bx)\}_{j \ge 1}$ of eigenfunctions of $\cA$ and manipulations of series.}
 
Next, we need the technical assumption that $\widehat f \ge 0$. This may not be the case in general, 
but the assumption can be achieved with simple processing as follows. Suppose that the probing pulse 
is actually some wavelet $\varphi(t)$ that is known or can be estimated \cite[]{pratt1999seismic}. 
Then, the measured wave convolved with $\varphi(-t)$ is the same as the solution of equation \eqref{eq:I1} 
evaluated at the receivers, with
\begin{equation}
f(t) = \varphi(t) \star_t \varphi(-t).
\end{equation}
Such $f(t)$ is obviously an even function, with Fourier transform $
\widehat f(\om) = |\widehat \varphi(\om)|^2 \ge 0,$ 
that is analytic by the Paley-Wiener-Schwartz theorem \cite[~Chapter VII]{Hormander}.

Analytic functions of $\cA$ commute, {as can be checked using power series}, so we can factor the right hand side in equation \eqref{eq:TD1} as 
\begin{equation}
\cos\big(t \sqrt{\cA}\big) \widehat f \big(\sqrt{\cA} \big) \delta_{\bx_s}(\bx) = 
\widehat f^\12\big(\sqrt{\cA} \big) u^\ss(t,\bx),
\label{eq:TD2}
\end{equation}
where 
\begin{equation}
u^\ss(t,\bx) = \cos \big( t \sqrt{\cA} \big) u_0^\ss(\bx),
\label{eq:cosA}
\end{equation}
is our new wave, with initial state 
\begin{equation}
u_0^\ss(\bx) = \widehat f^\12\big(\sqrt{\cA} \big)\delta_{\bx_s}(\bx).
\label{eq:u0}
\end{equation}
Note that $u^\ss(t,\bx)$ is just like the wave written in equation \eqref{eq:TD1}. 
The only difference is that it corresponds to a different pulse, with Fourier transform $\widehat f^\12$
instead of $\widehat f$. 

There are two important consequences of working with  $u^\ss(t,\bx)$. The first  is that 
by  the definition of $\cos\big( t \sqrt{\cA} \big)$, we can use the trigonometric identity 
\begin{equation}
\cos((t+\Delta t)\alpha) = 2 \cos(\Delta t \alpha) \cos(t \alpha) - \cos((t-\Delta t) \alpha), 
\end{equation}
for $\alpha = \sqrt{\la_j}$, with $j \ge 1$, to evolve the wave defined in equation \eqref{eq:cosA} over any  interval $\Delta t$,
\begin{align}
u^\ss(t+\Delta t,\bx) &= 2 \cos \big(\Delta t\sqrt{\cA} \big) u^\ss(t,\bx) - u^\ss(t-\Delta t,\bx).
\label{eq:tStep}
\end{align}
The second consequence is that 
the entries of $\bD(t)$, defined in equation \eqref{eq:defD}, admit a useful symmetric inner product expression
\begin{align}
D^{(r,s)}(t) &= p^\ss(t,\bx_r) + p^\ss(-t,\bx_r) \nonumber \\
&= \int_{\Omega} d \bx \, \delta_{\bx_r} (\bx) \widehat f^\12\big(\sqrt{\cA} \big) u^\ss(t,\bx)
\nonumber \\
&= \int_{\Omega} d \bx \, \big[\widehat f^\12\big(\sqrt{\cA} \big) \delta_{\bx_r}(\bx)\big]  u^\ss(t,\bx) \nonumber \\
&= \int_{\Omega} d \bx \, u_0^{(r)}(\bx) u^\ss(t,\bx) \nonumber \\
&=  \int_{\Omega} d \bx \, u_0^{(r)}(\bx)  \cos \big(t \sqrt{\cA}\big) u^\ss_0(\bx), \label{eq:innD}
\end{align}
for $1 \le r,s \le \m$. 
The second equality in this equation is from equation \eqref{eq:TD1} and the assumption $c(\bx_r) = \bar c$, 
the third equality is because $\cA$ and therefore $\widehat f^\12\big(\sqrt{\cA} \big)$ are self-adjoint 
operators that commute and the last equalities follow from {equations}~\eqref{eq:cosA} and \eqref{eq:u0}.
We also have 
\begin{align}
&\frac{d^2 D^{(r,s)}(t)}{dt^2} =  \int_{\Omega} d \bx \, u_0^{(r)}(\bx) \partial_t^2 u^\ss(t,\bx) \nonumber \\
& \hspace{0.2in}= - \int_{\Omega} d \bx \, u_0^{(r)}(\bx) \cA  u^\ss(t,\bx), \quad 1 \le r,s \le \m.
\label{eq:inddotD}
\end{align}

Now we can describe how we use equations \eqref{eq:tStep}--\eqref{eq:inddotD} to complete the second step in 
the flow chart of Figure~\ref{Fig.FC}. With the notation 
\begin{equation}
\langle \boldsymbol{\phi},\boldsymbol{\psi} \rangle = \int_\Omega d \bx \,
\boldsymbol{\phi}^T(\bx) \boldsymbol{\psi}(\bx)
\end{equation}
for the integral of the outer product of any two functions $\boldsymbol{\phi}(\bx)$ and 
$\boldsymbol{\psi}(\bx)$ with values in $\RR^{1\times \m}$, and from the definition  in equation 
\eqref{eq:defM}, we compute the $\m \times \m$ blocks of the mass matrix as
\begin{align}
\bM_{i,j} &= \langle \bu_i, \bu_j \rangle
= \langle \cos  \big(i \tau \sqrt{{\cal A}} \big) \bu_0, 
\cos \big(j \tau \sqrt{{\cal A}} \big) \bu_0 \rangle \nonumber \\
&= \langle \bu_0, \cos  \big(i \tau \sqrt{{\cal A}} \big) 
\cos \big(j \tau \sqrt{{\cal A}} \big) \bu_0 \rangle \nonumber \\
&= \frac{1}{2} \langle \bu_0, \big[ \cos  \big((i+j) \tau \sqrt{{\cal A}} \big) + 
\cos \big(|i-j| \tau \sqrt{{\cal A}} \big) \big] \bu_0 \rangle \nonumber \\
&= \frac{1}{2} \left( \bD_{i+j} + \bD_{|i-j|} \right), \quad 0 \le i,j \le \n-1. \label{eq:M32}
\end{align}
The second line in this equation is because $\cA$ and therefore $\cos\big(i \tau \sqrt{{\cal A}} \big)$ 
are self-adjoint operators that commute, the third line is due to equation \eqref{eq:tStep}, 
evaluated at $t = i \tau$ and $\Delta t = j \tau$, and the last line is by equation \eqref{eq:innD}.
The blocks of the stiffness matrix defined in equation \eqref{eq:defS} are 
\begin{align}
 \bS_{i,j} &= \langle \bu_i, {\cal A} \bu_j \rangle 
= \langle \cos  \big(i \tau \sqrt{{\cal A}} \big) \bu_0, 
{\cal A} \cos  \big(j \tau \sqrt{{\cal A}} \big) \bu_0 \rangle \nonumber \\
&= \langle \bu_0, {\cal A} \cos \big(i \tau \sqrt{{\cal A}} \big) 
\cos \big(j \tau \sqrt{{\cal A}} \big) \bu_0 \rangle \nonumber \\
&= \frac{1}{2}\langle \bu_0, {\cal A} \bu_{i+j} + {\cal A} \bu_{|i-j|}\rangle \nonumber \\
&= -\frac{1}{2} \left(\ddot \bD_{i+j} + \ddot \bD_{|i-j|} \right), \quad 0 \le i,j \le \n-1,
\label{eq:M43}
\end{align}
where we used again the self-adjointness of $\cA$, and equation \eqref{eq:tStep} evaluated at 
$t = i \tau$ and $\Delta t = j \tau$. The last equality is by equation \eqref{eq:inddotD}.
The block structure of the  matrices $\bM$ and $\bS$ is sketched in Figure~\ref{Fig.FC} for the
case $\n=5$.

The remaining two steps in the flow chart in Figure~\ref{Fig.FC} are self-explanatory and have 
been motivated in the previous subsection. We summarize the computation of $\bA^\RM$
in the following algorithm.

\begin{algorithm}\textbf{\emph{(Data-driven ROM operator)}}
\label{alg:arom}
\vspace{0.04in}

\vspace{0.04in} \noindent \textbf{Input:} 
The matrix $\boldsymbol{\cal{M}}(t)$ of measurements given by equation \eqref{eq:I3}, 
at time {instances} $t = j \tau$, for 
{$j = -N_{\rm f}, \ldots, 2\n-2$, with $N_{\rm f} = [t_{\rm f}/\tau]$}. 
{We have $\boldsymbol{\cal{M}}(j \tau) =0$ for $j< -N_{\rm f}$.}

\vspace{0.04in} \noindent 1. Compute 
\begin{equation*}
\bD_j = \boldsymbol{\cal{M}}(j \tau) +\boldsymbol{\cal{M}}(-j \tau), \quad 0 \le j \le 2\n-2.
\end{equation*}

\vspace{0.04in} \noindent 2. Compute $\{\ddot \bD_j\}_{j=0}^{2\n-2}$  using, 
e.g., the Fourier transform (see Appendix \ref{app:numdata}).

\vspace{0.04in} \noindent 3. Calculate
$\bM,  \bS \in \mathbb{R}^{\m\n \times \m\n}$ with the block entries
\begin{align*}
\bM_{i,j}  & =  \frac{1}{2} \big( \bD_{i+j} + \bD_{|i-j|} \big)
\in \mathbb{R}^{\m \times \m}, \\
 \bS_{i,j} & =  - \frac{1}{2} \big( \ddot{\bD}_{i+j} + \ddot{\bD}_{|i-j|} \big)
\in \mathbb{R}^{\m \times \m},
\end{align*}
for $0 \le i,j  \le \n-1$.

\vspace{0.04in} \noindent 4. Perform the block Cholesky factorization $ \bM = \bR^T \bR $ using 
\cite[~Algorithm 5.2]{druskin2018nonlinear}.

\vspace{0.04in} 
\noindent \textbf{Output:}  $\bA^\RM = \bR^{-T} \bS \bR^{-1}$.
\end{algorithm}

\subsection{ROM based velocity estimation}

We estimate $c(\bx)$ by minimizing the misfit of the ROM, as in equation \eqref{eq:ROM_obj}. 
The computation of the term ${\cal R} \big(\cF[v]\big)$ in that equation 
involves two steps. The first step is to solve the wave equation \eqref{eq:I1} with $c(\bx)$ replaced 
by the search velocity $v(\bx)$. The solution evaluated at the receivers gives $\cF[v](t)$. 
The second step is to apply Algorithm \ref{alg:arom} with input $\cF[v](t)$.
In an abuse of notation, we let henceforth 
\begin{equation}
\bA^\RM(v) = {\cal R} \big(\cF[v]\big).
\end{equation}

The search space  $\mathcal{C}$,  where $v(\bx)$ lies, is parametrized using some appropriate 
basis functions $\{\phi_l(\bx)\}_{l=1}^N$
\begin{equation}
{v}(\bx;\bet) = c_o(\bx) + \sum_{l = 1}^N \eta_l \phi_l(\bx),
\label{eq:IM2}
\end{equation}
where $c_o(\bx)$ is the initial guess.
The optimization is then $N$-dimensional, for the vector $\bet = (\eta_1, \ldots, \eta_N)^T$ 
of coefficients in  {equation}~\eqref{eq:IM2}. 

%\hyphenation{com-pu-ta-tion}

The causality of the ROM (Appendix~\ref{app:Causal}) allows us to carry out the inversion in a layer 
stripping fashion, from the data at time {instances} $\{t_j = j \tau\}_{j=0}^{2 \k -2}$, with $ \k \le \n$. 
To do so, we replace $\bA^\RM(v)$ and $\bA^\RM$ in the objective function by the upper left $\k\m \times \k\m$ 
blocks of these matrices, denoted by $\big[\bA^\RM(v)\big]_\k$ and $\big[\bA^\RM\big]_\k$, respectively.

Since $\bA^\RM$ and thus $\big[\bA^\RM\big]_\k$ are symmetric matrices, it is enough to consider their 
block upper triangular part in the optimization. As shown in Appendix~\ref{app:algrom}, the entries of 
$\bA^\RM$ decay away from the diagonal. Thus, we can ease the computational burden by including 
only the first few $d \m$ diagonals in the objective function, where $d$ is an integer between $1$ and $\k$. 
For this purpose, we denote by 
\begin{equation}
{\rm Rest}_{d,\k}:\RR^{\k\m \times \k\m} \mapsto \RR^{d \m (\k\m- (d \m-1)/2)}
\label{eq:restdk}
\end{equation}
the mapping that takes a $\k\m \times \k\m $ matrix, keeps only its  first $d \m$ upper diagonals, 
including the main one, and puts their entries into a column vector, of length
\begin{equation}
\sum_{j=0}^{ d \m-1} (\k\m -j) = d \m [\k\m - ( d \m-1)/2].
\end{equation}
The objective function that takes into account both the time windowing and the restriction of the ROM 
to a few diagonals is denoted henceforth by 
\begin{equation}
\mathcal{O}_{d,\k}({v}) = \left \|\mbox{Rest}_{d,\k} \big(\left[\bA^\RM(v)-\bA^\RM\right]_\k\big)\right\|_2^2,
\label{eq:newObj}
\end{equation}
where $\| \cdot \|_2$ is the vector Euclidean norm.

\begin{algorithm}\textbf{\emph{({ROM based} velocity estimation)}}
\label{alg:prowi}

\vspace{0.04in} \noindent  \textbf{Input:} The data driven $\bA^\RM$. 

\vspace{0.04in} \noindent  1. Set the number of layers for the layer stripping approach to 
$\ell$ and the number of iterations per layer to $\q$.

\vspace{0.04in} \noindent  2. Choose $\ell$ natural numbers $\{k_l\}_{l=1}^\ell$, satisfying
\begin{equation*}1 \le k_1 \leq k_2 \leq  \cdots \leq  k_{\ell} = \n.\end{equation*} The data subset for the $l^{\rm th}$ layer is 
$\{\bD_j, \ddot \bD_j\}_{j = 0}^{2k_l-2}$.

\vspace{0.04in} \noindent  3.  Starting with the initial vector $\bet^{(0)}= {\bf 0}$, proceed: 
\item For $l = 1,2,\ldots, \ell$, and $j = 1,\ldots,\q$, set the update index $i = (l-1)\q+j$. 
Compute $\bet^{(i)}$ as a Gauss-Newton update for minimizing the functional
\begin{align}
\hspace{-0.0in}\mathcal{L}_i(\bet) & = \mathcal{O}_{d,k_l} \big(v(\cdot; \bet)) + \mathcal{L}_i^{\rm reg}(\bet),
\label{eqn:regObj}
\end{align}
linearized about $\bet^{(i-1)}$. The term $\mathcal{L}_i^{\rm reg}(\bet)$ introduces a user defined 
regularization penalty in the optimization.

\vspace{0.04in} \noindent  \textbf{Output:} The velocity estimate $c^{\rm est}(\bx) = {v}(\bx;\bet^{(\ell \q)})$.
\end{algorithm}

The details on our implementation of Algorithm \ref{alg:prowi} and the regularization penalty 
are provided in Appendix \ref{app:implem}.

\subsection{Computational cost} %{TO BE EDITED}}

Since our Algorithm 2 for ROM based velocity estimation uses a Gauss-Newton iteration to minimize the 
objective function {in equation} \eqref{eqn:regObj}, we compare its cost to that of the Gauss-Newton 
method for minimizing the FWI objective function  in {equation}~\eqref{eq:FWI_obj}. 
The same parametrization of the search velocity is assumed for both approaches.

{The numerical examples considered below are for two-dimensional media $\Omega \subset \mathbb{R}^2$
with a relatively modest number  $\m$ of {colocated sources/receivers}, not exceeding $60$. In such settings the} cost of 
each Gauss-Newton step is dominated by the computation of the Jacobian of the objective function. 
This computation requires solving the forward problem for all $\m$ sources. The ROM based 
approach requires, in addition, the computation of $\bA^\RM$ {and its derivatives}. 
We compare next the cost of solving the forward problem with that of computing  the ROM with 
Algorithm \ref{alg:arom}.

We solve the forward problem (equations \eqref{eq:I1}--\eqref{eq:I2}) in a rectangular domain $\Omega$, 
with homogeneous Dirichlet boundary conditions at $\partial \Omega$, using explicit time stepping, a 
three point finite difference approximation of $\partial_t^2$ with step $\tau_\cf$, and a five point finite 
difference discretization of the Laplacian on a uniform mesh with $N_\cf$ points. 
To write down the order of $N_\cf$, let $\bar \lambda$ be the reference wavelength,  
calculated with the constant reference speed $\bar c$ and at the central frequency of the probing signal $f(t)$. 
An accurate and stable forward solver requires a mesh size $h$ that is a small fraction of the wavelength 
and does not exceed $\bar c \tau_\cf$. The number of mesh points is therefore 
\begin{equation}
N_\cf = \frac{\mbox{area}(\Omega)}{h^2} \gg \m \n,
\end{equation}
where the inequality is because the {colocated sources/receivers} are at $O(\bar \lambda)$ distance, the array length is 
$O(m \bar \lambda)$ which is usually much smaller than the width of $\Omega$, 
and the time sample $\tau$ used in the ROM construction is much larger than $\tau_\cf$. 
Each time step requires multiplying an $N_\cf \times N_\cf$ sparse matrix with a vector in 
$\RR^{N_\cf}$, at an $O(N_\cf)$ cost. Thus, the cost of solving the forward problem, 
for the $\m$ sources and up to time $T$, is 
\begin{equation}
\mbox{cost}(\cF) = O(\m n_{\cf} N_{\cf}),
\end{equation} 
where
$n_\cf = T/\tau_\cf \gg \n$. Recall that $\cF$ denotes the forward map.

The computational  cost of running Algorithm~\ref{alg:arom} lies mainly in the block Cholesky 
factorization (see {equation}~\eqref{eq:Cholesky}) and the operator ROM computation from equation \eqref{eq:ROMA}, 
where $\bR^{-1}$ can be calculated by block-wise backward substitution. 
Therefore, the cost of computing  $\bA^\RM$ is estimated at 
\begin{equation}
\mbox{cost}(\bA^\RM) = O(\m^3 \n^3),
\end{equation}
and it is typically smaller than $\mbox{cost}(\cF)$ if the array is not too large and {we  sample in time  at about the Nyquist rate, 
as explained below, after equation \ref{eq:pulse}}. The bulk of the computational 
cost of derivatives of $\bA^\RM$ is in the differentiation of the block Cholesky factors $\bR$. 
This cost is essentially the same as that of the block Cholesky factorization itself, since the 
derivatives of $\bR$ can be computed by a similar factorization algorithm, as described in detail in 
\cite[~Appendix A]{borcea2014model}.}

\label{sect:sum}

{For three-dimensional media $\Omega \subset \mathbb{R}^3$ and settings with large $\m$,
the dominant computational cost is not in {the} Jacobian calculation itself, but in solving the regularized
normal equations for the Gauss-Newton update direction for the objective function in equation \ref{eqn:regObj}. 
While small-scale examples allow for direct computation of the update direction using, e.g., 
equation \ref{eqn:updatedir}, large-scale settings call for iterative approaches like the Conjugate Gradient
method. Note, however, that in such settings the computational cost difference between the conventional 
FWI and ROM based velocity estimation virtually disappears, since the sizes of the Jacobians of both methods
can be made essentially identical by an appropriate choice of parameter $d$ in equation \ref{eq:newObj}.}

\section{Numerical illustration}

In this section we give two numerical illustrations of the benefits of the velocity estimation with the ROM operator 
vs. FWI. We assume, as in the theory section above, knowledge of the noiseless array response matrix 
$\boldsymbol{\cal M}(t)$. Noisy measurements  and the approximation of $\boldsymbol{\cal M}(t)$ 
from towed-streamer data are considered in the next section.

The first illustration is for  a two-parameter velocity model, where we can plot the objective function 
over the search space. The second is for  the ``Camembert example" introduced in \cite[]{gauthier1986two} 
to demonstrate the challenge of velocity estimation with FWI. We also display components of $\bU(\bx)$ and 
$\bV(\bx)$ for the Camembert example, to illustrate the properties of the projection basis discussed in the theory section.

All the numerical results are for the source pulse
\begin{equation}
f(t) = \cos(\om_o t) \exp \Big[-\frac{(2 \pi B)^2 t^2}{2}\Big],
\label{eq:pulse}
\end{equation}
with central frequency $\om_o/(2 \pi) = 6$~Hz and bandwidth $B=4$~Hz.  
See Appendix~\ref{app:numdata} for details on the numerically simulated data. 
To choose $\tau$, we use $\om_o/(2 \pi) + B = 10$~Hz as the Nyquist frequency. 
Thus, for $\tau = 1/(2.3 \cdot 10~\mbox{Hz}) = 0.0435~\mbox{s}$, 
the data are sampled at $2.3$ points per wavelength.

The array of $\m$ sensors is at $150$~m below the top boundary. The sensor spacing is $160.3$~m for 
the two-parameter velocity model and $155.5$~m for the Camembert example.
For each simulation we specify $\m$, the size of the rectangular domain $\Omega$, 
the data sampling interval $\tau$ and the number $\n$ of snapshots that define the approximation space. 

\subsection{Topography of the objective function}
\begin{figure*}[ht]
\begin{center}
\begin{tabular}{ccc}
(a) Velocity (m/s) &
(b) Log of FWI misfit &
(c) Log of  ROM misfit \\
\hspace{-0.01\textwidth}\includegraphics[width=0.32\textwidth]{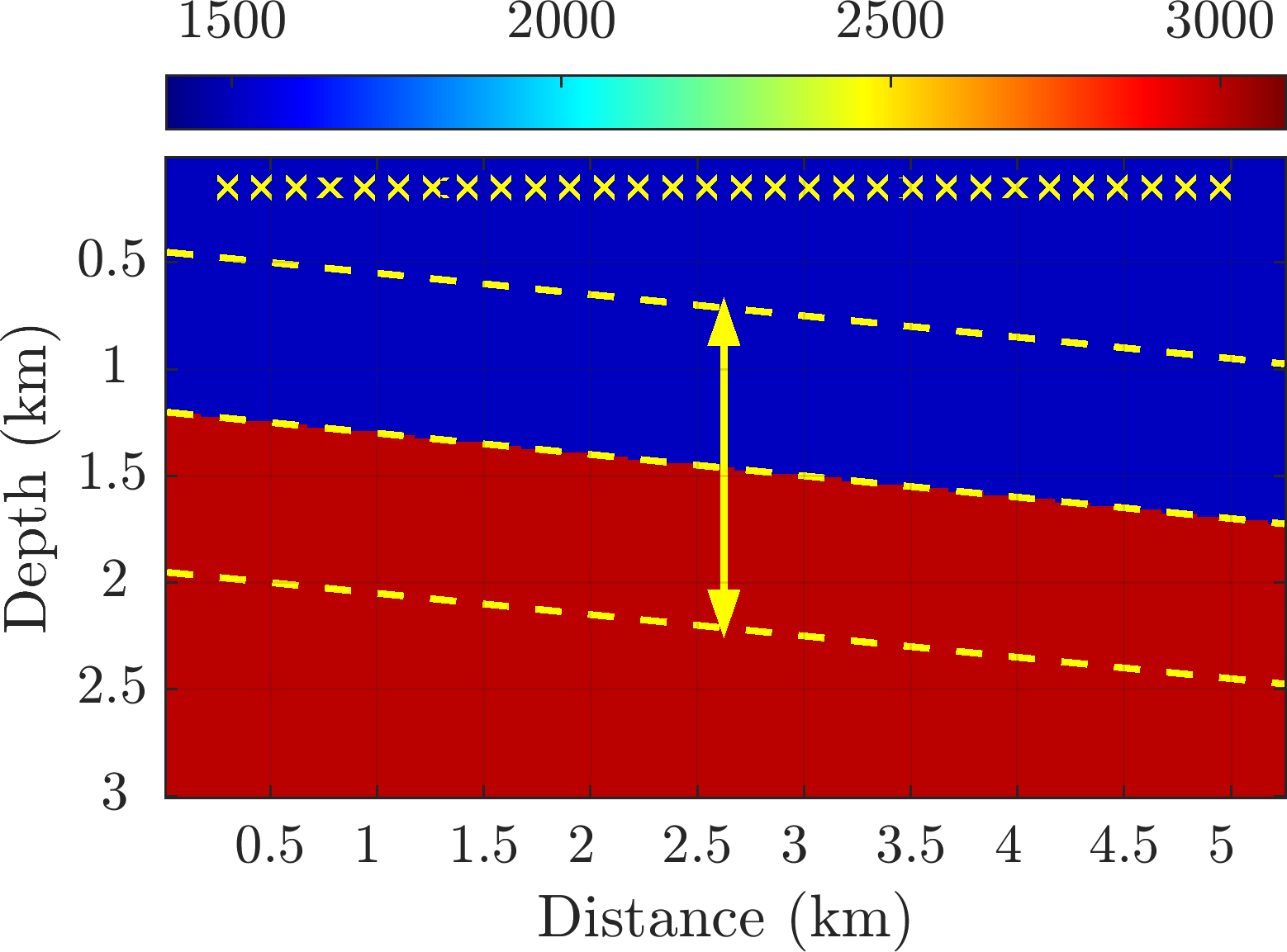} &
\hspace{-0.015\textwidth}\includegraphics[width=0.331\textwidth]{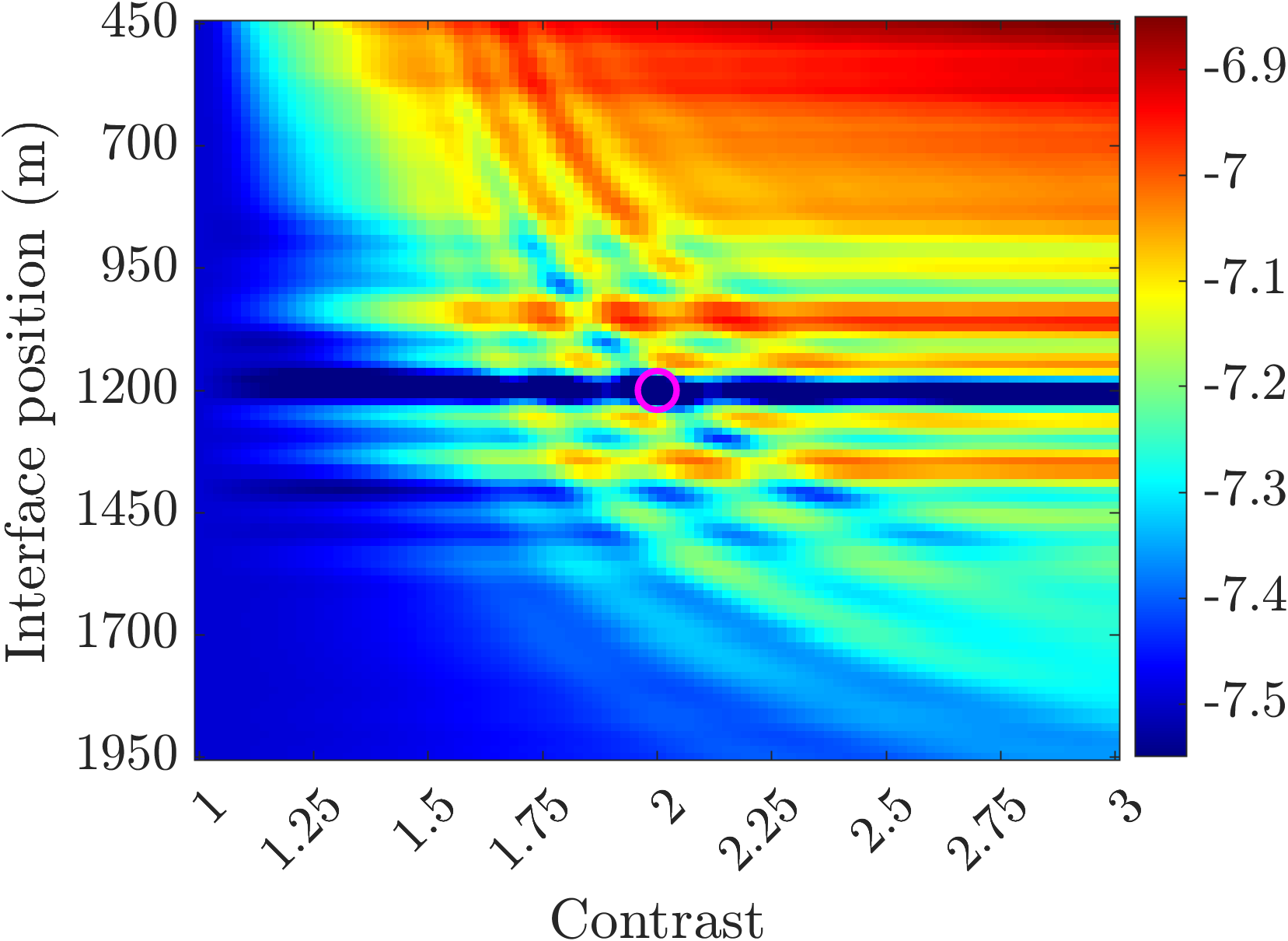} &
\hspace{-0.015\textwidth}\includegraphics[width=0.331\textwidth]{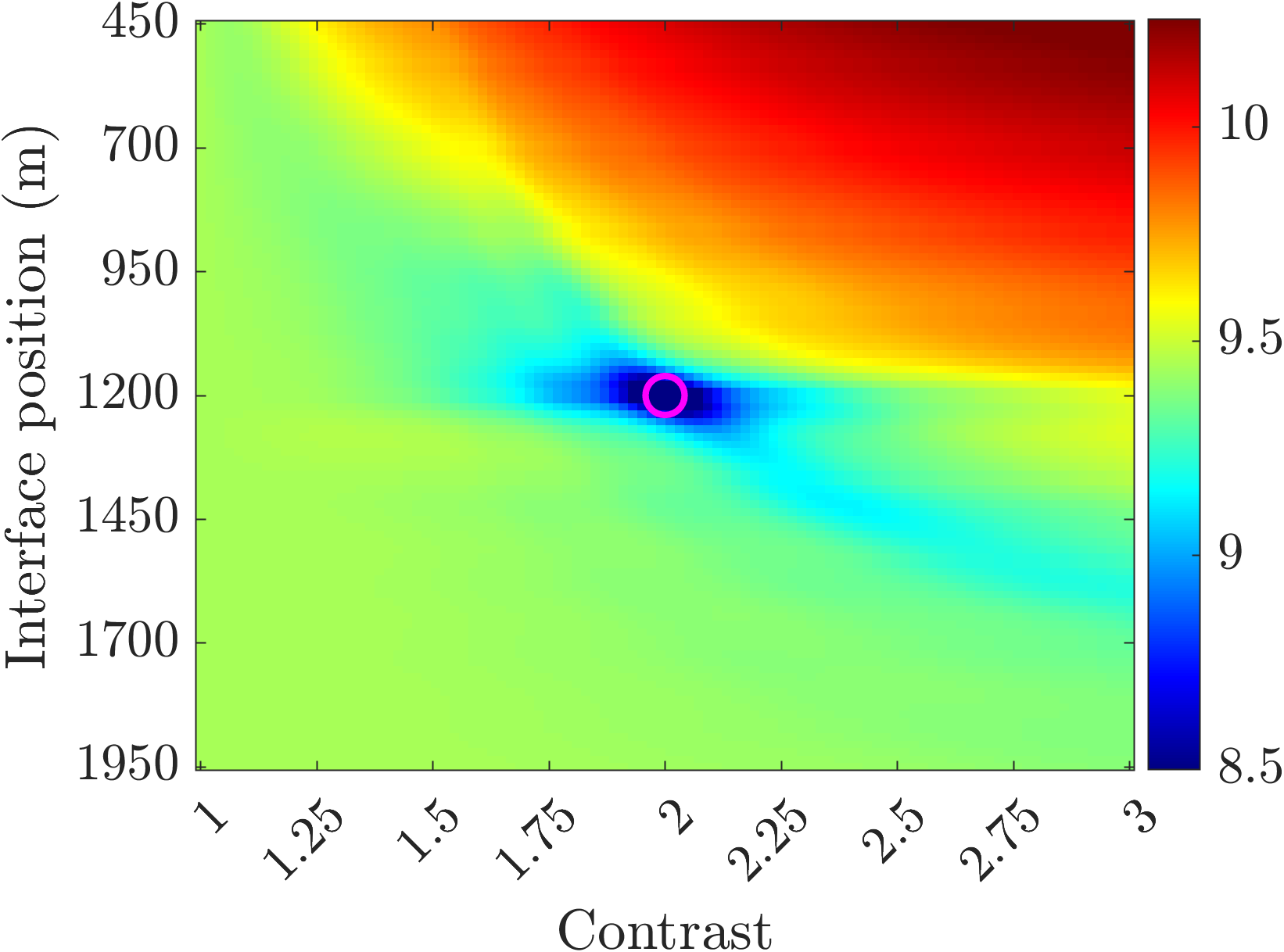}
\end{tabular}
\end{center}
\vspace{-0.2in}
\caption{Objective functions topography study:
(a) Velocity model used in objective topography study.
The middle dashed line shows the actual interface location, while the top and bottom 
dashed lines show the extent of the interface location parameter sweep. 
All $\m=30$ {colocated sources/receivers} are shown as yellow $\times$. Velocity colorbar is in $\rm{m/s}$;
(b)--(c) Decimal logarithms of the objective functions \eqref{eq:FWIObj}--\eqref{eq:ROMObj}, 
vs. the interface position and velocity contrast. The actual position and contrast parameters are 
indicated by $\textcolor{magenta}{\bigcirc}$. {These true values are not included in the search space.}}
\label{fig:topo}
\end{figure*}

Consider the velocity model displayed in Figure~\ref{fig:topo}a, in the domain 
$\Omega = [0, 5\text{ km}] \times [0, 3\text{ km}]$. It consists of two homogeneous regions separated 
by a slanted interface. The top region has the slower velocity $c_{\rm t} = 1500 \text{ m/s}$, while the bottom
region has the faster velocity $c_{\rm b} = 3000 \text{ m/s}$. The purpose of this example is to visualize 
the objective function, so we do not run Algorithm \ref{alg:prowi} and we do not use a search velocity 
of the form given in equation \eqref{eq:IM2}. Instead, we sweep a two-parameter search space: The first parameter is the 
interface position in the search interval $[0.47\text{ km}, 1.95\text{ km}]$, measured as the depth of the 
leftmost point of the interface. The actual position is $1.2\text{ km}$. The second parameter is the 
contrast $c_{\rm b}/c_{\rm t}$ in the interval $[1,3]$. The actual contrast is {two}. The angle of the interface is kept 
constant and equal to the actual angle.

In Figures~\ref{fig:topo}b--\ref{fig:topo}c we display the decimal logarithms of two objective functions, calculated 
for $\m = 30$ {colocated sources/receivers} and   $\n = 39$  time samples at interval  $\tau = 0.0435~\mbox{s}$. 
The first objective function is for the FWI approach, 
\begin{equation}
\mathcal{O}^\FWI({v}) = \sum_{k=0}^{2\n-1} \left\| {\rm Triu} \big(\bD_\k({v})-\bD_\k \big)\right\|_2^2,
\label{eq:FWIObj}
\end{equation}
where $\bD_\k({v})$ are the $\m\times \m$ data matrices for the search velocity  $v(\bx)$ and
${\rm Triu}:\RR^{\m \times \m} \mapsto \RR^{\m(\m+1)/2}$ is the mapping that takes a symmetric 
$\m\times \m$ matrix, extracts its upper triangular part, including the main diagonal,
and arranges its entries into a $\m(\m+1)/2$-dimensional column vector.
The second objective function measures the misfit of the ROM 
\begin{equation}
\mathcal{O}^\RM({v}) = \left\|{\rm Triu} \big(\bA^\RM(v)-\bA^\RM\big)\right\|_2^2.
\label{eq:ROMObj}
\end{equation}
This corresponds to the particular case $d = \k = \n$ of the objective function in equation \eqref{eq:newObj}.

We observe in Figure \ref{fig:topo}b that the FWI objective function displays numerous local minima, 
at points in the search space that are far from the true one, marked in the plots by the magenta circle. 
There is no minimum at this circle because the exact values of the interface position and contrast are not 
in our parameter grid search space. The clearly visible horizontal stripes in Figure \ref{fig:topo}b 
are manifestations of cycle skipping. The ROM operator misfit shown in Figure \ref{fig:topo}c 
is smooth and has a single minimum, at the true interface position and contrast.

\subsection{The ``Camembert" example}

\begin{figure}[h!]
\begin{center}
\begin{tabular}{c}
\includegraphics[width=0.3\textwidth]{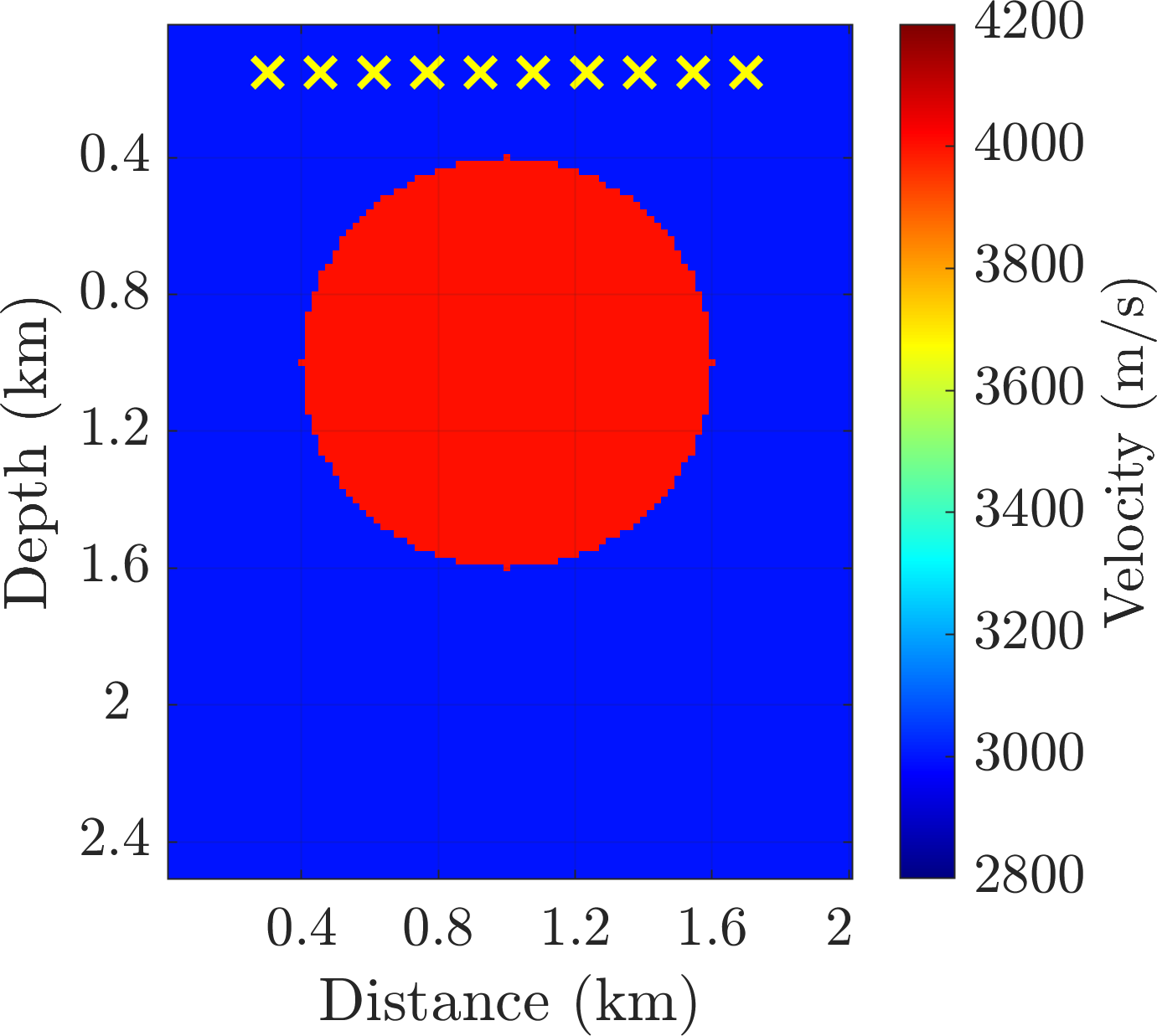}
\end{tabular}
\end{center}
\vspace{-0.2in}
\caption{Camembert velocity model. All $\m=10$ {colocated sources/receivers} are shown as yellow $\times$. 
Velocity colorbar is in $\rm{m/s}$.}
\label{fig:Camembert}
\end{figure}

\begin{figure*}[ht]
\begin{center}
\begin{tabular}{l@{\hspace{0.55em}} l @{\hspace{0.5em}}  l @{\hspace{0.55em}}  l @{\hspace{0.55em}}  }
(a) {\footnotesize ROM estimate iter. $10$} & (b) {\footnotesize ROM estimate iter. $20$} & 
(c) {\footnotesize ROM estimate iter. $40$} & (d) {\footnotesize ROM estimate iter. $60$} \\ 
\hspace{-0.01\textwidth}\includegraphics[width=0.21\textwidth]
{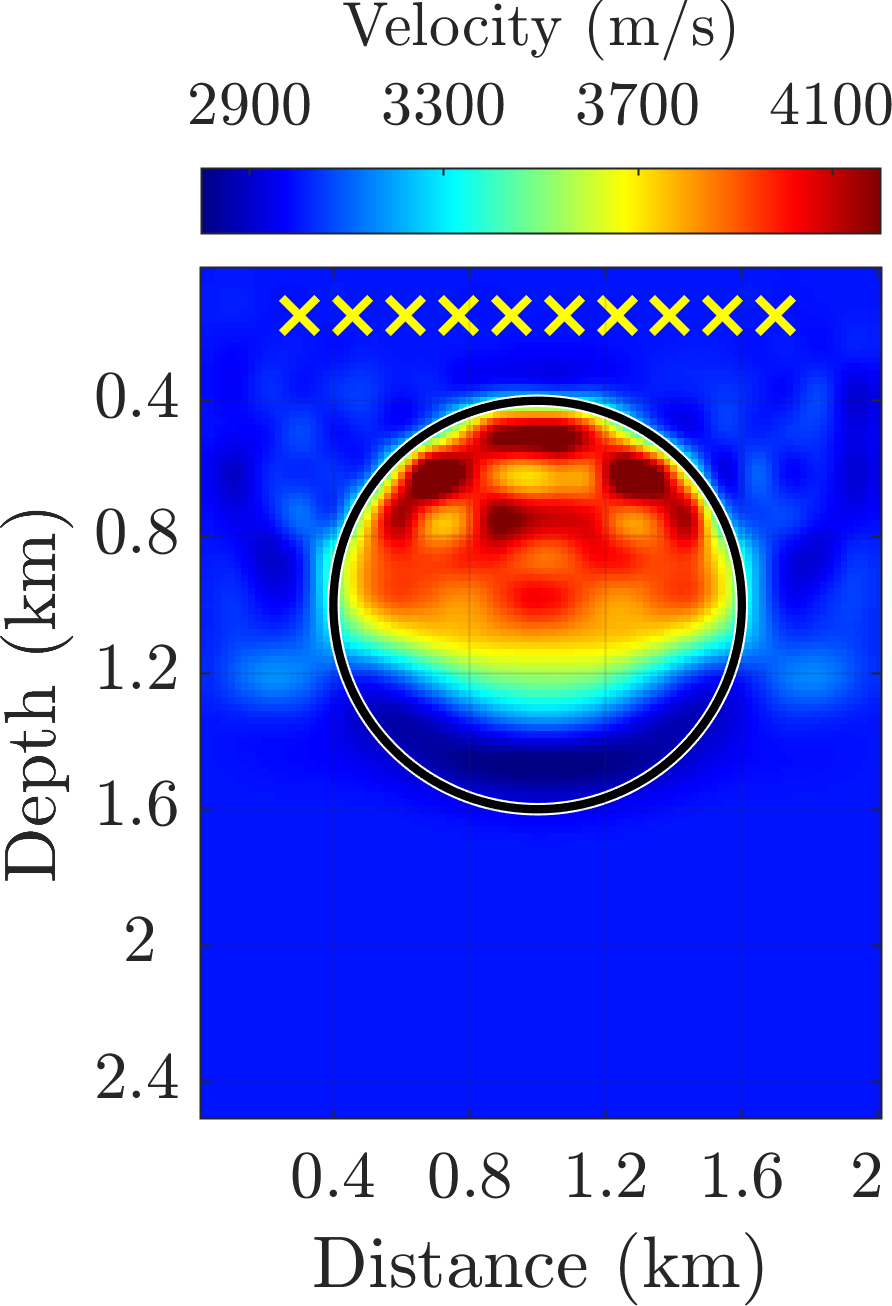} &
\hspace{-0.01\textwidth}\includegraphics[width=0.21\textwidth] 
{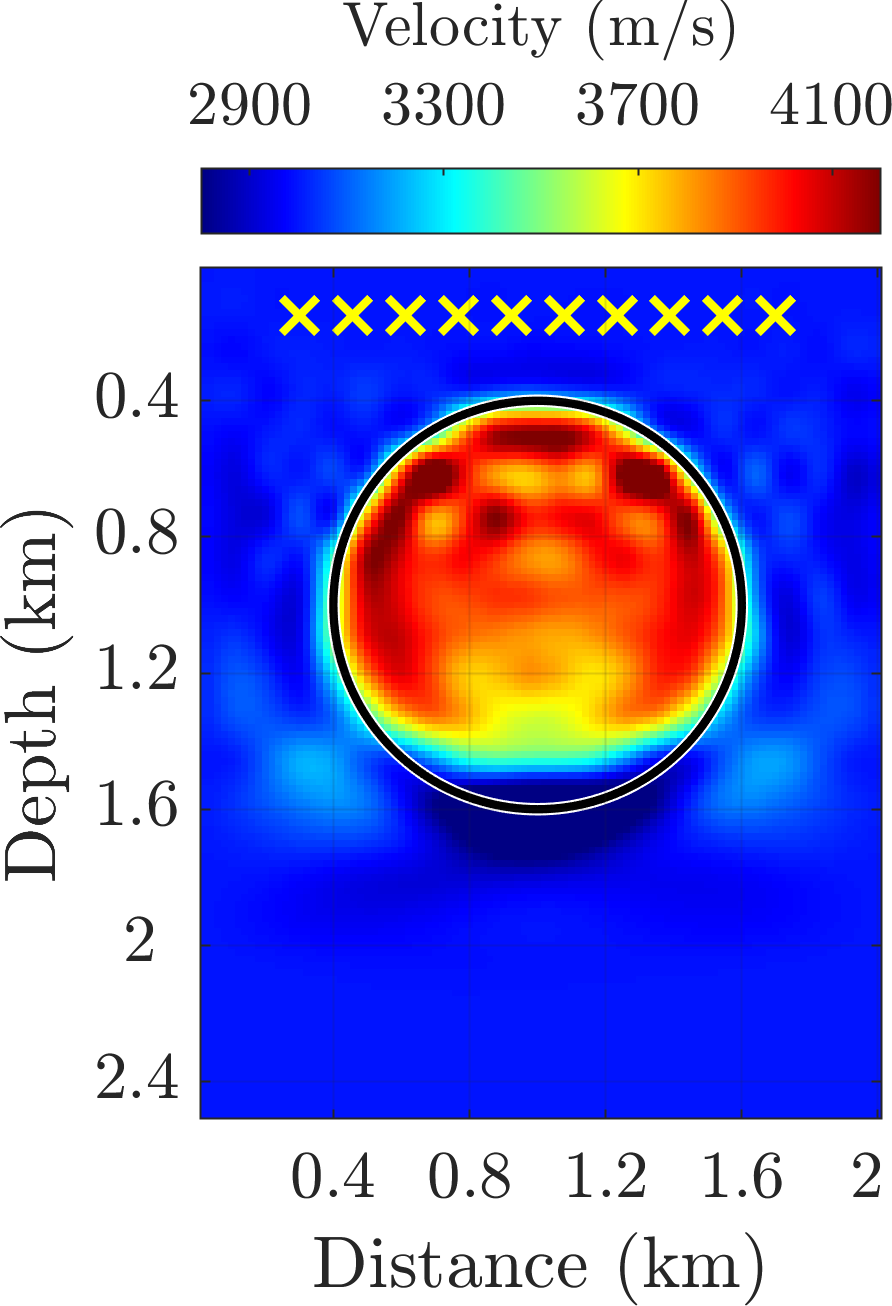} &
\hspace{-0.01\textwidth}\includegraphics[width=0.21\textwidth]
{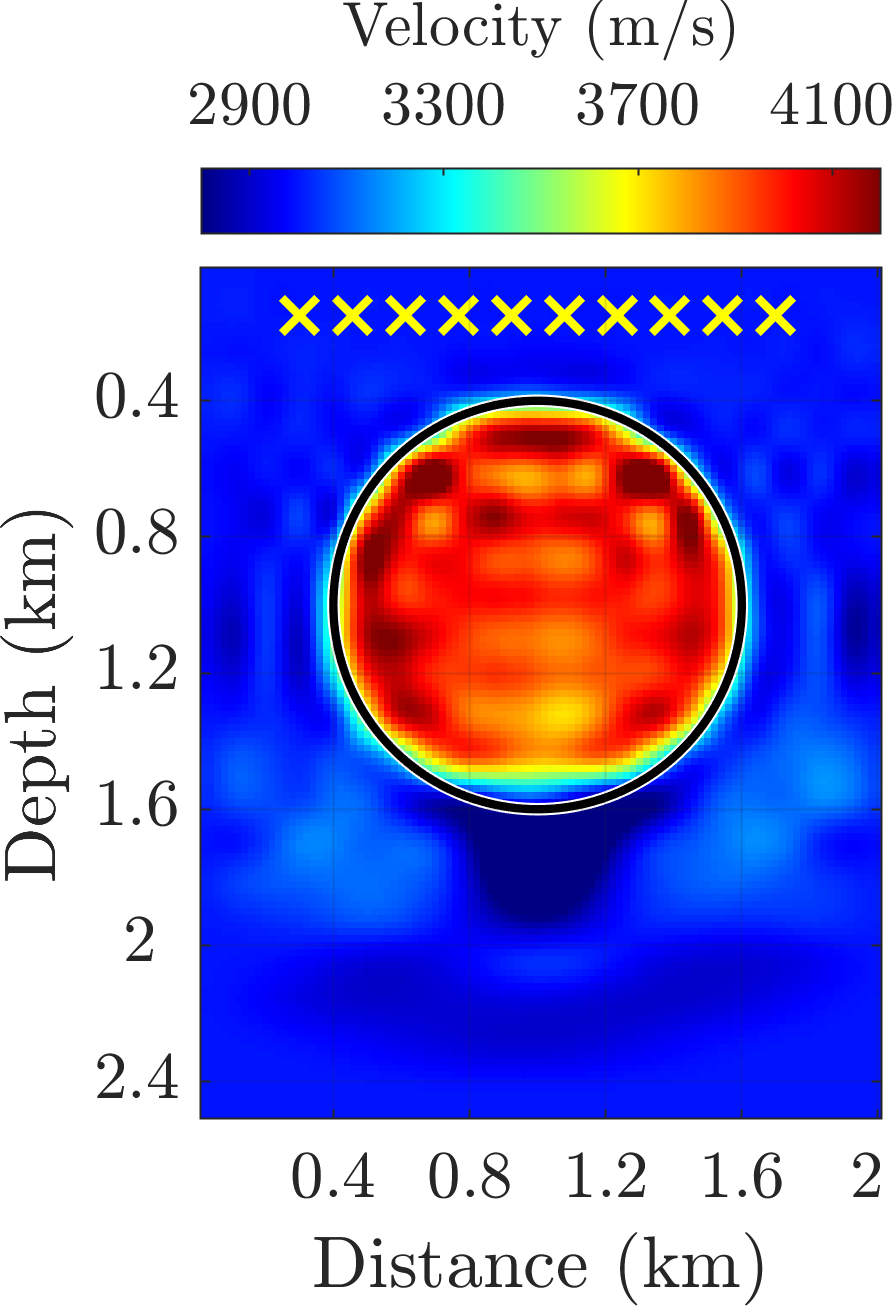} &
\hspace{-0.01\textwidth}\includegraphics[width=0.21\textwidth]
{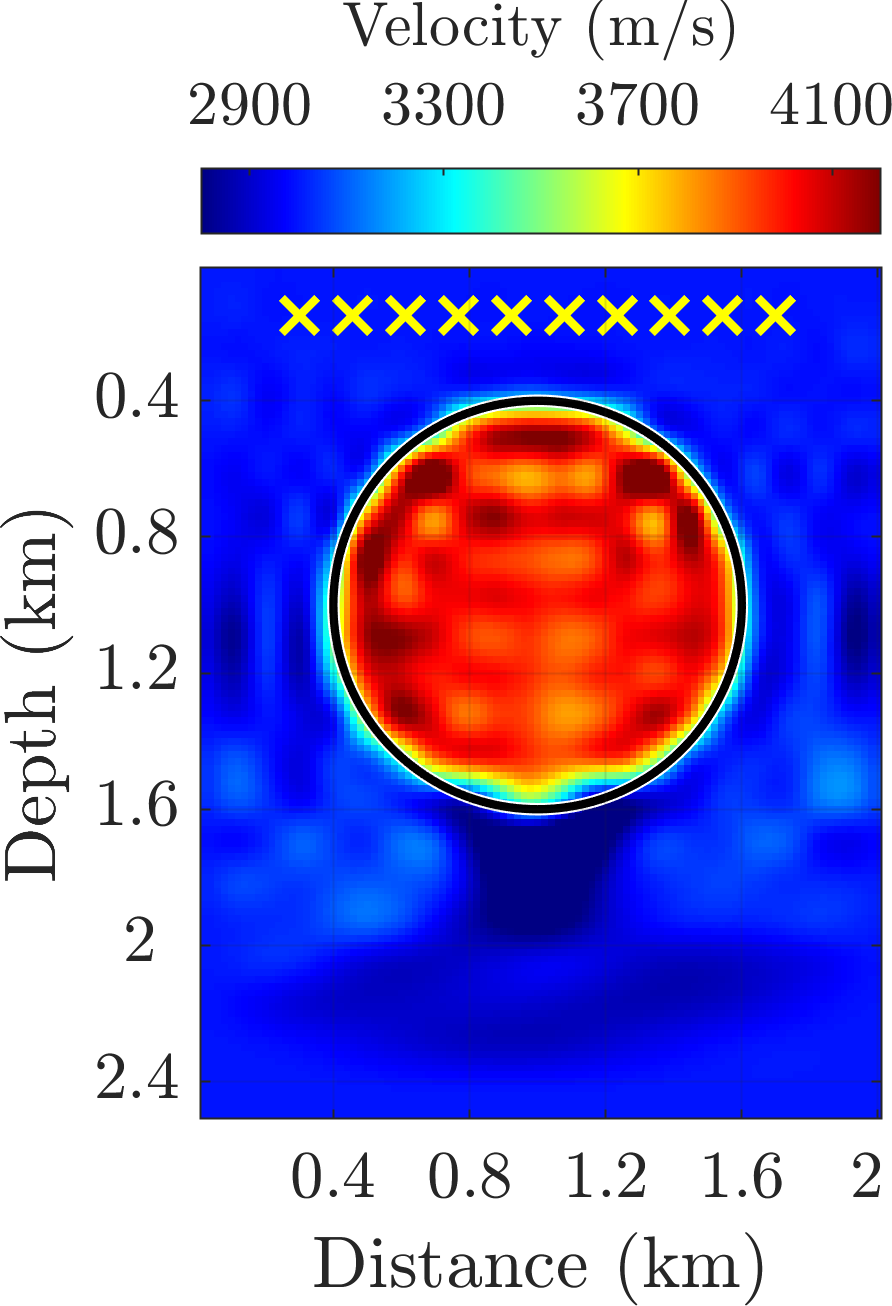} \\
(e) {\footnotesize FWI estimate iter. $10$} & (f) {\footnotesize FWI estimate iter. $20$} & 
(g) {\footnotesize FWI estimate iter. $40$} & (h) {\footnotesize FWI estimate iter. $60$} \\
\vspace{-0.2in}
\hspace{-0.01\textwidth}\includegraphics[width=0.21\textwidth]
{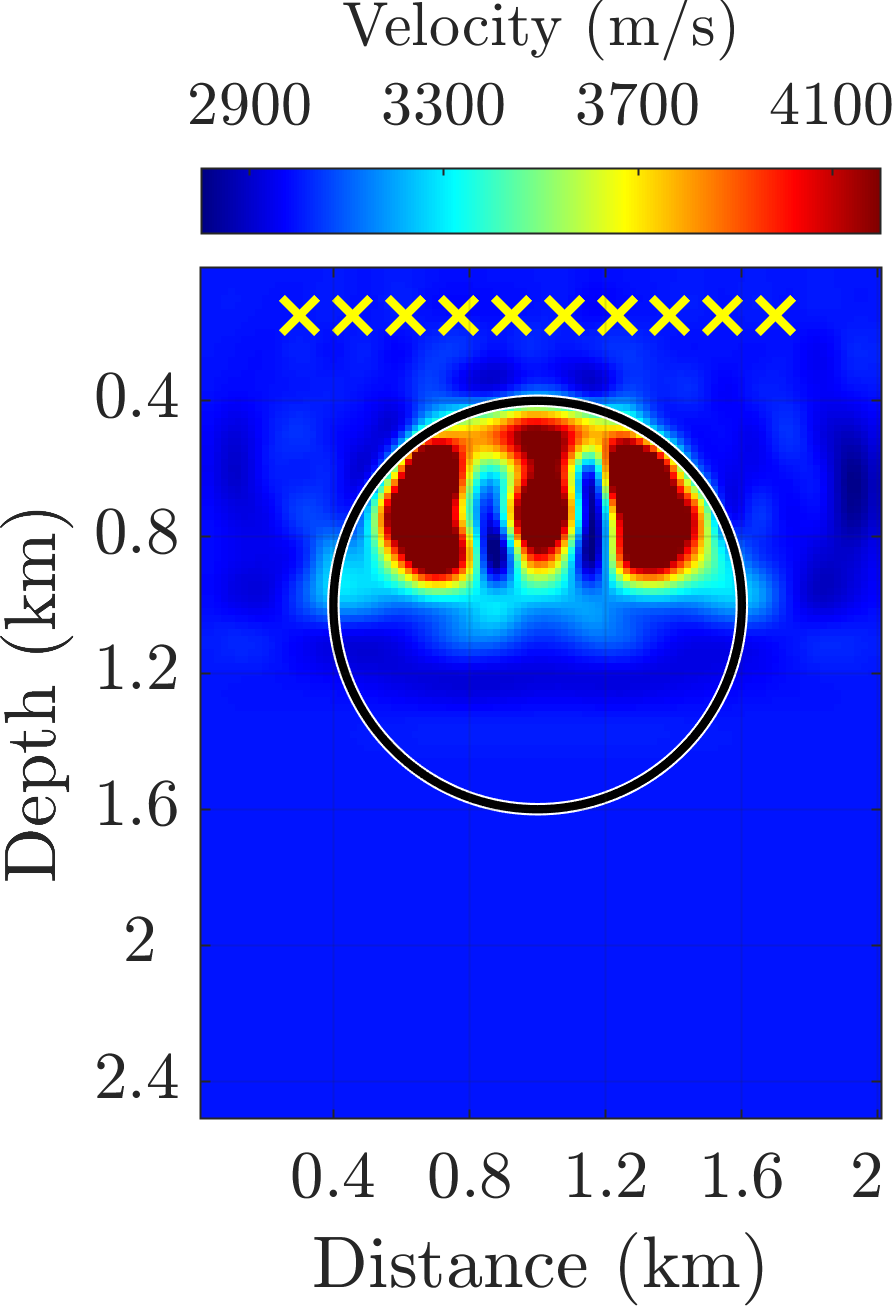} &
\hspace{-0.01\textwidth}\includegraphics[width=0.21\textwidth]
{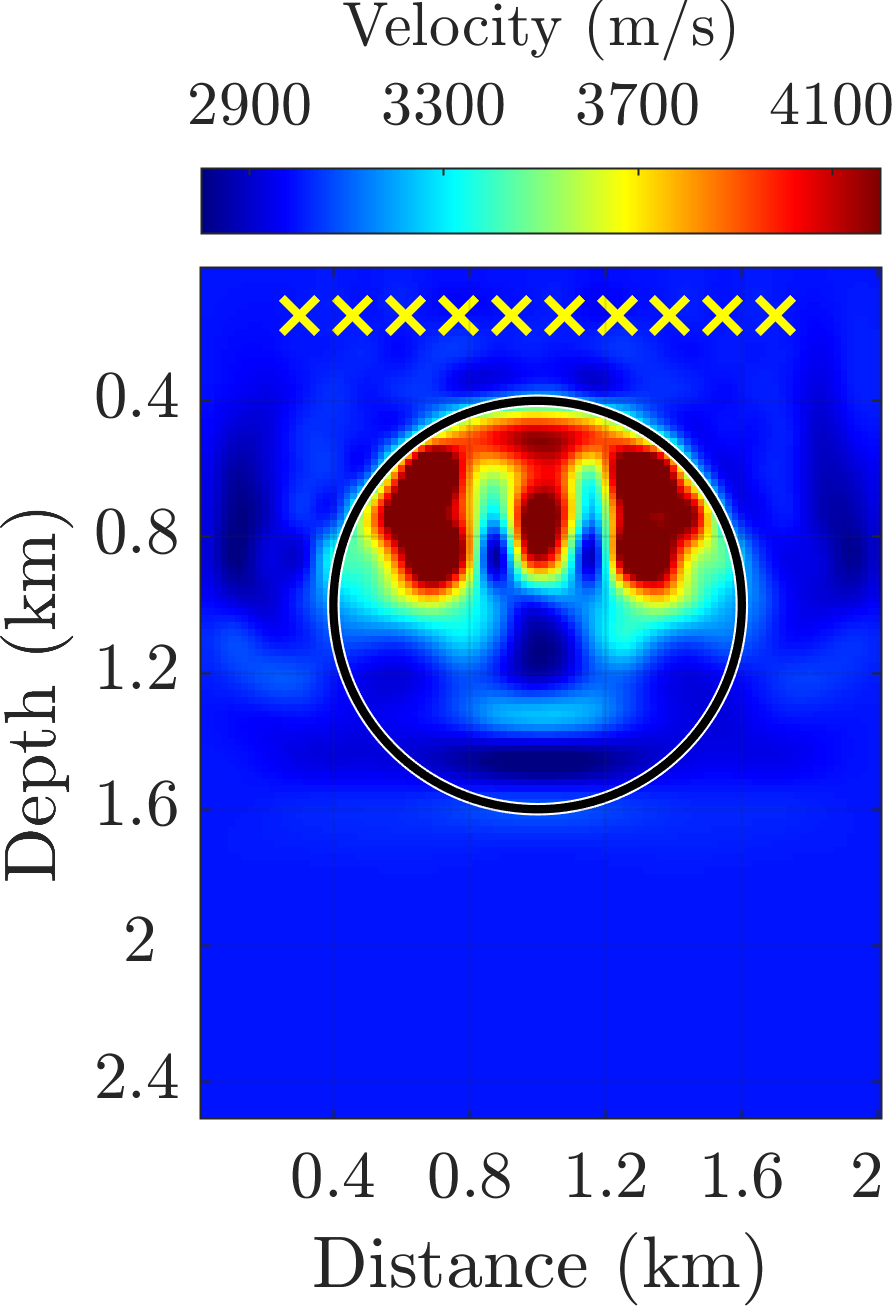} &
\hspace{-0.01\textwidth}\includegraphics[width=0.21\textwidth]
{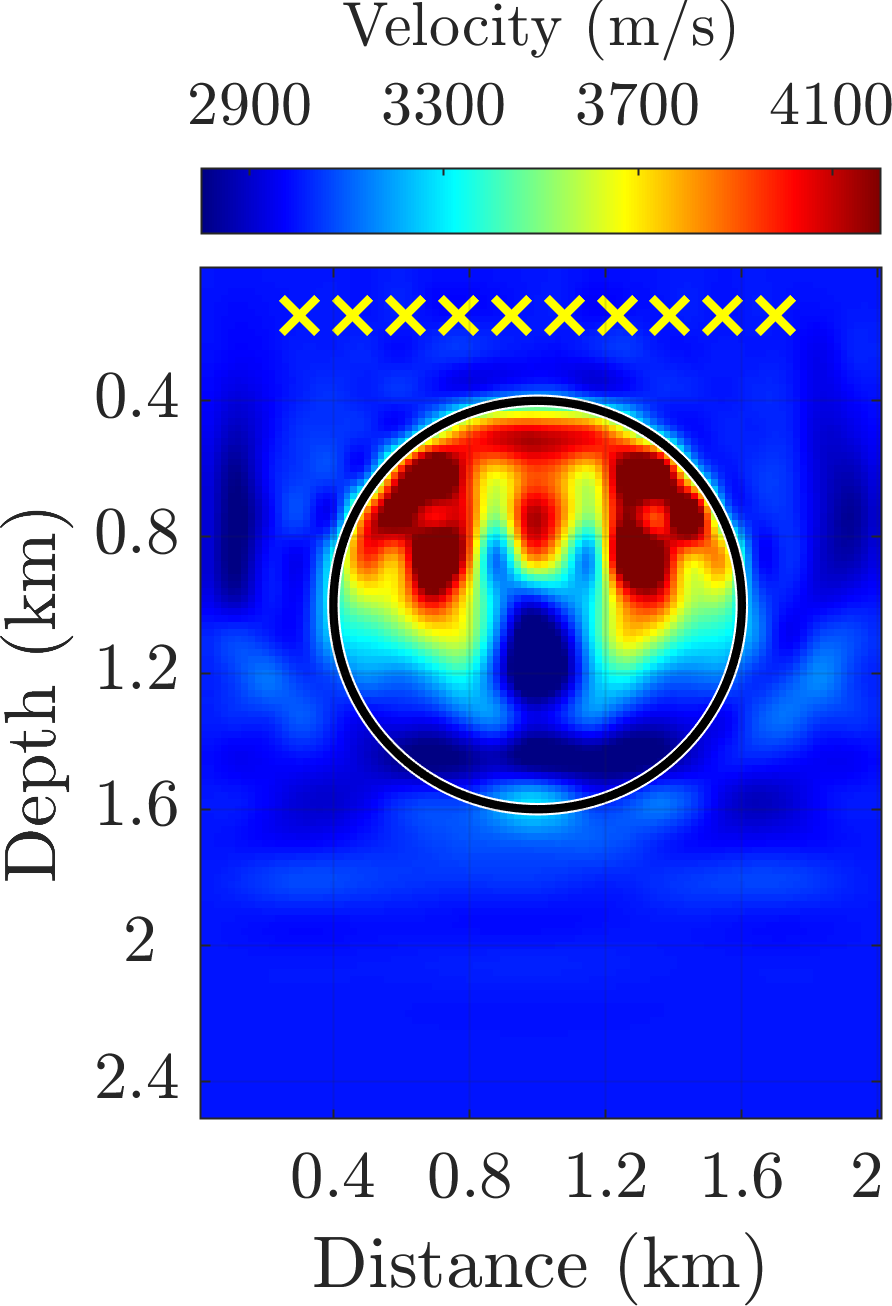} &
\hspace{-0.01\textwidth}\includegraphics[width=0.21\textwidth]
{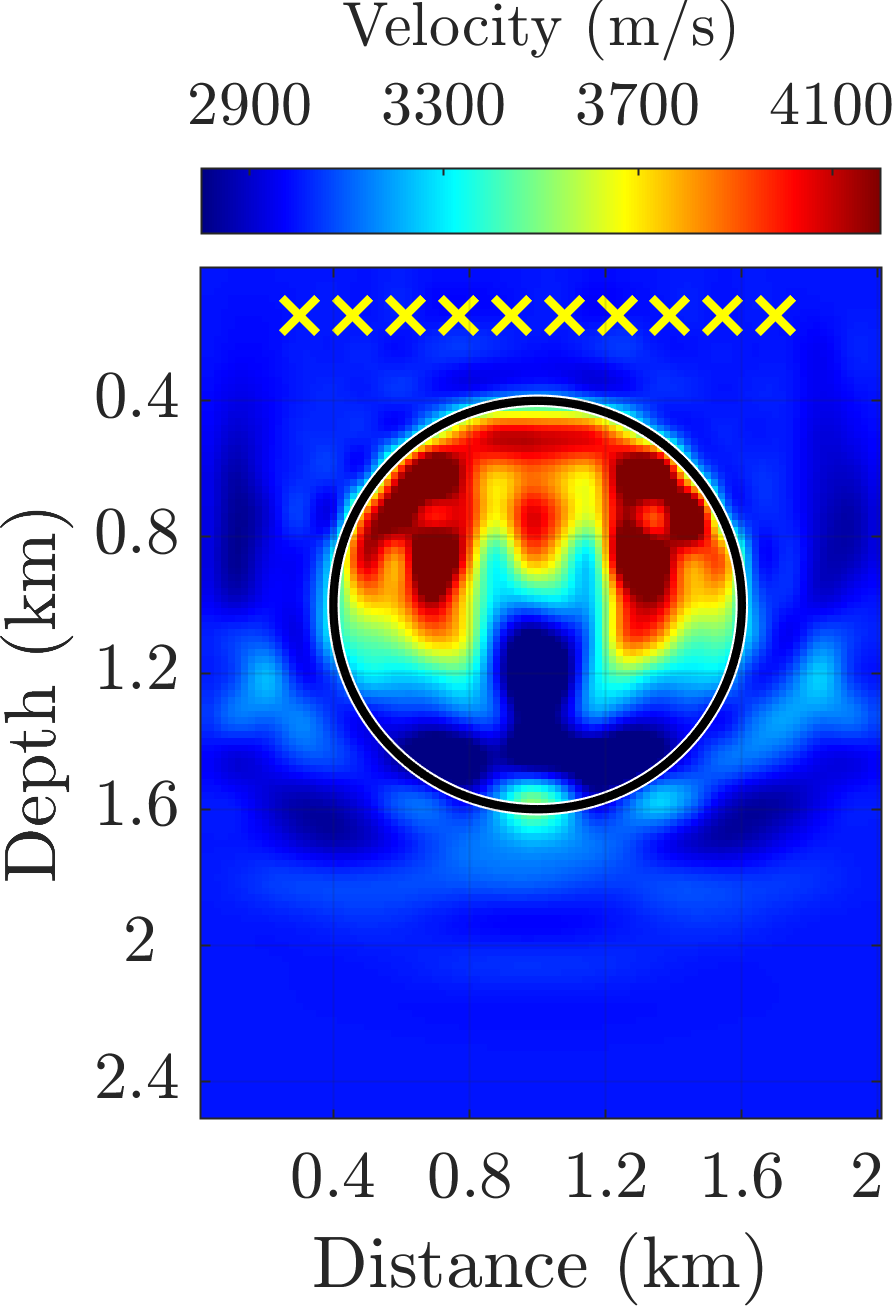} 
\end{tabular}
\end{center}
\caption{Estimated velocity after $10-60$ Gauss-Newton iterations: 
(a)--(d) ROM based velocity estimates;
(e)--(h) FWI velocity estimates.
The true inclusion boundary is shown as a black circle. 
All $\m=10$ {colocated sources/receivers} are shown as yellow $\times$. 
Velocity colorbars are in $\rm{m/s}$, all plots share the same color scale. 
%{Remove colorbar in  (e)-(h) to save space??}
}
\label{fig:CamembertRes}
\end{figure*}

We follow \cite[]{yang2018application} and model the ``Camembert" inclusion as a disk with 
radius of $600$~m, centered at point $(1\text{ km},1\text{ km})$ in the domain 
$\Omega = [0, 2\text{ km}] \times [0, 2.5\text{ km}]$. The setup is illustrated in Figure~\ref{fig:Camembert}, 
where $c(\bx)$ equals $4000\text{ m/s}$ in the inclusion and $3000\text{ m/s}$ outside. 
The data sampling interval is $\tau = 0.0435~\mbox{s}$, 
$\m = 10$ and $\n = 16$.

The search space ${\cal C}$ has dimension $N = 20 \times 20 = 400$, and the velocity is 
parametrized as in equation \eqref{eq:IM2}, with the constant initial guess 
$c_o(\bx) = \bar c = 3000\text{ m/s}$ and the Gaussian basis functions
\begin{equation}
\label{eq:Gaussphi}
\phi_{l}(\bx) = \frac{1}{2 \pi \sigma_\phi \sigma_\phi^\perp} 
\exp \Big[-\frac{(x^\perp-x_l^\perp)^2}{2 (\sigma_\phi^\perp)^2} 
-\frac{(x-x_l)^2}{2 \sigma_\phi^2}\Big],
\end{equation}
with standard deviation {$\sigma_\phi^\perp = 55.5$~m} in the {horizontal (distance)
direction  and {$\sigma_\phi = 69.4$~m} in depth.  
Here we use the system of coordinates $\bx = (x^\perp,x)$, with depth coordinate $x$ and distance coordinate $x^\perp$ orthogonal to it. }
The centers of the Gaussians are at the locations $\bx_l = (x_l^\perp,x_l)$ on a uniform $20\times 20$ 
grid that discretizes the imaging domain 
$\Omega_{\rm im} = [95\text{ m}, 1905\text{ m}] \times [119 \text{ m}, 2381\text{ m}] \subset \Omega.$
Note that $2 \sigma_\phi$ and $2 \sigma_\phi^\perp$ are smaller than half the wavelength 
$\bar c / (10 \text{ Hz}) = 300$~m corresponding to the essential Nyquist frequency. 
Hence, the velocity is over-parametrized and we stabilize the inversion with the adaptive Tikhonov 
regularization  described in Appendix \ref{app:implem}. 

We show in Figure~\ref{fig:CamembertRes}a--\ref{fig:CamembertRes}d the velocity estimates obtained with 
Algorithm~\ref{alg:prowi}, implemented with $\ell = 9$, the number of iterations per layer $\q = 4$, 
and with the restriction parameter $d = \n$. 
The plots in Figure~\ref{fig:CamembertRes}e--\ref{fig:CamembertRes}h are the velocity estimates obtained with  the FWI approach, 
which minimizes the objective function 
\begin{align}
\mathcal{L}^\FWI_i(\bet) & = \mathcal{O}^\FWI \big({v}(\cdot; \bet) \big)+ \mu_i^\FWI \|\bet\|^2_2,
\label{eq:FWIobjE}
\end{align}
with the same time windowing of the data as in the ROM based estimation. 
The Tikhonov regularization parameter $\mu_i^\FWI$ is computed as explained in Appendix \ref{app:implem}.

The results show that the ROM approach gives a much better estimate of $c(\bx)$. 
This estimate improves as we iterate, and by the time we reach the $60^{\rm th}$ step, 
the circular inclusion is reconstructed well. The FWI approach does not improve much after the 
$10^{\rm th}$ step, indicating that the optimization is stuck in a local minimum. 
While the top and arguably the bottom of the inclusion are correctly located, FWI fails to fill in the inclusion
with the correct velocity, overestimating it in the upper half of the disk and underestimating it in the lower half.

\begin{figure*}[ht]
\begin{center}
\begin{tabular}{cccc}
(a) $u^{(5)}(4\tau,\bx)$ & 
(b) ${\rm v}^{(5)}(4 \tau,\bx)$ &
(c) $u^{(5)}(4\tau,\bx)$ for $c(\bx) = \bar c$ & 
(d) ${\rm v}^{(5)}(4 \tau,\bx)$  for $c(\bx) = \bar c$ \\
\includegraphics[width=0.227\textwidth]{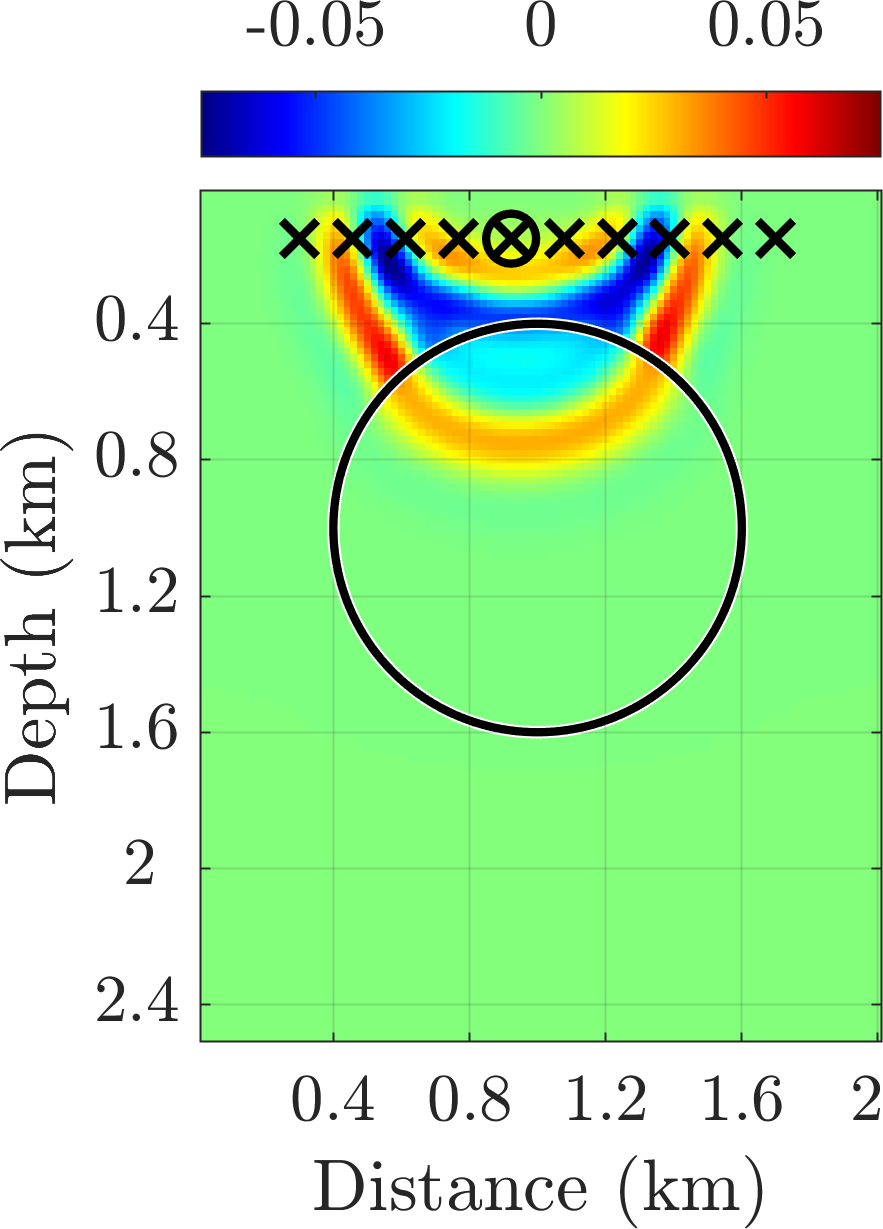} &
\includegraphics[width=0.227\textwidth]{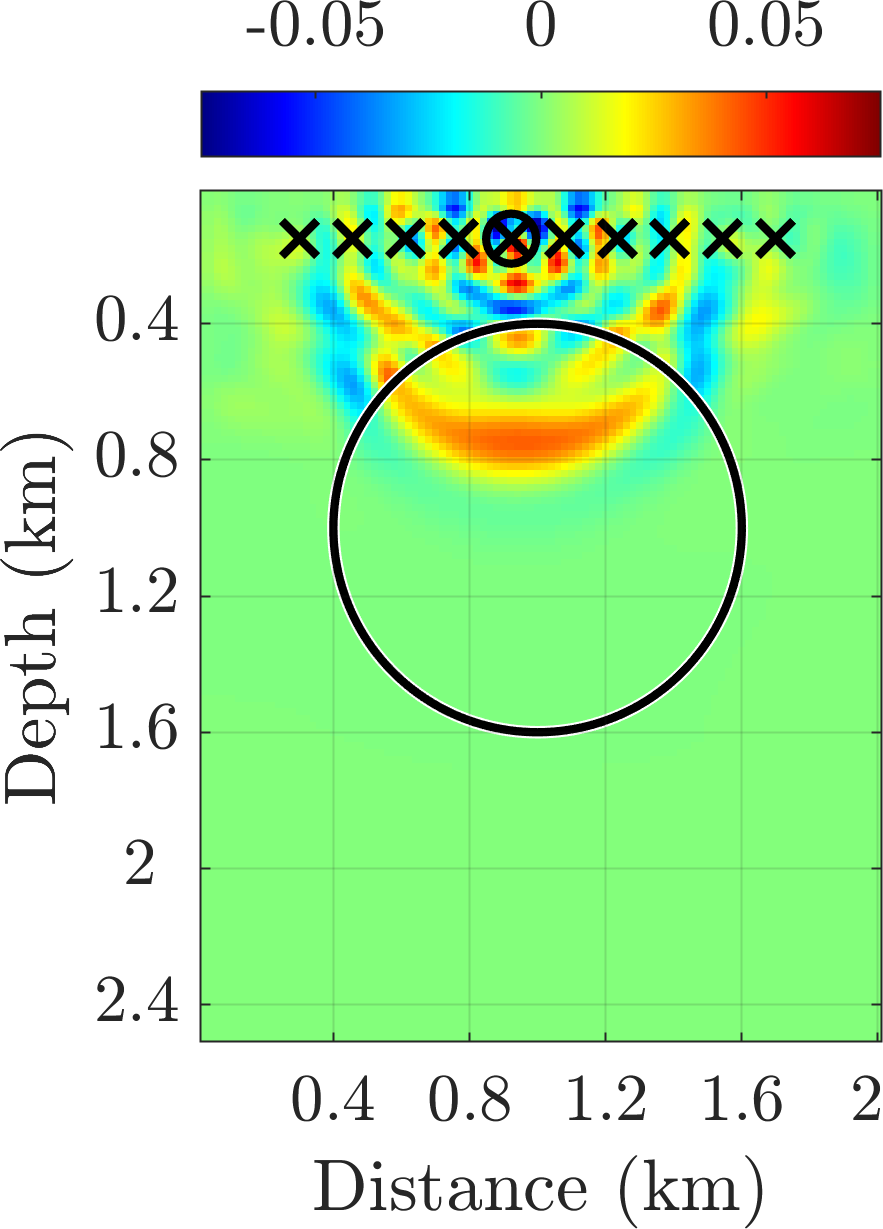} &
\includegraphics[width=0.227\textwidth]{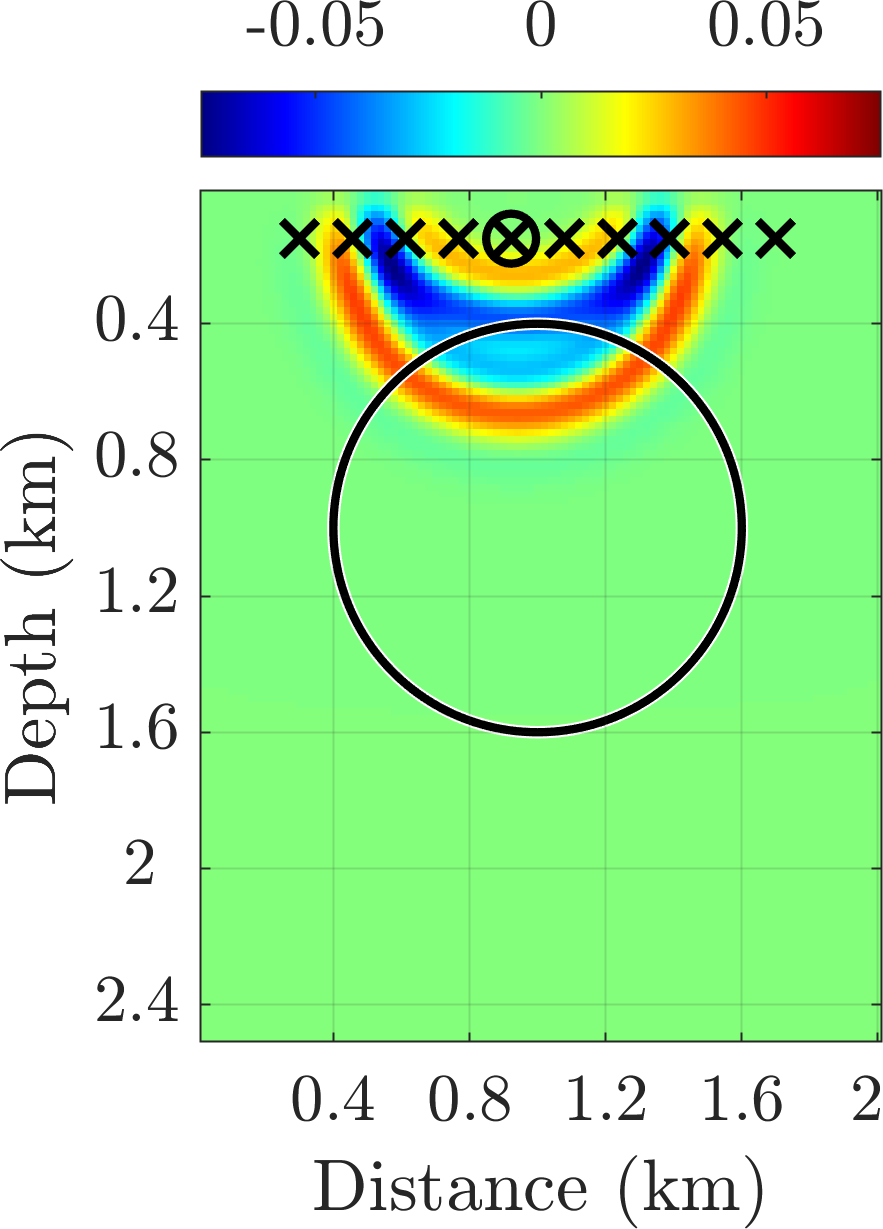} &
\includegraphics[width=0.227\textwidth]{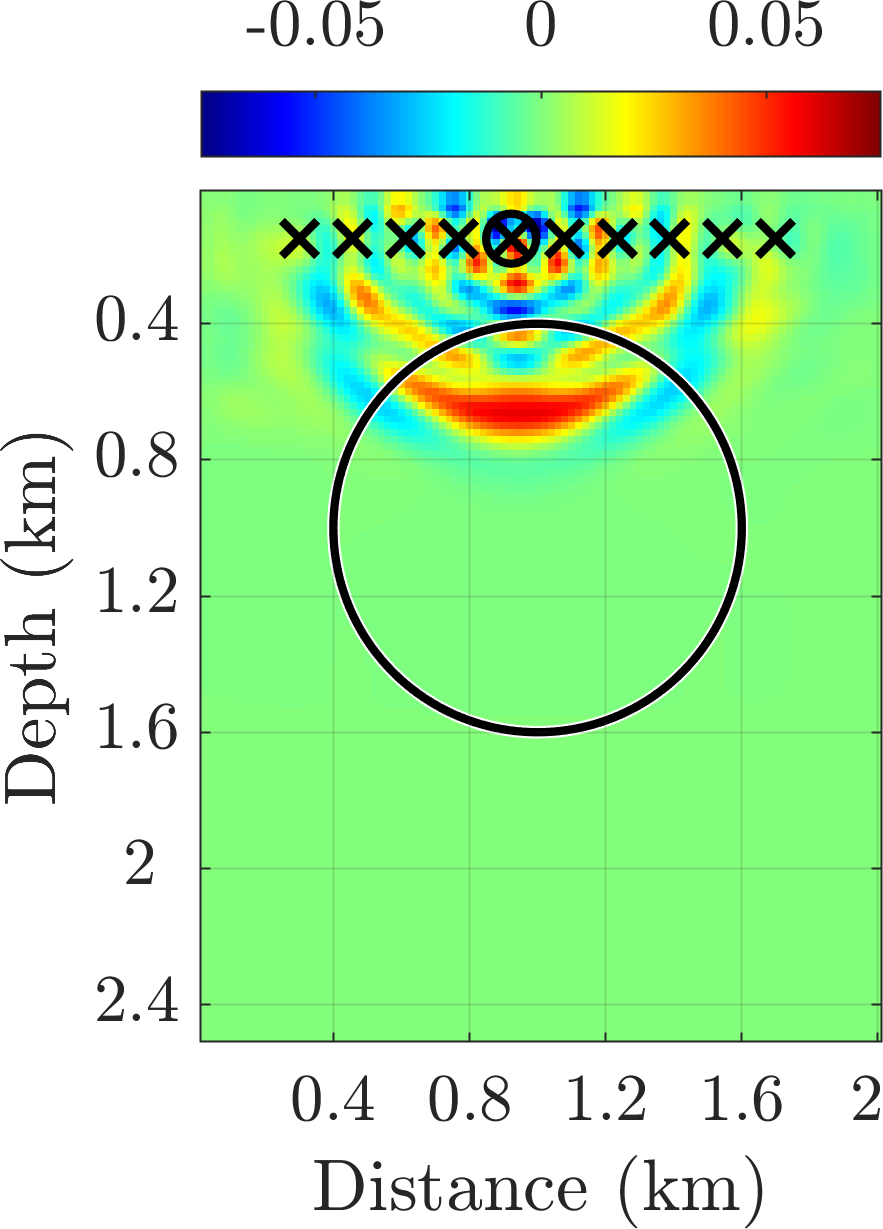} 
\end{tabular}
\end{center}
\vspace{-0.2in}
\caption{Wavefield snapshots and orthonormal basis components at time {instance}  $t = 4 \tau$, 
corresponding to the center left source, indexed by $s = 5$, shown as a black $\bigcirc$:
(a)--(b) plots for the true velocity $c(\bx)$ displayed in Figure~\ref{fig:Camembert}; 
(c)--(d) plots for the reference medium, with $c(\bx) \equiv \bar c = 3000$~m/s. 
The $\m=10$ {colocated sources/receivers} are shown as black $\times$. All the plots share the same color scale.}
\label{fig:Basis}
\end{figure*}

\subsection{Illustration of the orthonormal basis}

We display in Figure~\ref{fig:Basis}a the snapshot $u^\ss(4 \tau,\bx)$ in the medium with the 
Camembert inclusion and in Figure~\ref{fig:Basis}c the snapshot computed with the reference, constant 
velocity $\bar c = 3000$~m/s. The source is in the middle of the array, indicated in the 
plots by the circle, and indexed by $s = 5$. Obviously, the snapshot in the true medium is different 
from the one in the reference medium. In the reference medium, the wave is simply a spherical 
wave emitted by the point source and reflected by the top surface modeled as a sound soft boundary. 
In the true Camembert model medium, the wave is scattered at the boundary and at the top of the inclusion, 
and it travels further down for the same  $t = 4\tau$, due to the fast inclusion.

The corresponding components of the orthonormal basis stored in $\bV(\bx)$, called 
${\rm v}^{(5)}(4 \tau,\bx)$, are shown in Figures~\ref{fig:Basis}b and \ref{fig:Basis}d. 
They illustrate the second and third attributes of the orthonormal basis, stated in the outline of our velocity estimation 
method. Indeed, the basis function in the true and reference medium are very similar. 
They both have a localized peak near the deepest point reached by the wave at {instance}  $t = 4 \tau$ 
and they are oscillatory away from it. The scattering at the top of the inclusion does not have a strong 
effect on the basis function, but the kinematics makes a difference. As mentioned above, the wave penetration at 
$t = 4 \tau$ is deeper in the true medium, due to the fast inclusion, so the localized peaks are in different locations. 

\section{Velocity estimation with noisy and towed-streamer data}
\label{sect:noisy}

\begin{figure*}[ht]
\begin{center}
\begin{tabular}{ccc}
(a) Marmousi model & (b) ROM estimate, Gaussian basis & (c) ROM estimate, hat basis \\
\hspace{-0.01\textwidth}\includegraphics[width=0.33\textwidth]
{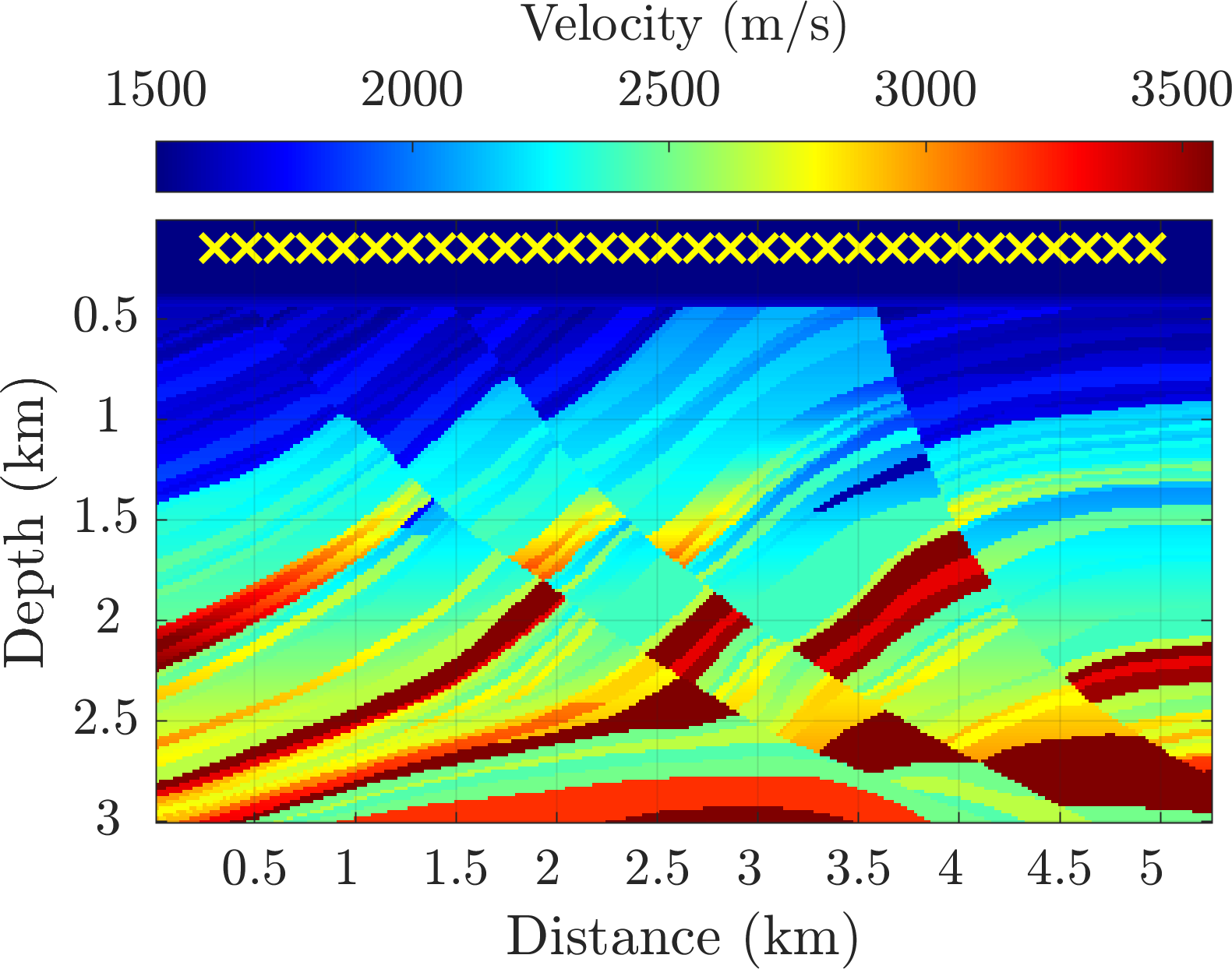} &
\hspace{-0.02\textwidth}\includegraphics[width=0.33\textwidth]
{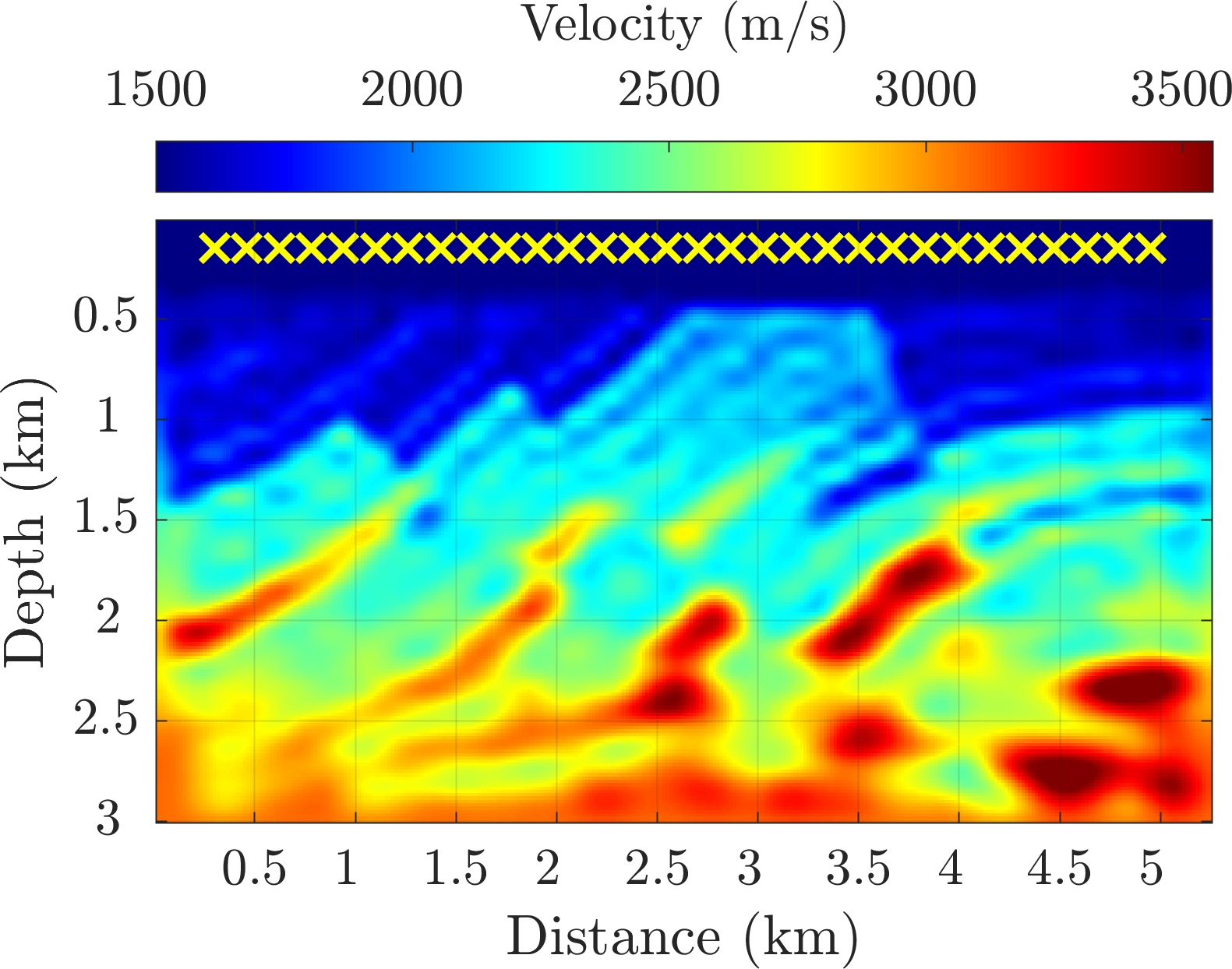} &
\hspace{-0.02\textwidth}\includegraphics[width=0.33\textwidth]
{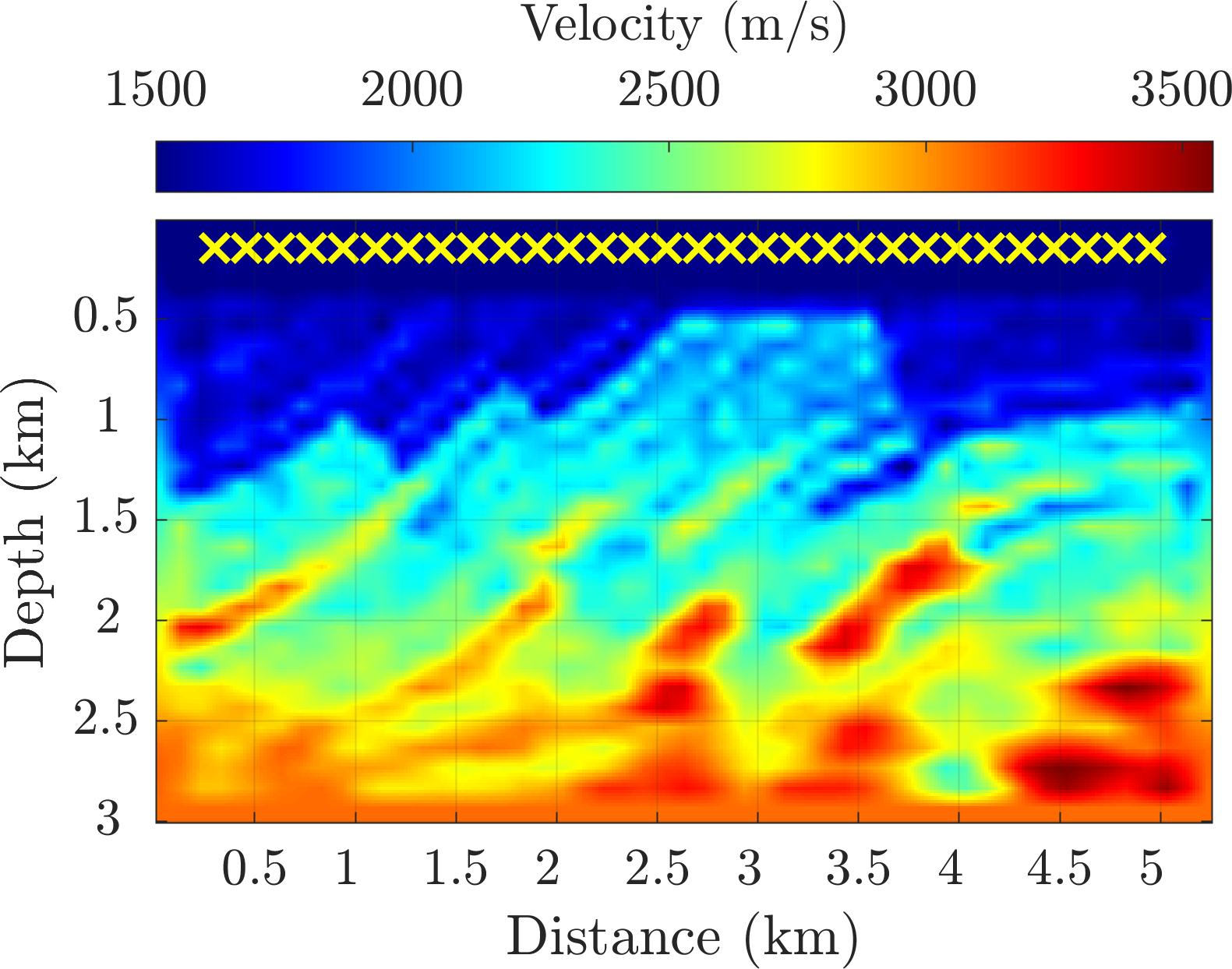} \\
(d) Initial model & (e)  Towed-streamer ROM estimate & (f) Refined ROM estimate \\
\hspace{-0.01\textwidth}\includegraphics[width=0.33\textwidth]
{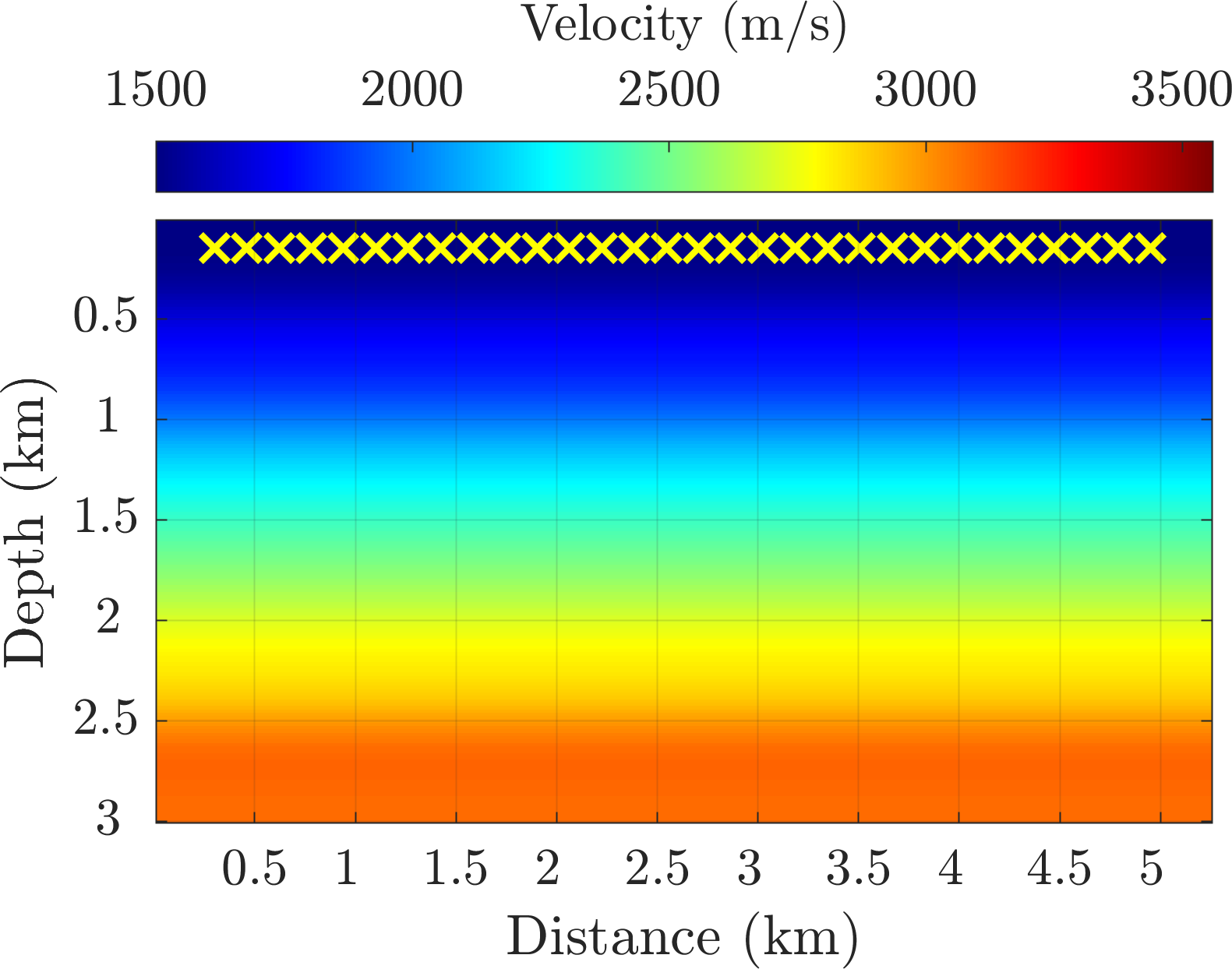} &
\hspace{-0.02\textwidth}\includegraphics[width=0.33\textwidth]
{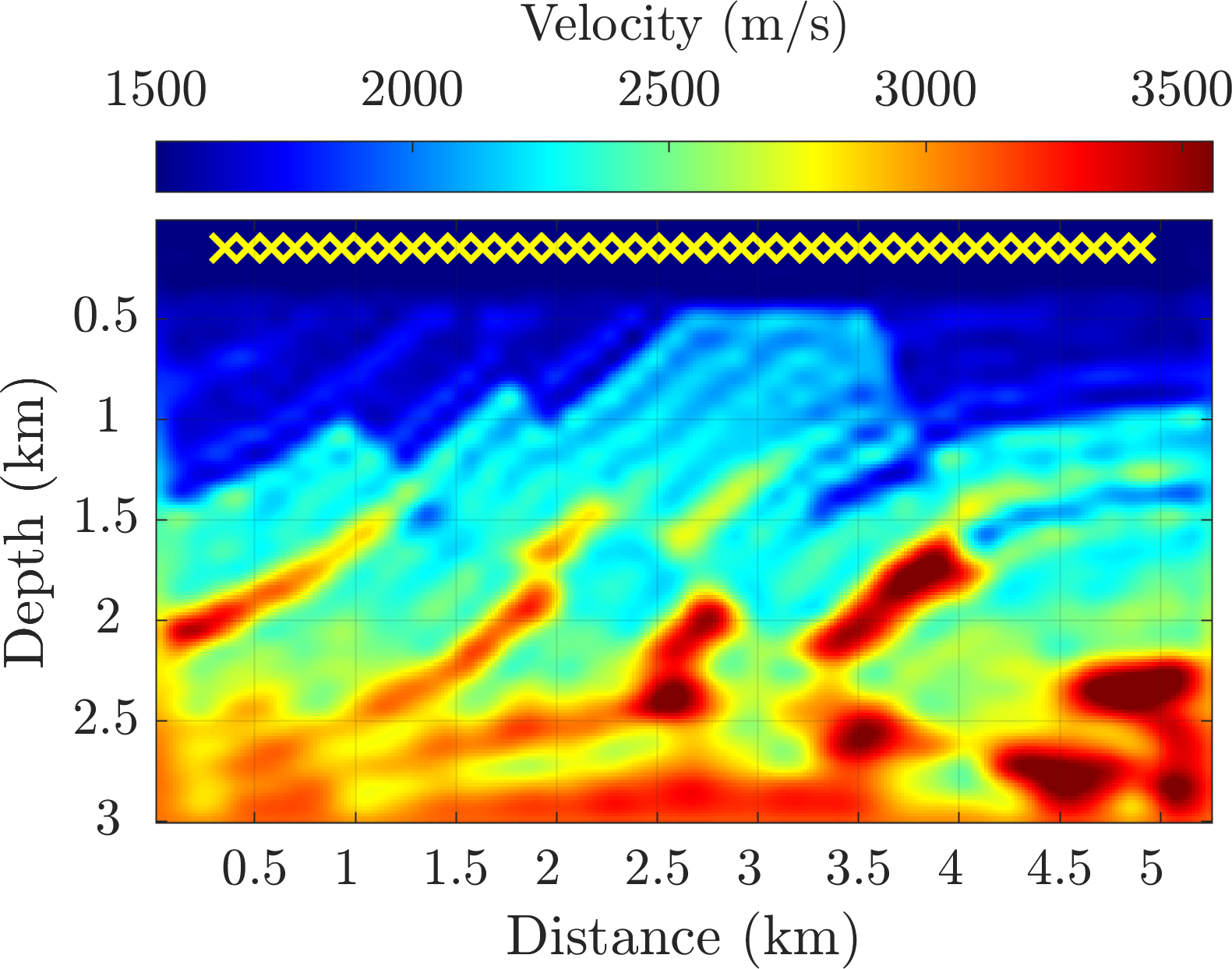} &
\hspace{-0.02\textwidth}\includegraphics[width=0.33\textwidth]
{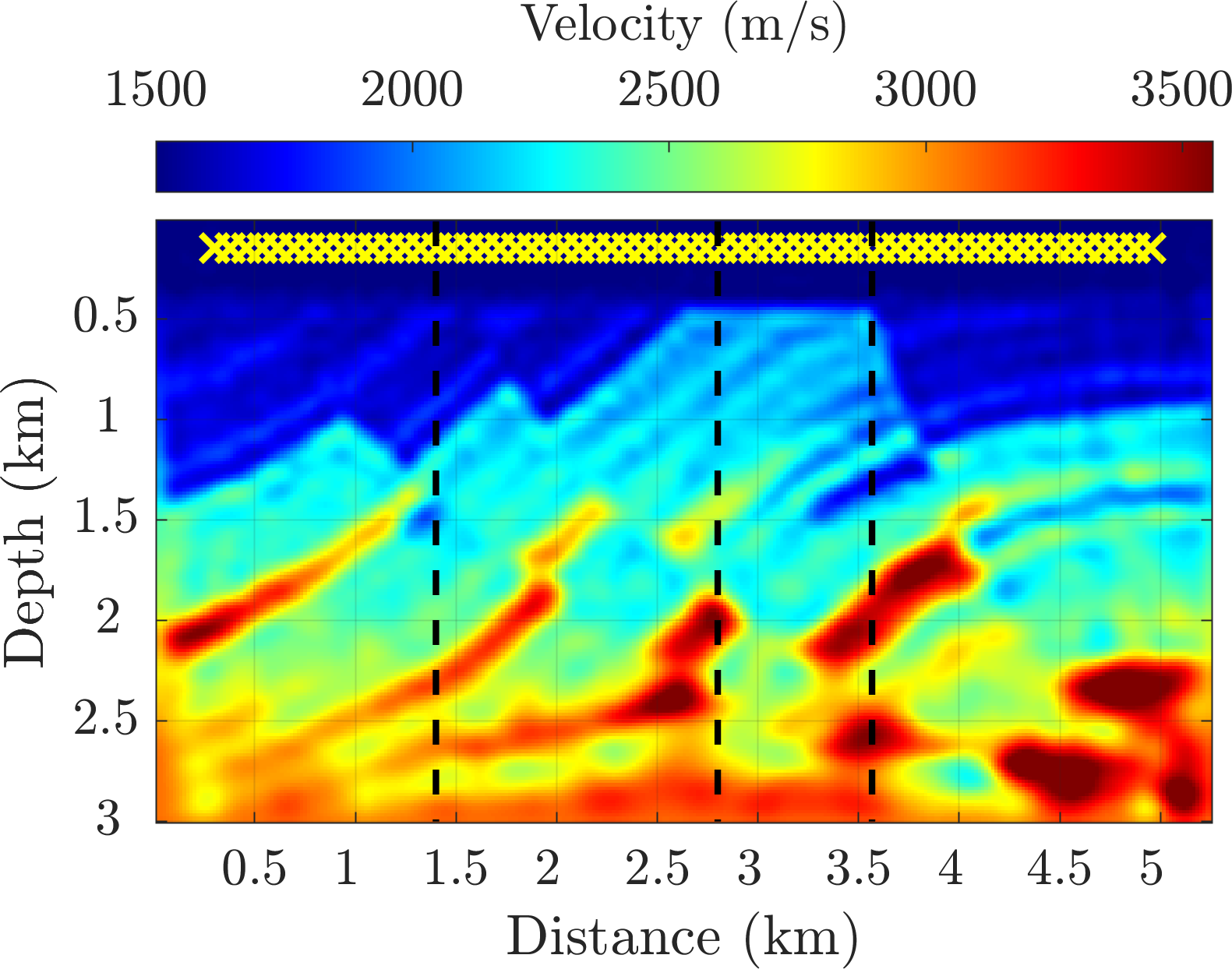} \\
\end{tabular}
\end{center}
\vspace{-0.2in}
\caption{ROM based velocity estimates for Marmousi model with noisy data and towed-streamer measurements:
(a) The section of the Marmousi model; 
(b) Velocity estimate from noisy data with Gaussian basis functions parametrization of $v$;
(c) Velocity estimate from noisy data with hat basis functions parametrization of $v$;
(d) Initial guess model $c_o(\bx)$;
(e) Velocity estimate from towed-streamer measurements;
(f) Velocity estimate refinement from data gathered on a dense array sensors and at small time interval $\tau$.
All the {sources/receivers}, $\m=30$ in (a)--(d), $\m=40$ in (e) and $\m=60$ in (f), are shown as yellow $\times$. 
Velocity colorbars are in $\rm{m/s}$. All plots share the same color scale. 
%{Take out the colorbars in the second row to save space?}
}
\label{fig:Marmousi}
\end{figure*}

In this section we present velocity estimation results with noisy measurements and with the 
array response matrix $\boldsymbol{\cal M}(t)$ assembled from towed-streamer type measurements.  
In both cases we have uncertainty of the data, which affects the computation of $\bA^\RM$. 
There are two critical steps in Algorithm \ref{alg:arom} that must be addressed, 
and they both involve the mass matrix $\bM$ computed at step 3, which will likely be 
neither symmetric nor positive definite. These properties are needed for the computation of the 
Cholesky square root $\bR$ at step 4 and the inverse $\bR^{-1}$ that gives the output of the algorithm. 
The lack of symmetry is easy to fix, but to ensure the positive definiteness, 
we need a regularization procedure that involves a spectral projection of $\bM$ on the space of its 
leading eigenvectors, corresponding to the significant eigenvalues. These eigenvectors and eigenvalues 
are least affected by the uncertainty. The regularization procedure is not straightforward, because we 
must preserve the causality of $\bA^\RM$ in order for the velocity estimation to succeed. 
We explain it in detail in Appendix \ref{app:RegROM}. 

To assemble the matrix $\boldsymbol{\cal M}(t)$ from towed-streamer measurements, 
we use source-receiver reciprocity on-the-fly to fill in the missing off-diagonal entries in 
$\boldsymbol{\cal M}(t)$. To compute the diagonal entries, corresponding to the source being 
also a receiver, we use interpolation of the values at nearby measurement locations, 
two on the left and two on the right. We use Lagrange polynomial interpolation in the Fourier 
(frequency) domain, for 
\begin{equation}
\int_{\RR} d t \, e^{i \om t} [\boldsymbol{\cal M}(t) - \cF[\bar c](t)].
\end{equation}
Then, we inverse Fourier transform to get $\boldsymbol{\cal M}(t)$.
\subsection{Numerical results}

We do not show the Camembert estimation for uncertain measurements, because the information 
needed to get the good result in Figure~\ref{fig:CamembertRes} requires accurate knowledge of 
$\boldsymbol{\cal M}(t)$. This is not the fault of the inversion method. It is due to the fact that the 
bottom part of the Camembert inclusion gives very weak signal at the array, which is accounted for 
in the small eigenvalues of the mass matrix. Any uncertainty of the data will perturb significantly 
these eigenvalues and the associated eigenvectors, so the ROM inversion is not better than that with FWI.

We present instead velocity estimation results for a section of the Marmousi model shown in 
Figure~\ref{fig:Marmousi}a, where we exclude the portion of the water down to depth $266${~m}. 
The domain is $\Omega = [0,5.25\text{ km}]\times [0,3\text{ km}]$. 
The data sampling for the ROM construction is $\tau = 0.0435$~{\rm s} and the number of snapshots 
that span the approximation space is $\n = 40$. The {colocated sources/receivers} are located underwater at depth $150$~m and 
they emit the same pulse given in equation \eqref{eq:pulse}. We present results in two settings. First, when working with 
noisy data, we employ an array of $\m = 30$ {colocated sources/receivers}, {separated by the distance $166.66$~m}. Second, when working with data approximated 
from towed-streamer type measurements, we use {closely spaced receivers, at $16.66$~m apart, to carry out the interpolation 
of the measurements and fill in the missing zero offset data. Then, we subsample the result before we input it in Algorithm \ref{alg:arom},  by keeping $\m = 40$ {sources/receivers} separated by the distance $116.66$~m.
}

\begin{figure*}[ht]
\begin{center}
\begin{tabular}{ccc}
(a) ROM estimate, iteration $6$ & (b) ROM estimate, iteration $12$ & (c) ROM estimate, iteration $18$ \\
\hspace{-0.01\textwidth}\includegraphics[width=0.33\textwidth]
{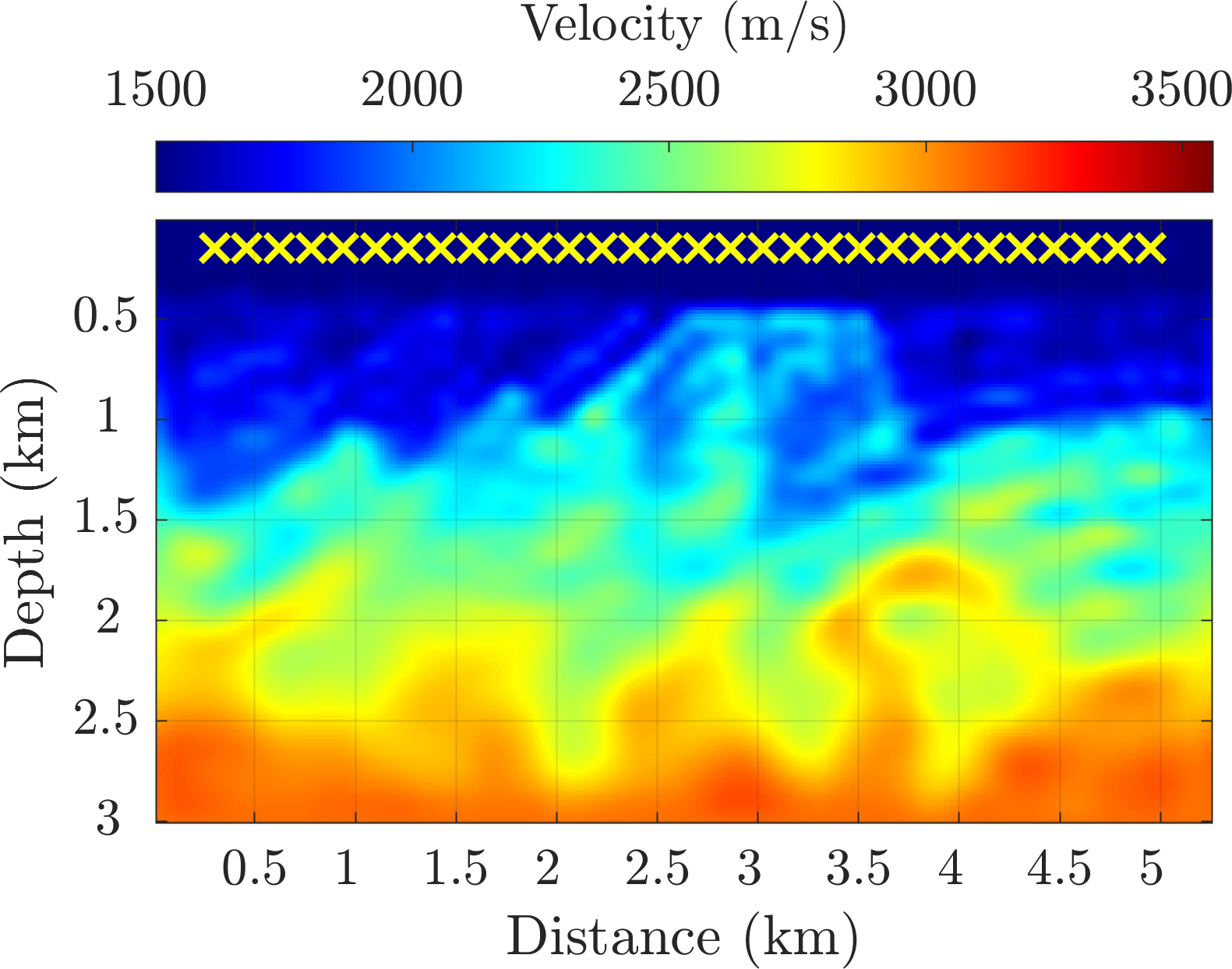} &
\hspace{-0.02\textwidth}\includegraphics[width=0.33\textwidth]
{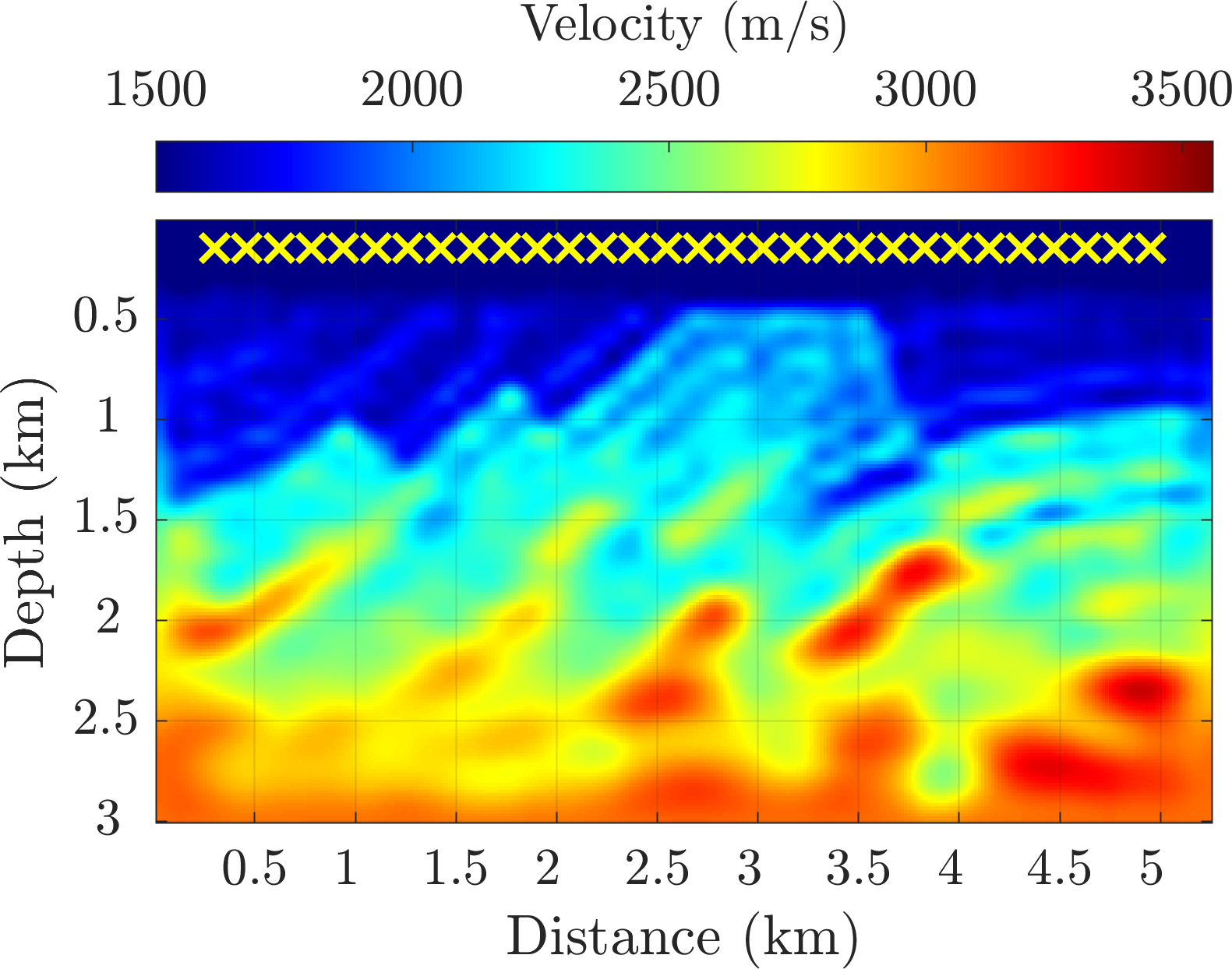} &
\hspace{-0.02\textwidth}\includegraphics[width=0.33\textwidth]
{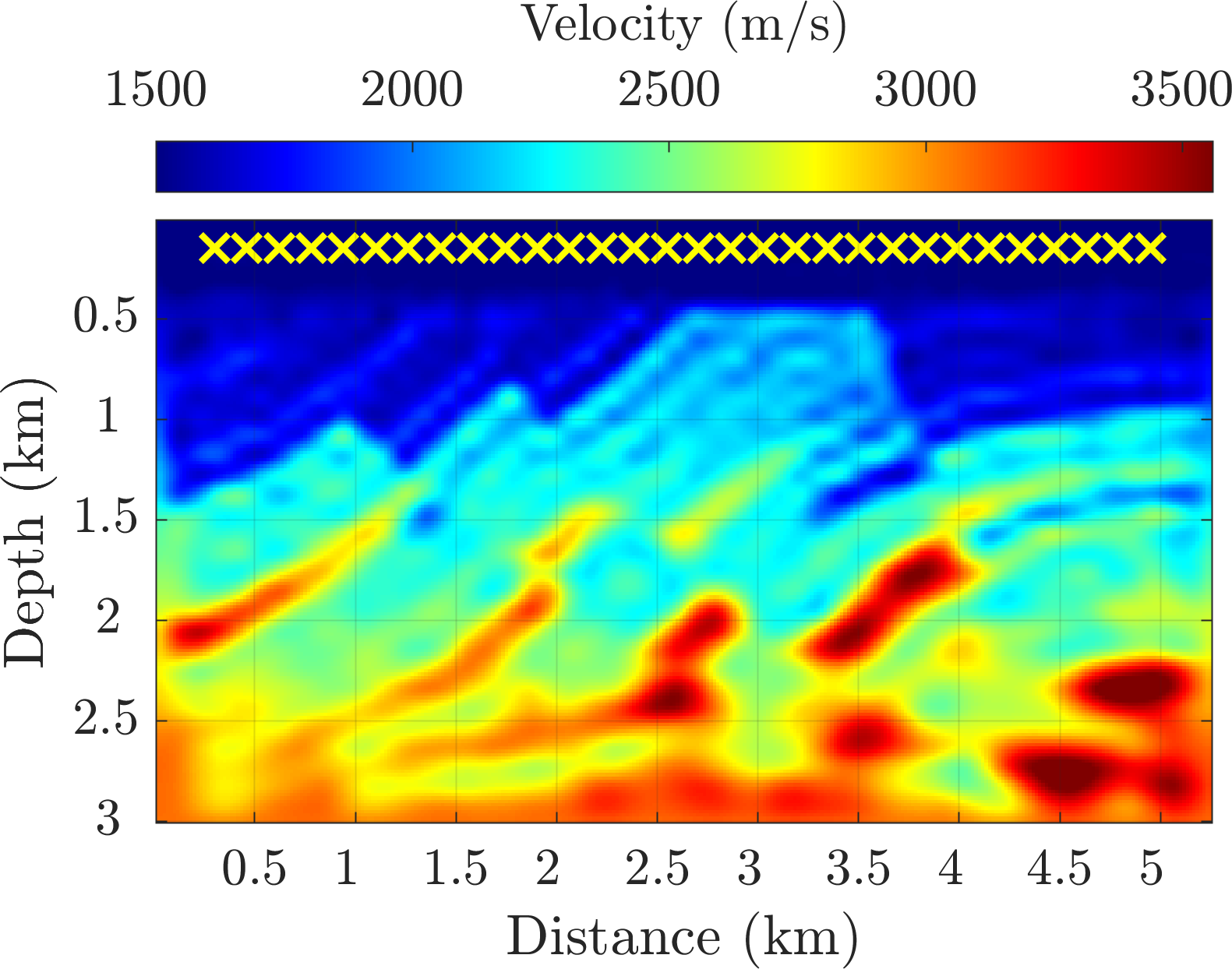} \\
(d) FWI estimate, iteration $6$ & (e) FWI estimate, iteration $12$ & (f) FWI estimate, iteration $18$ \\
\hspace{-0.01\textwidth}\includegraphics[width=0.33\textwidth]
{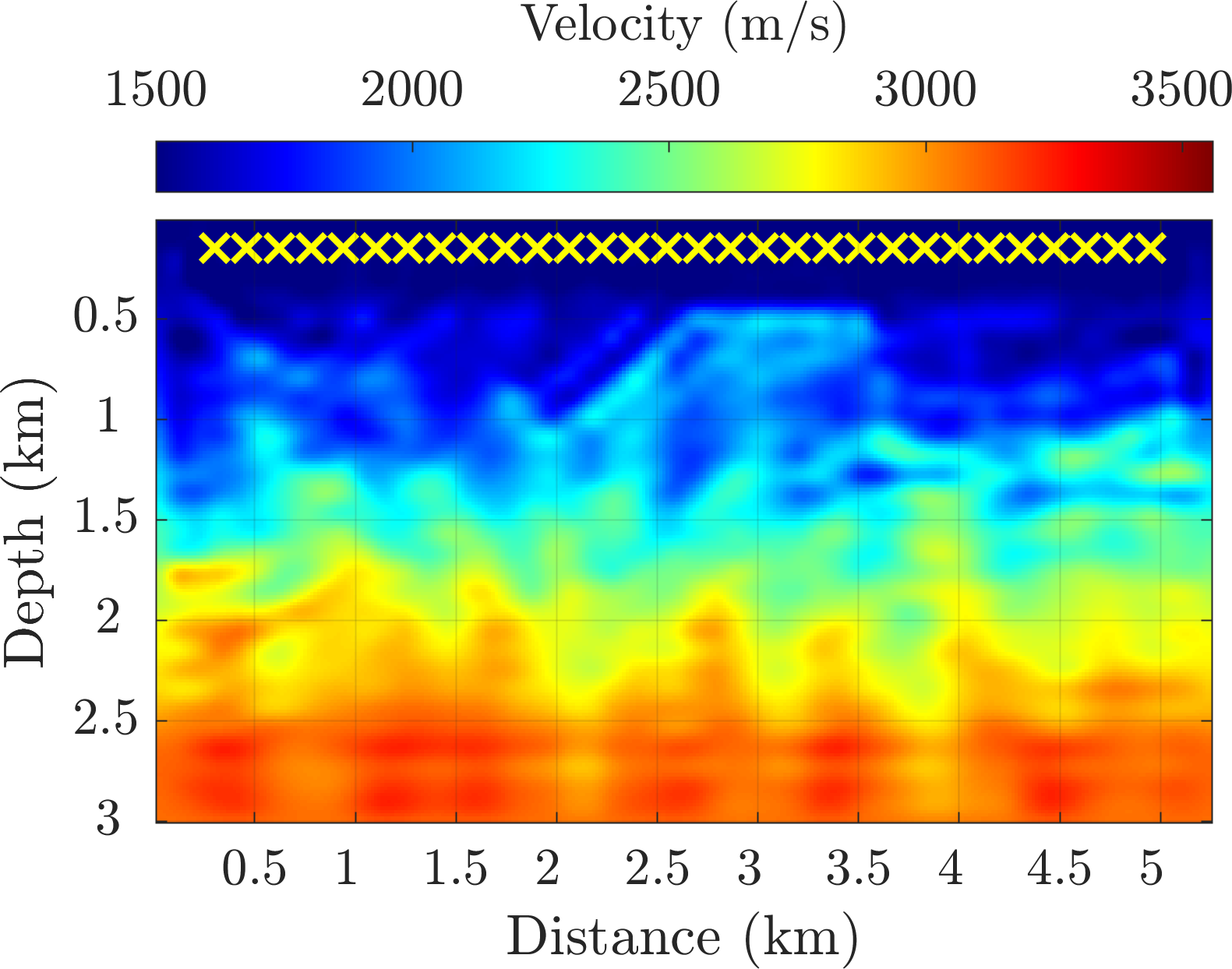} &
\hspace{-0.02\textwidth}\includegraphics[width=0.33\textwidth]
{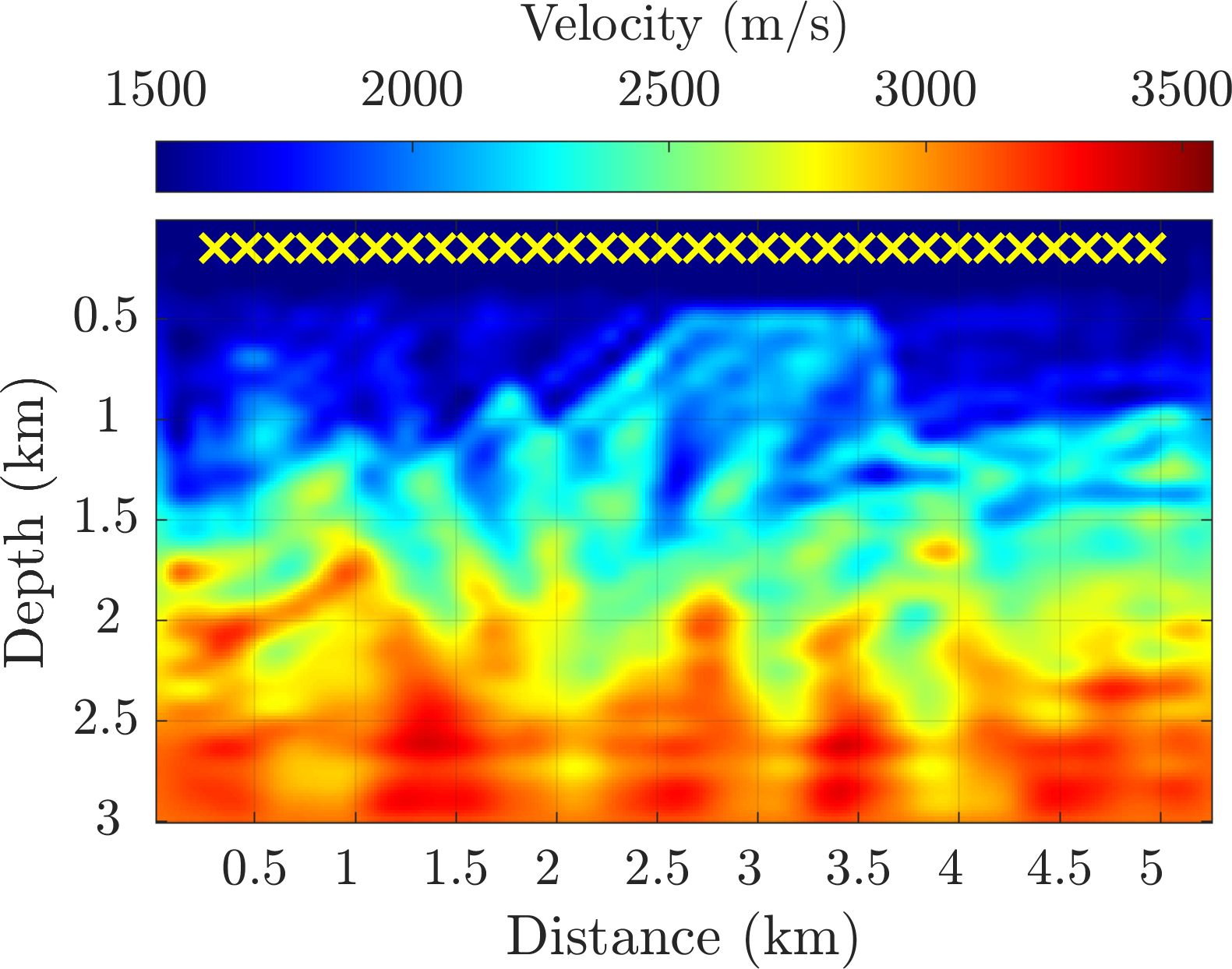} &
\hspace{-0.02\textwidth}\includegraphics[width=0.33\textwidth]
{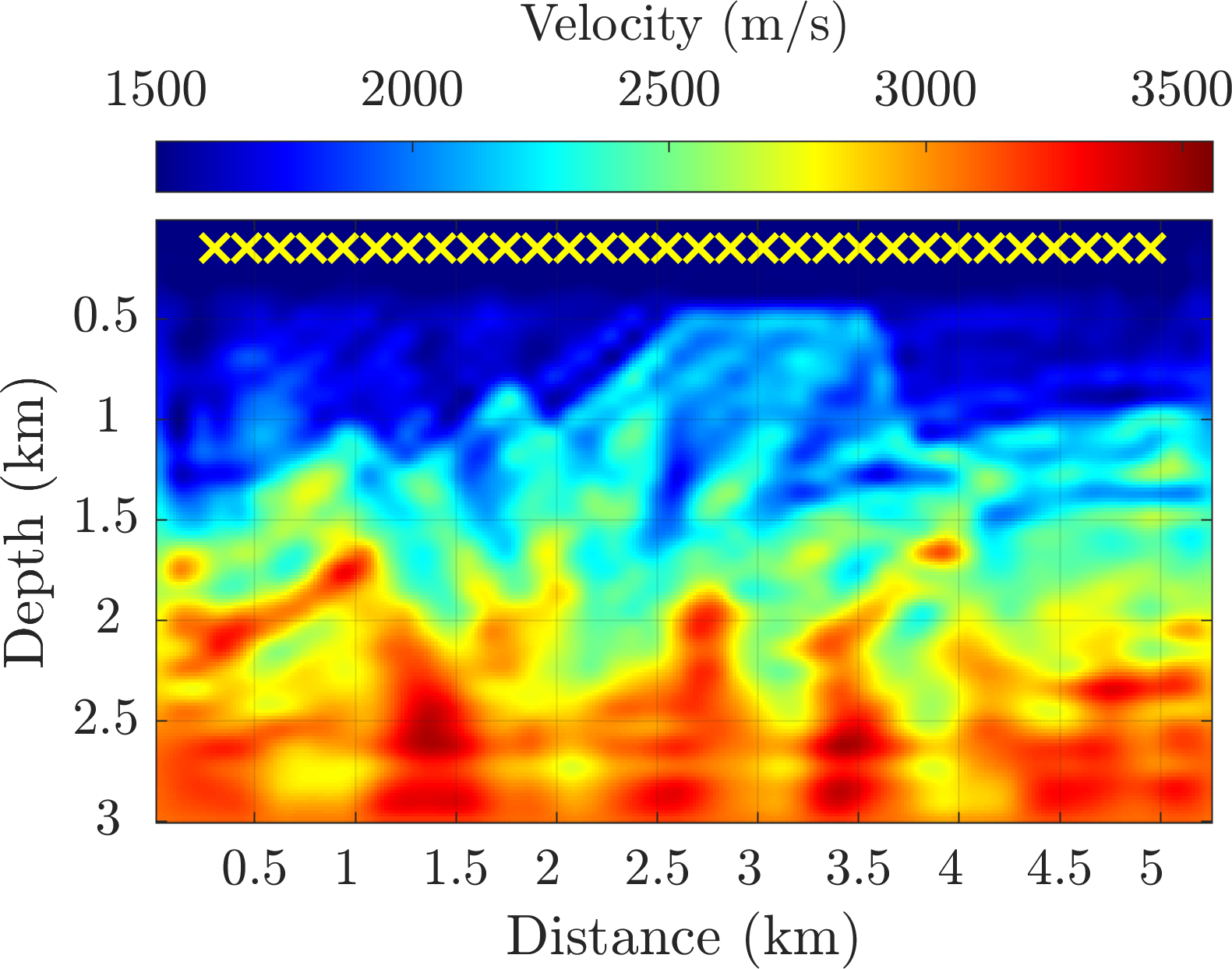} 
\end{tabular}
\end{center}
\vspace{-0.2in}
\caption{Velocity estimates for the Marmousi model with noisy data after $6$, $12$ and $18$ 
Gauss-Newton iterations: (a)--(c) ROM based approach; (d)--(f) FWI approach.
The $\m=30$ {colocated sources/receivers} are shown as yellow $\times$.
Velocity colorbars are in $\rm{m/s}$ and all plots share the same color scale. 
%{Take out the colorbars in the second row to save space?}
}
\label{fig:MarmRes}
\end{figure*}

\begin{figure*}[ht]
\begin{center}
\begin{tabular}{ccc}
(a) Distance $1.4\text{ km}$ & (b) Distance $2.8\text{ km}$ & (c) Distance $3.566\text{ km}$ \\
\hspace{-0.01\textwidth}\includegraphics[width=0.33\textwidth] 
{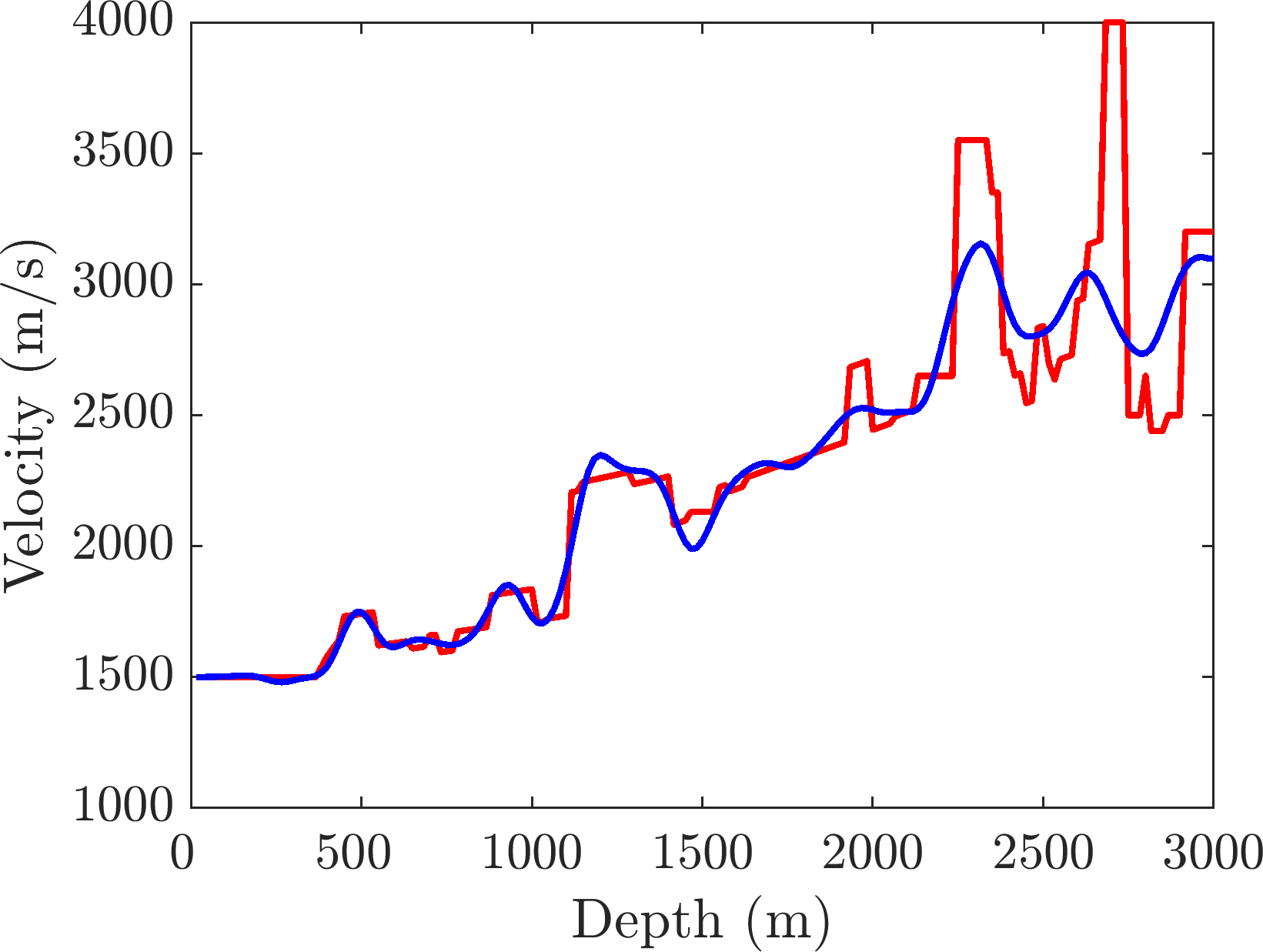} &
\hspace{-0.02\textwidth}\includegraphics[width=0.33\textwidth] 
{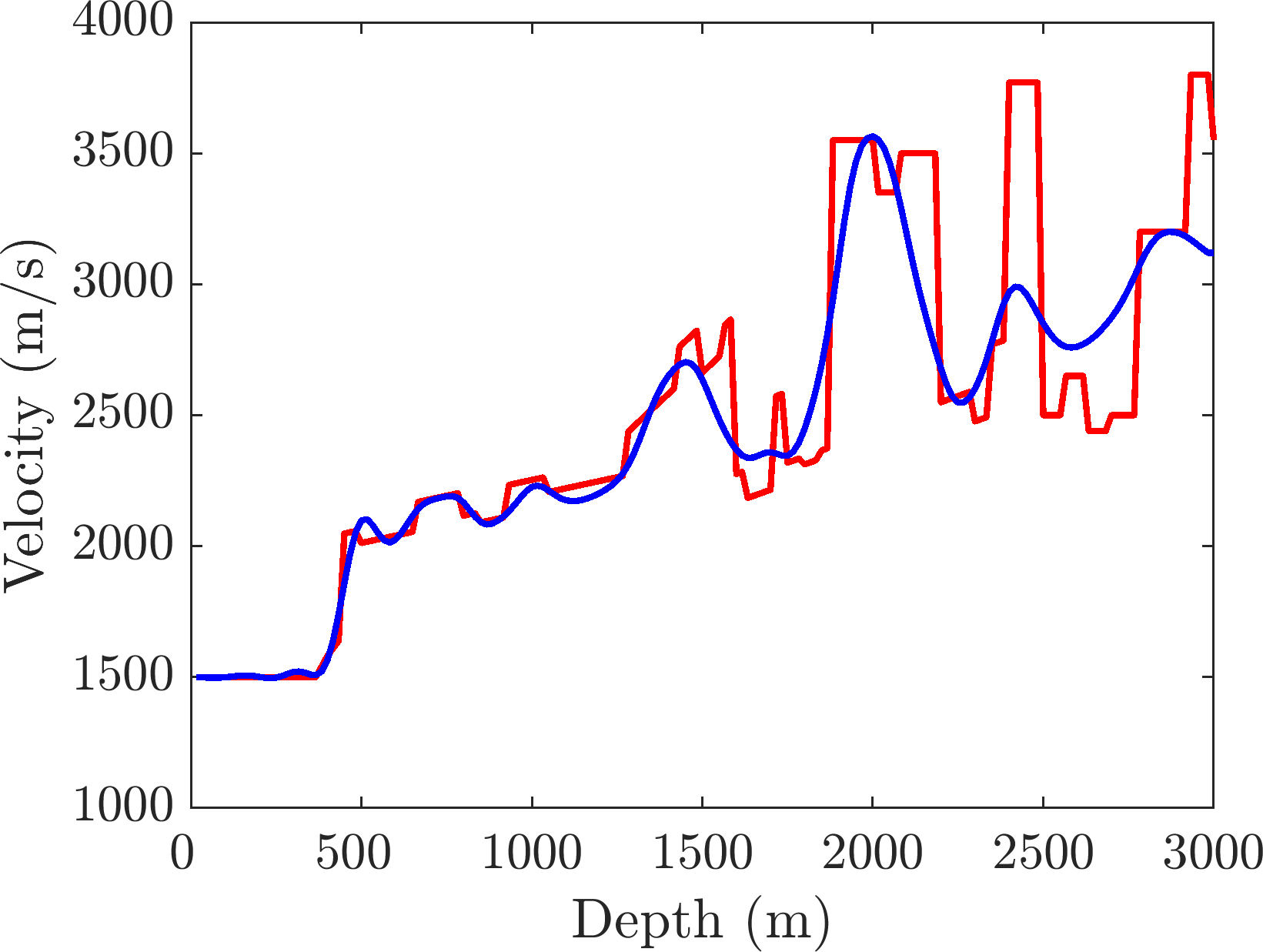} &
\hspace{-0.02\textwidth}\includegraphics[width=0.33\textwidth] 
{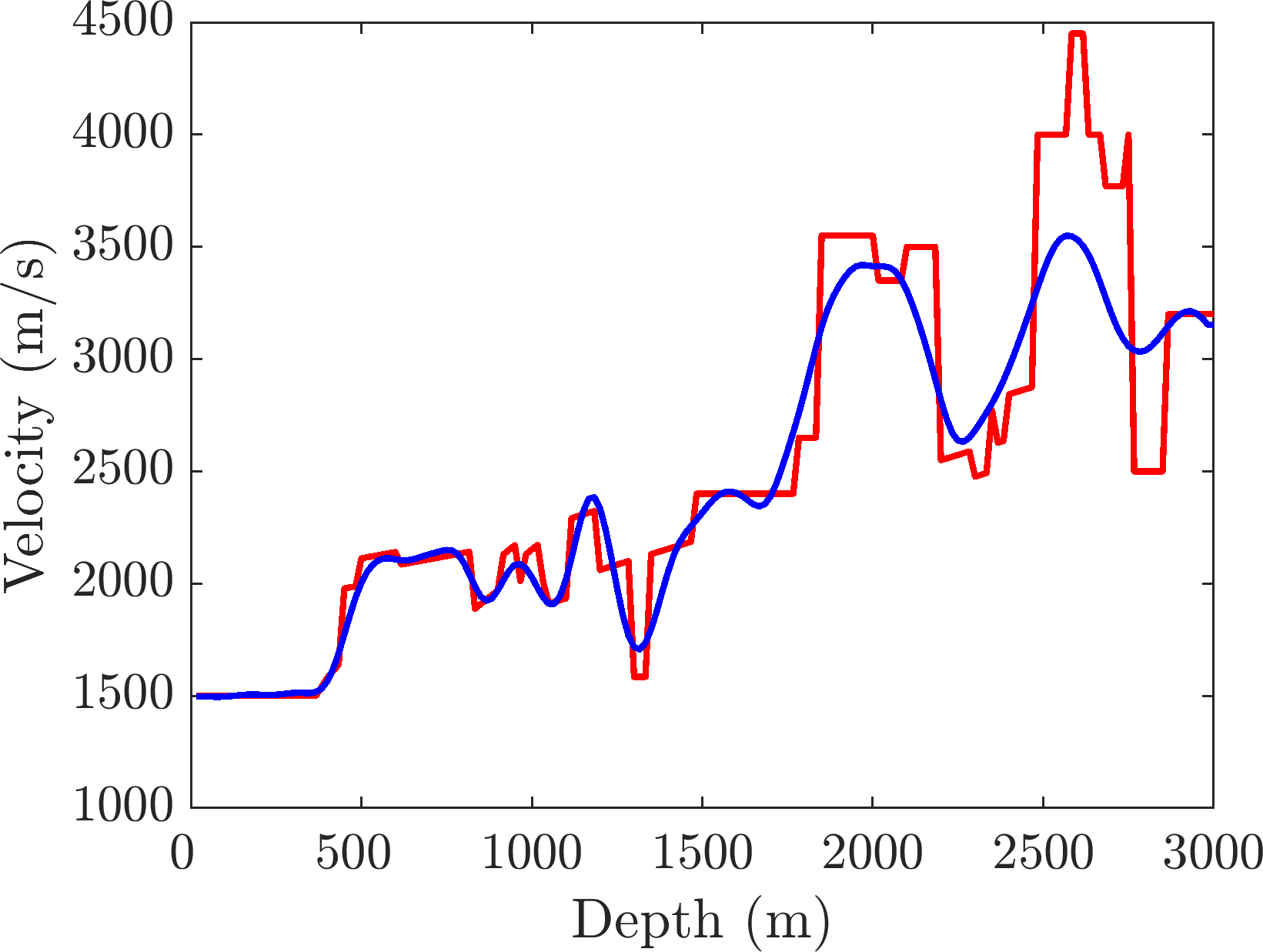}
\end{tabular}
\end{center}
\vspace{-0.2in}
\caption{Vertical slices of the Marmousi model velocity (red lines) and its refined ROM estimate (blue lines) 
at distances shown as dashed lines in Figure~\ref{fig:Marmousi}f.}
\label{fig:slices}
\end{figure*}

In Figures~\ref{fig:MarmRes}a--\ref{fig:MarmRes}c we show the ROM based inversion results obtained from data 
contaminated with $1\%$ additive noise described in Appendix \ref{app:numdata}. We used $\ell = 6$ 
layers in Algorithm \ref{alg:prowi}, with $\q = 3$ iterations per layer, and the restriction parameter $d = 10$. 
The ROM operator is regularized as explained in Appendix~\ref{app:RegROM} with the spectral 
threshold parameter set to $r =\n- 9 = 31$. The velocity is parametrized as in equation \eqref{eq:IM2}, 
with the initial guess $c_o(\bx)$ displayed in Figure~\ref{fig:Marmousi}d.  
We used  $N = 50 \times 30 = 1500$ Gaussian basis functions {defined as in equation} \eqref{eq:Gaussphi}, with standard deviations 
$\sigma_\phi^\perp = 60$~m, and $\sigma_\phi = 56.4$~m. The peaks of the Gaussians are on a uniform 
$50 \times 30$ grid discretizing the imaging domain 
$\Omega_{\rm im} = [103\text{ m}, 5147\text{ m}]\times[97\text{ m}, 2903\text{ m}]$ contained in $ \Omega$.
Figure \ref{fig:MarmRes}d--\ref{fig:MarmRes}f shows the FWI results computed for noiseless data. 
We use the same parametrization of the search velocity and invert in $\ell  = 6$ layers with the same data 
windowing as in the ROM based inversion.

We observe in Figure \ref{fig:MarmRes} that the ROM based {velocity estimation} captures correctly 
many features of the Marmousi model, and continues to improve with the iterations. 
The imaging near the bottom boundary can be improved further by extending the duration of 
the measurements and the depth of the domain $\Omega$, so that the artificial bottom boundary 
has no effect. We also note that the FWI approach recovers the top features of the Marmousi model. 
However, the velocity estimate does not improve much after the $12^{\rm th}$ iteration and the 
result is far from the true model. Effectively, FWI is stuck in a local minimum.

In Figures \ref{fig:Marmousi}b--\ref{fig:Marmousi}c we compare the ROM based estimates obtained with two 
different choices of the basis functions in the parametrization from expression \eqref{eq:IM2} of the search velocity. 
The Gaussian ones given in equation \eqref{eq:Gaussphi} and the commonly used piecewise linear hat functions, 
which interpolate between the values of {zero and one}  on the same $50 \times 30$ inversion grid. 
The estimate with the Gaussian basis looks smoother, as expected, but the point of this comparison 
is to illustrate that the inversion is {very mildly} sensitive to the parametrization of the search velocity, 
once the inversion grid is fixed.  

We show in Figure~\ref{fig:Marmousi}f how the velocity estimation improves if we double the 
number of {colocated sources/receivers} to $\m = 60$, decrease the time sampling to $\tau = 0.0333\text{ s}$ 
and increase $\n$ to $50$, while also setting $r =\n- 17 = 33$.  
The inversion is carried out as above, except that the parametrization of the velocity is with 
$N = 75\times 38 = 2850$ Gaussian functions with $\sigma_\phi^\perp = 40.2$~m, and $\sigma_\phi = 44.8$~m. 
We use the estimate from Figure~\ref{fig:Marmousi}b as an initial guess. Since this initial velocity 
estimate is already very good, it is sufficient to perform $\q=4$ Gauss-Newton iterations for a single layer 
$\ell=1$ using all the available data, i.e., $k_1 = r$. We note that the resulting refined velocity estimate 
sharpens the boundaries of the features and improves their contrast.

To illustrate better the quality of the refined ROM estimate in Figure~\ref{fig:Marmousi}f, 
we display in Figure~\ref{fig:slices} the true and refined estimated velocity for three vertical slices, 
at distances $1.4\text{ km}, 2.8\text{ km}$ and $3.566\text{ km}$. We note again that the reconstruction 
is accurate away from the bottom boundary, where the results can be improved by extending 
the depth of the domain $\Omega$ and the recording time, as explained above.

We end the section with the velocity estimate obtained with the array response matrix estimated 
from towed-streamer type measurements, which is displayed in Figure~\ref{fig:Marmousi}e. 
We observe that this estimate is practically the same as the one in Figure~\ref{fig:Marmousi}b.

\section{Conclusion}

We introduced a novel approach for velocity estimation based on a reduced order model (ROM) of the wave operator.
The ROM is computed from the data gathered by an array of {colocated  sources and receivers}. 
Such data can be approximated in geophysics applications from towed-streamer type measurements. 
No prior information of the medium is used, except for the assumption that the velocity is known in the immediate 
vicinity of the sensors. While the mapping from the data to the ROM is nonlinear, we can compute it using efficient 
numerical linear algebra algorithms. We explain that the ROM is an approximation of the wave operator on a space 
defined by the snapshots of the wavefield at uniformly spaced time {steps}. This space is not known and 
neither is the wave operator. Yet, we can compute its approximation, the ROM, from the data.
We describe the properties of the ROM and formulate a velocity estimation algorithm that minimizes  the ROM misfit.
We also explain how to regularize the ROM in order to mitigate additive noise.
We demonstrate with numerical simulations that the ROM misfit objective function is better than the nonlinear 
least-squares data misfit used in full waveform inversion (FWI). In particular, for a low-dimensional velocity model 
where we can plot the objective functions, we obtain that the ROM misfit objective function 
%is convex, 
{has convexity properties}
while the FWI objective function displays multiple local minima. We present velocity estimation results for two 
well known models where FWI is known to fail in the absence of an excellent initial guess: 
the ``Camembert" model and the Marmousi model.

{Our ROM construction uses that the data matrices have a symmetric inner-product mathematical expression. 
This requirement is the main impediment to having a straightforward extension of  the methodology to more general data acquisition setups, with sources and receivers  at very different locations. We hope to address this open challenge in future research.}

\section*{ACKNOWLEDGMENTS}
We thank the associate editor and the anonymous referees for the valuable feedback on how to improve the exposition. 
This material is based upon research supported in part by the U.S. Office of Naval Research 
under award number N00014-21-1-2370 to Borcea and Mamonov. Borcea, Garnier and Zimmerling 
also acknowledge support from the AFOSR awards FA9550-21-1-0166 and FA9550-22-1-0077.
Zimmerling also acknowledges support from the National Science Foundation under Grant No. 2110265.

\append[app:numdata]{Numerically simulated data}

The data for  the numerical experiments are computed  with a time-domain wave equation 
solver for equations \eqref{eq:I1}--\eqref{eq:I2}, with Laplacian discretized on a uniform grid  with a five 
point finite difference stencil. We use homogeneous Dirichlet boundary conditions at $\partial \Omega$.
The second time derivative is approximated by a three point finite difference scheme, on a fine time grid 
with step $\tau_\cf = \tau / 20$.  Using {equation}~\eqref{eq:defD}, we get the finely 
sampled data $\bD^\cf_\k$, for $k=0,1,\ldots,n_\cf$, where {$n_\cf = 20(2 \n
-1)$}. 

The noisy data are computed as follows. Define 
\begin{equation}
\beta = \frac{b}{\m \sqrt{n_\cf + 1}} \left( \sum_{k=0}^{n_\cf} \| \bD^\cf_\k \|^2_F \right)^{1/2},
\end{equation}
where $b$ is the desired noise level, e.g., $b = 10^{-2}$ for $1\%$ noise. Then, the contaminated
finely sampled data is obtained by adding to $\bD^\cf_\k$ a realization of an $\m \times \m$
random matrix with independent, normally distributed entries with mean zero and standard deviation 
$\beta$ for each $k=1,\ldots,n_\cf$. Since the data at time zero is computed in the known medium 
near the {colocated sources/receivers},  we exclude $k=0$. To simplify notation, hereafter we denote by $\bD^\cf_\k$ 
both the noiseless and the noise contaminated, {finely sampled data}.

We now explain how we compute the 
second derivative data matrices. We begin by extending the finely sampled data evenly in discrete time
to get $\bD^\cfe_j$, $j = -n_\cf,\ldots,n_\cf$, with $\bD^\cf_\k = \bD^\cfe_{\pm k}$, $k=0,1,\ldots,n_\cf$. 
Then, we take the discrete Fourier transform of $\big( \bD^\cfe_j \big)_{j=-n_\cf}^{n_\cf}$ and differentiate 
in the Fourier domain after using a sharp cutoff low-pass filter intended to stabilize the calculation. 
The cutoff frequency is at $\om_o/(2 \pi) + 4 B = 22$~Hz. We take the  inverse Fourier transform to 
obtain $\ddot{\bD}^\cfe_j$, at $j = -n_\cf,\ldots,n_\cf$, the finely sampled second derivative data. 
Finally,  we sub-sample both $\bD^\cfe_j$ and $\ddot{\bD}^\cfe_j$ to get
\begin{equation}
\bD_\k = \bD^\cfe_{20\k}, \quad \ddot{\bD}_\k = \ddot{\bD}^\cfe_{20\k}, \quad \k = 0,1,\ldots,2\n
-1.
\end{equation}

\append[app:Causal]{Causal construction of the ROM}

Here we prove that the upper left $\k\m \times \k\m$ block of $\bA^\RM$, denoted by 
$\big[\bA^\RM\big]_\k$, is the ROM operator computed by Algorithm \ref{alg:arom} from the data subset 
$\{\bD_j,\ddot \bD_j\}_{j=0}^{2\k-2}$, for any $\k = 1, \ldots, \n$.

Let us begin by writing $\big[\bA^\RM\big]_\k$ from equation \eqref{eq:ROMA} 
\begin{align}
\big[\bA^\RM\big]_\k &= \begin{pmatrix} \bI_{\k\m} & {\bf 0} \end{pmatrix} 
\bR^{-T} {\bS} \bR^{-1}  \begin{pmatrix} \bI_{\k\m}  \\ {\bf 0} \end{pmatrix} 
\nonumber \\
&=\begin{pmatrix} \big[\bR \big]_{\k}^{-T} & {\bf 0} \end{pmatrix} \bS
\begin{pmatrix}   \big[ \bR \big]_{\k}^{-1} \\ {\bf 0} \end{pmatrix} \nonumber \\
&= \big[\bR\big]_{\k}^{-T} \big[ \bS \big]_\k  \big[\bR\big]_{\k}^{-1}   , 
\label{eq:B1}
\end{align}
where $\bI_{\k\m}$ is the $\k\m \times \k\m$ identity matrix and 
$\big[\bS \big]_\k  $ and $\big[\bR \big]_{\k}$ are the upper left $\k\m \times \k\m$ blocks 
of $\bS $ and $\bR$, respectively. Here we used that $\bR$ is block  upper triangular, and so is its inverse. 
Moreover, the upper left $\k\m \times \k\m$ block of $\bR^{-1}$ is the same as 
the inverse of $\big[ \bR \big]_{\k}$.

At step 3, Algorithm \ref{alg:arom} computes from  $\{\bD_j,\ddot \bD_j\}_{j=0}^{2\k-2}$ 
the upper left $\k\m \times \k\m$ block of $\bM$, denoted by $\big[\bM\big]_\k$,  and also $\big[ \bS\big]_\k $.  
The Cholesky factorization in equation \eqref{eq:Cholesky} and the block upper triangular
structure of $\bR$ give
\begin{align}
\big[\bM\big]_\k &= \begin{pmatrix} \bI_{k\m} & {\bf 0} \end{pmatrix} 
\bR^T \bR \begin{pmatrix} \bI_{\k\m}  \\ {\bf 0} \end{pmatrix}  = \big[\bR\big]_{\k}^T \big[\bR\big]_{\k}.
\end{align}
This shows that  $\big[\bR\big]_{\k}$ is the Cholesky square root of $\big[\bM\big]_\k$, computed in Algorithm 
\ref{alg:arom}. The result follows from equation \eqref{eq:B1}.

\append[app:algrom]{Algebraic structure of the ROM}

We explain here that the entries of the ROM operator $\bA^\RM$ decay away from the main diagonal, 
which is why we can use the restriction mapping $\mbox{Rest}_{d,\k}$ {defined in equation \ref{eq:restdk}} to reduce the computational 
cost of inversion. Let us write 
\begin{equation}
\bV(\bx) = \big( \bv_0(\bx), \ldots, \bv_{\n-1}(\bx)\big),
\end{equation}  
where $\bv_j(\bx) \in \RR^{1 \times \m}$, for $j = 0, \ldots, \n-1.$
We obtain from equation \eqref{eq:ROMAProj} that the $\m\times \m$ blocks of $\bA^\RM$ are 
\begin{equation}
\bA^\RM_{i,j} = \int_{\Omega} d \bx \, \bv_i^T(\bx) {\cal A} \bv_j(x) , \qquad i,j = 0, \ldots, \n-1.
\label{eq:B2}
\end{equation}
Moreover, the Gram-Schmidt orthogonalization from equation \eqref{eq:GS} gives
\begin{equation}
\bu_j(\bx) = \sum_{q = 0}^{j} \bv_q(\bx) \bR_{q,j},
\end{equation}
and conversely
\begin{equation}
\bv_{j}(\bx) = \sum_{q = 0}^{j} \bu_q(\bx) \bGa_{q,j},
\label{eq:B4}
\end{equation}
where 
\begin{equation}
\bGa  = \bR^{-1} = \begin{pmatrix} \bGa_{0,0} & \bGa_{0,1} & \ldots & \bGa_{0,\n-1} \\
{\bf 0} &  \bGa_{1,1} & \ldots & \bGa_{1,\n-1} \\
\vdots & \vdots & \vdots & \bGa_{\n-1,\n-1} \end{pmatrix}
\label{eq:B3}
\end{equation}
{is block upper triangular,  like  $\bR$}.

Now let  us substitute equation \eqref{eq:B4} into equation \eqref{eq:B2}, to obtain
\begin{align}
\bA^\RM_{i,j} 
%&= \int_{\Omega} d \bx \, \bv_j^T(\bx) A \sum_{\q = 0}^{l} \bu_q(\bx) \bGa_{q,l} \nonumber \\
&= \sum_{q=0}^j \int_{\Omega} d \bx \, \bv_i^T(\bx) {\cal A} \bu_q(\bx) \bGa_{q,j} \nonumber \\
&= -\sum_{q=0}^j \int_{\Omega} d \bx \, \bv_i^T(\bx) \partial_t^2 \bu(q\tau,\bx) \bGa_{q,j}.
\label{eq:B6}
\end{align}
We use next 
the Whittaker-Shannon interpolation formula, which says that if $\tau$ satisfies the Nyquist criterion, then
\begin{equation}
u(t,\bx) = \sum_{s = -\infty}^\infty u_{|s|}(\bx) \mbox{sinc} \Big[ \frac{\pi (t- s\tau)}{\tau}\Big].
\label{eq:WS1}
\end{equation}
Differentiating twice and evaluating at $t = q \tau$, we get
\begin{align}
\tau^2 \partial_t^2 u(q\tau,\bx) =\sum_{s=-\infty, s\neq 0}^\infty  \frac{2(-1)^{s+1}}{s^2}u_{|q-s|}(\bx) 
%\nonumber \\
- \frac{\pi^2}{3} u_q(\bx),
\end{align}
and substituting into equation \eqref{eq:B6}, we obtain
\begin{align}
\nonumber
\bA_{i,j}^\RM & = \frac{1}{\tau^2} \sum_{q=0}^j \bGa_{q,j} 
\Big\{ \int_{\Omega} d \bx \, \bv_i^T(\bx) \bu_{q}(\bx) \\
\nonumber
& \quad -\sum_{s=-\infty, s\neq 0}^\infty  \frac{2(-1)^{s+1}}{s^2} 
\int_{\Omega} d \bx \, \bv_i^T(\bx) \bu_{|q-s|}(\bx)\Big\} \\
& = \frac{1}{\tau^2} \sum_{q=0}^j \bGa_{q,j}  \Big\{\bR_{i,q} 
- \sum_{s=-\infty, s\neq 0}^\infty \frac{2(-1)^{s+1}}{s^2}  \bR_{i,|q-s|}\Big\}.
\label{eq:B8}
\end{align}
To avoid boundary terms,  we have assumed in this formula a large $\n$ so we can take  $\n \to \infty$. 

Since $\bGa_{q,j} = 0$ for $q > j$, and $\bR_{i,q} = 0$ for $i > q$, the first term on the right-hand side 
of equation \eqref{eq:B8} is zero for $i > j$. But we are interested only in the block upper triangular part of 
$\bA^\RM$ (i.e., $i \le j$), due to symmetry, so this first term contributes only to the main block diagonal.
The other block diagonals are due to the series in equation \eqref{eq:B8}. Each term in this series adds an 
$s^{\rm th}$ diagonal, whose entries decay as $1/s^2$. Thus, only the first few block diagonals are large.

\append[app:implem]{Implementation of the inversion}

In principle, the optimization at step 3 of Algorithm \ref{alg:prowi} could have a constraint on 
$\bet$ to ensure that the search velocity in equation \eqref{eq:IM2} is positive. 
We did not need such a constraint in our numerical simulations, 
{as the velocity has stayed positive throughout the iterations}. 

There are many possible regularization penalties. For simplicity, we use the adaptive Tikhonov regularization 
\begin{equation}
\mathcal{L}_i^{\rm reg}(\bet) = \mu_i \|\bet\|_2^2,
\label{eq:regul}
\end{equation}
where $\| \cdot \|_2$ is the Euclidean norm and $\mu_i$ is chosen adaptively 
with the following procedure. Let
\begin{equation}
\mathcal{E}(\bet;d,k_l) = \mbox{Rest}_{d,k_l} \left(\left[\bA^\RM({v}(\cdot;\bet))-\bA^\RM\right]_{k_l}\right)
\label{eq:residual}
\end{equation}
be the $\m d(2\k-d+1)/2$-dimensional  residual vector, whose Euclidean norm squared appears in equation \eqref{eq:newObj}. The  Jacobian of the objective function evaluated at $\bet = \bet^{(i-1)}$ is the matrix
\begin{equation}
\bJ^{(i)} = \nabla_\bet  \mathcal{E}(\bet^{(i-1)};d,k_l) \in \RR^{d\m (\k\m- (d\m-1)/2) \times N}.
\end{equation}
% is a rectangular matrix with $\md(2\k-d+1)/2$ rows and $N$ columns.  
We always choose the parametrization of the velocity from equation \eqref{eq:IM2} so that the Jacobian 
has more rows than columns.
Let $\sigma_1^{(i)} \geq \sigma_2^{(i)} \geq \cdots \geq \sigma_{N}^{(i)}$
be the singular values of $\bJ^{(i)}$. For a fixed parameter $\gamma \in (0, 1)$, 
with smaller values corresponding to stronger regularization, we set
\begin{equation}
\mu_i = \big( \sigma^{(i)}_{\lfloor \gamma N \rfloor} \big)^2.
\end{equation}
The choice of $\gamma$ depends on the parametrization in equation \eqref{eq:IM2}. 
Since it is not clear what is the resolution of the inversion, we choose to over-parametrize 
the velocity, and stabilize the inversion with a small $\gamma$, in the range $(0.2,0.4)$. 
For the results presented in the paper we used $\gamma = 0.25$.

The Gauss-Newton update direction for the objective function {in equation \eqref{eqn:regObj}, 
regularized with $\mathcal{L}_i^{\rm reg}(\bet)$ from equation \eqref{eq:regul}}, is
\begin{equation}
\bd^{(i)} = - \left( \big( \bJ^{(i)} \big)^T \bJ^{(i)} + \mu_i \bI_N \right)^{-1}
\big( \bJ^{(i)} \big)^T \br^{(i)},
\label{eqn:updatedir}
\end{equation}
where $\bI_N$ is the $N \times N$ identity matrix and $\br^{(i)}$ is the residual vector 
in equation \eqref{eq:residual} evaluated at $\bet^{(i-1)}$. {Note that equation \eqref{eqn:updatedir} is the
same as the Levenberg-Marquardt update direction for equation \eqref{eq:residual} with damping $\mu_i$.}

Given the update direction $\bd^{(i)}$, we use a line search
\begin{equation}
\alpha^{(i)} = \mathop{\mbox{argmin}}\limits_{\alpha \in (0, \alpha_{\max})}
\mathcal{L}_i\big( \bet^{(i-1)} + \alpha \bd^{(i)} \big)
\end{equation}
to compute the step length $\alpha^{(i)}$, where we take $\alpha_{\max} = 3$.
Then, the Gauss-Newton update  is
\begin{equation}
\bet^{(i)} = \bet^{(i-1)} + \alpha^{(i)} \bd^{(i)}.
\end{equation}

 We use a similar regularization strategy for the FWI objective function in equation \eqref{eq:FWIobjE}:
If we let $\mathcal{E}^\FWI(\bet)$  be the $\n \m (\m+1)$ dimensional residual vector, with entries 
\begin{equation}
\big(\mathcal{E}_j^\FWI(\bet)\big)_{ j = \k\m(\m+1)/2 +1}^{(\k+1)\m(\m+1)/2} 
= {\rm Triu} \big(\bD_\k({v})-\bD_\k\big), 
\end{equation}
for $k = 0, \ldots, 2 \n
-1$, then  its Jacobian evaluated at $\bet = \bet^{(i-1)}$ is 
%the $\n\m(\m+1) \times N$ matrix
\begin{equation}
\bJ^{\FWI,(i)} = \nabla_\bet  \mathcal{R}^\FWI(\bet^{(i-1)}) \in \RR^{\n\m(\m+1) \times N},
\end{equation}
where we assume $N \leq \n\m(\m+1)$. Furthermore, for the same fixed parameter
$\gamma$ used in the ROM approach,  we set
$\mu_i ^\FWI= \big( \sigma^{\FWI,(i)}_{\lfloor \gamma N \rfloor} \big)^2$,
where $\big\{\sigma_j^{\FWI,(i)}\big\}_{j=1}^N$ are the singular values of $ \bJ^{\FWI,(i)}$, 
sorted in decreasing order.

\append[app:RegROM]{Regularization of the ROM }

Let us denote by $\{\bD_j^\cN\}_{j = 0}^{2\n-1}$ the uncertain data matrices that are either 
contaminated with noise (see Appendix~\ref{app:numdata}) or approximated from the towed-streamer 
measurements. Source-receiver reciprocity is built into the approximation for towed-streamer measurements, 
but it does not hold for noisy array  data. To ensure symmetry, we transform
$\bD_j^\cN$ into $ \frac{1}{2} ( \bD_j^\cN + \bD_j^{\cN^T} ) $.

The mass and stiffness matrices computed at step 3 of Algorithm \ref{alg:arom} are denoted by 
$\bM^\cN$ and $\bS^\cN$. In theory, they should be positive definite matrices, but they will have a 
number of eigenvalues that are negative or zero. This is critical in the case of $\bM^\cN$, 
because we need the inverse of its block Cholesky square root to compute $\bA^\RM$.

A natural way of regularizing $\bM^\cN$ is via projection on the space spanned by the leading eigenvectors. 
Thus, let 
\begin{equation}
\bM^\cN = \bZ^\cN \bLa^\cN (\bZ^\cN)^T
\label{eq:R1}
\end{equation}
be the eigendecomposition of $\bM^\cN$, where $\bZ^\cN$ is the orthogonal matrix of eigenvectors and 
$\bLa^\cN = {\rm diag}(\la_1^\cN, \ldots, \la_{\n\m}^\cN\big)$ is the diagonal matrix of eigenvalues, 
in descending order. We wish to keep the eigenvalues  that are larger than the 
noise contribution (see Appendix~\ref{app:rreg}). Since we work with $\m \times \m$ blocks, 
we choose the cut-off at index $r \m$, for integer $r$ satisfying $1 \le r < \n$, and use the first $r \m$ 
eigenvectors, stored in 
\begin{equation}
\label{eq:restrZ}
%\mathcal{Z} 
 \bZ^{\cN,r}
 = (Z_{jl}^\cN)_{1\leq j \leq \n\m, 1\leq l \leq r\m}
\in \RR^{\n\m \times r\m}
\end{equation}
to define the projected mass matrix
\begin{align}
\label{eq:projMass}
\bLa^{\cN,r} &= 
%\mathcal{Z}^T
( \bZ^{\cN,r})^T
 \bM^\cN 
 %\mathcal{Z}
  \bZ^{\cN,r}
  =  \mbox{diag}\big(\la_1^\cN, \ldots, \la_{r\m}^\cN\big).
\end{align}
The resulting $\bLa^{\cN,r}$ is well-conditioned, but it does not have the block Hankel + Toeplitz 
structure deduced from the causal propagation of the wave (recall equation \eqref{eq:M32}). 
Thus, we need an additional transformation to recover  causality. The desired transformation cannot be obtained by looking at the 
ROM operator construction alone, because all we know about the algebraic structure of 
$\bA^\RM$ is that its entries decay away from the main diagonal. However, we can get the 
transformation using another ROM, for the ``propagator" operator
\cite[]{borcea2018untangling,borcea2021reduced},
\begin{equation}
{\cal P} =  \cos \big(\tau \sqrt{\cA} \big).
\end{equation}

The ROM propagator is obtained from the Galerkin approximation of the time stepping equation
\begin{align}
\bu_{j+1}(\bx) = 2 {\cal P}  \bu_j(\bx) - \bu_{|j-1|}(\bx), \quad j \ge 0,
\label{eqn:timestepp}
\end{align}
obtained from  equation \eqref{eq:tStep} evaluated at 
$t = j \tau$ and $\Delta t = \tau$, for $ j \ge 0$. 
The approximation space is the same as in the computation of $\bA^\RM$, i.e., 
$\mbox{range} \big(\bU(\bx)\big)$. If we let $\tilde \bu_{\rm G,j} = \bU(\bx) \tilde \bg_j$ 
be the Galerkin approximation at {instance}  $t = j \tau$, this satisfies the algebraic system of equations
\begin{align}
\underbrace{\int_{\Omega} d \bx \, \bU^T(\bx) \bU(\bx)}_{\bM} &
\big( \tilde \bg_{j+1} + \tilde \bg_{|j-1|}\big) \nonumber \\
& = 2  \underbrace{ \int_{\Omega} d \bx \, \bU^T(\bx) {\cal P} \bU(\bx)}_{\tilde \bS} \tilde \bg_j,
\end{align}
for $j \ge 0$. Note how the same data driven mass matrix $\bM$ appears in this equation. 
The propagator stiffness matrix $\widetilde \bS$ is also data driven, with $\m \times \m$ blocks given by 
\begin{align}
\widetilde \bS_{i,j} &= \langle \bu_i, \cos \big(\tau \sqrt{\cA}\big) \bu_j \rangle  \nonumber \\
&= \frac{1}{2} \langle \bu_i, \bu_{j+1} + \bu_{|j-1|} \rangle \nonumber \\
&= \frac{1}{4} \left( \bD_{i+j+1} + \bD_{|i-j-1|} + \bD_{|i+j-1|} + \bD_{|i-j+1|} \right),
\end{align}
{where $ 0 \leq i,j \leq \n-1$}. Using the same block Cholesky factorization of $\bM$ from equation \eqref{eq:Cholesky} 
 and multiplying the Galerkin equation above by $\bR^{-1}$, we get the time stepping 
scheme in the ROM space
\begin{equation}
\tilde \bu_{j+1}^\RM = 2 \cP^\RM \tilde \bu_{j}^\RM - \tilde \bu^\RM_{|j-1|}, 
\end{equation}
{the algebraic analogue of equation \eqref{eqn:timestepp}}, where 
\begin{equation}
\tilde \bu_{j}^\RM = \bR \tilde \bg_j, \quad j \ge 0,
\end{equation}
are the ROM snapshots and 
\begin{equation}
\cP^\RM = \bR^{-T} \widetilde{\bS} \bR^{-1} = \int_{\Omega} d \bx \bV^T(\bx) {\cal P}  \bV(\bx),
\end{equation}
is the ROM propagator.  In the last equality we used the definition of $\widetilde \bS$ 
and the Gram-Schmidt orthogonalization equation \eqref{eq:GS}.

We refer the interested reader to \cite[]{borcea2020reduced} for a long and detailed analysis of $\cP^\RM$. 
For our purpose, it suffices to say that it is useful to look at it because, as proved in 
\cite[~Appendix C]{borcea2020reduced}, unlike $\bA^\RM$, the matrix $\cP^\RM$ is sparse, 
with block tridiagonal structure. Moreover, $\cP^\RM$ is determined by the same mass matrix as $\bA^\RM$. 
Thus, even though the regularized matrix in equation \eqref{eq:projMass}  is not in the right algebraic form, 
we can bring it in the right form by imposing the block tridiagonal structure of the resulting ROM propagator.
To do this, we use the block-Lanczos algorithm \cite[~Chapter 10]{golubVanLoan} 
that takes any symmetric matrix in $\RR^{\n\m \times \n\m}$ and computes an orthogonal basis of 
$\RR^{\n\m}$ that puts the matrix in block tridiagonal form. 

We can now describe the regularization procedure.  First, we compute the ROM propagator stiffness matrix 
$\widetilde \bS^\cN$, with blocks given as above, in terms of the uncertain data 
$\big\{\bD_j^\cN \big\}_{j=0}^{2\n -1}$. Then, we project this matrix onto the range of 
$\bZ^{\cN,r}$, defined in equation \eqref{eq:restrZ},
\begin{equation}
\widetilde \bS^{\cN,r} = ( \bZ^{\cN,r})^T \widetilde \bS^\cN \bZ^{\cN,r} \in \RR^{r\m \times r\m},
\label{eq:projStiff}
\end{equation}
and we compute 
\begin{equation}
\bP^{\cN,r} = (\bLa^{\cN,r})^{-1/2} \widetilde \bS^{\cN,r} (\bLa^{\cN,r})^{-1/2}  \in \RR^{r\m \times r\m}.
\label{eq:projROMprop}
\end{equation}
This is a symmetric, positive definite matrix that we put in block tridiagonal form using the 
block-Lanczos algorithm \cite[~Chapter 10]{golubVanLoan}, with starting 
{block $(\bLa^{\cN,r})^{-1/2} (\bZ^{\cN,r})^T \be_0 \in \RR^{r\m \times \m}$}.
This generates an orthogonal matrix $\bQ^{\cN,r} \in \RR^{r\m \times r\m}$ such that 
\begin{equation}
\cP^{\RM,r} = (\bQ^{\cN,r})^T \bP^{\cN,r} \bQ^{\cN,r}  \in \RR^{r\m \times r\m}
\label{eq:noisyROMProp}
\end{equation}
is a block tridiagonal matrix with $\m \times \m$ blocks, which we call the regularized ROM propagator.

The matrix $ \cP^{\RM,r}$ itself is irrelevant for our velocity estimation approach. It is the orthogonal
transformation given by $\bQ^{\cN,r}$ that we need, which restores the desired algebraic causality 
of the regularized mass matrix. Using this transformation we can obtain the regularized ROM operator 
with the following procedure: Compute the block Cholesky factorization of the transformed mass matrix
\begin{align}
(\bQ^{\cN,r})^T  \bLa^{\cN,r} \bQ^{\cN,r} & 
= (\boldsymbol{\Pi}^{\cN,r})^T  \bM^\cN \boldsymbol{\Pi}^{\cN,r} =  (\bR^{\cN,r})^T \bR^{\cN,r},
\label{eq:noisyR}
\end{align}
where 
\begin{equation}
\boldsymbol{\Pi}^{\cN,r} = \bZ^{\cN,r}  \bQ^{\cN,r} \in \RR^{\n\m \times r\m} .
\label{eq:defPi}
\end{equation}
and $\bR^{\cN,r} \in \RR^{r\m \times r\m}$ is block upper triangular and well conditioned, 
due to the spectral truncation in equation \eqref{eq:projMass}. 
Then, using the data driven stiffness matrix $\bS^\cN$ computed at step 3 of Algorithm \ref{alg:arom},
we obtain the regularized operator ROM as 
\begin{equation}
\bA^{\RM,r}= (\bR^{\cN,r})^{-T}  (\boldsymbol{\Pi}^{\cN,r})^T \bS^\cN  \boldsymbol{\Pi}^{\cN,r}  (\bR^{\cN,r})^{-1} .
\label{eq:noisyAROM}
\end{equation}

Equation \eqref{eq:noisyAROM} gives the regularization of the data driven ROM operator construction. 
For the inversion, we also need the ROM operator for the search velocity ${v}(\bx;\bet)$ computed 
via the same chain of transformations, using the same matrix from equation \eqref{eq:defPi}: 
Let $\bM({v})$ and $\bS({v})$ be the mass and stiffness matrices calculated as in 
step~3 of Algorithm~\ref{alg:arom} from the data computed numerically in the medium with 
velocity ${v}(\bx,\bet)$. We compute the block Cholesky factorization
\begin{equation}
(\boldsymbol{\Pi}^{\cN,r})^T  \bM({v})  \boldsymbol{\Pi}^{\cN,r} = \bR^r({v})^T \bR^r({v}),
\label{eq:noisyRV}
\end{equation}
where $r$ is an index (not a power).
Then,  the ROM operator at the search velocity $v$ is given by
\begin{equation}
\bA^{\RM,r}({v}) = 
\bR^r({v})^{-T} (\boldsymbol{\Pi}^{\cN,r})^T \bS({v}) \boldsymbol{\Pi}^{\cN,r} \bR^r({v})^{-1}.
\label{eq:noisyAROMV}
\end{equation}
The velocity inversion is carried out as in Algorithm~\ref{alg:prowi}, 
with $\bA^\RM$ and $\bA^\RM(v)$ in  equation \eqref{eq:newObj}
replaced by the regularized $\bA^{\RM,r}$ and $\bA^{\RM,r}({v})$.
Note that the matrix with orthogonal columns $\boldsymbol{\Pi}^{\cN,r}$ used in 
equations \eqref{eq:noisyRV} and \eqref{eq:noisyAROMV} is computed once using the uncertain data
and does not change over the course of velocity estimation.

\begin{figure*}[h!]
\begin{center}
\begin{tabular}{cc }
(a) Log of objective function, $r = \n$ 
& (b) Log of objective function, $r = \n-4$ \\ 
\includegraphics[width=0.32\textwidth]
{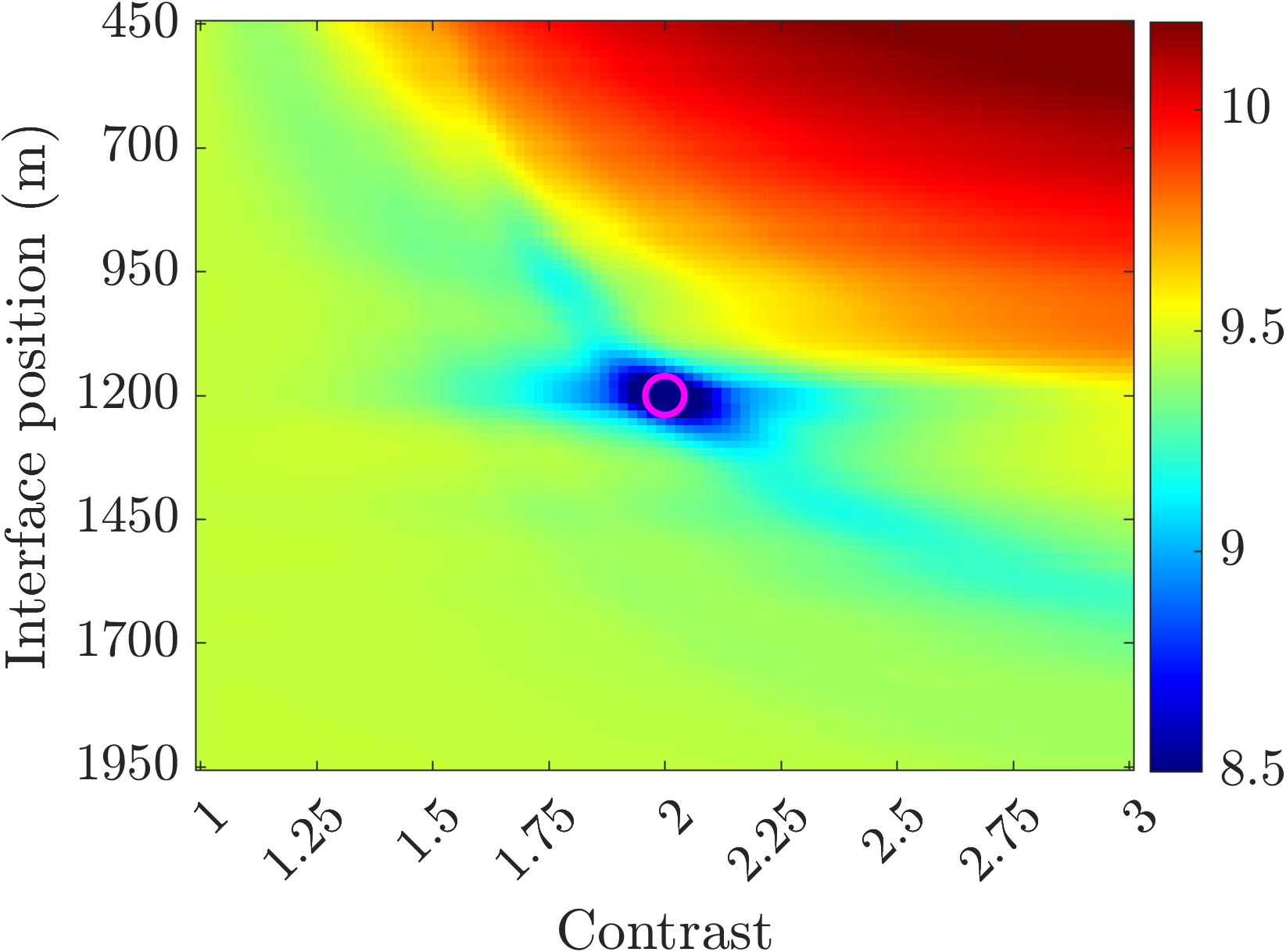} 
&
\includegraphics[width=0.32\textwidth]
{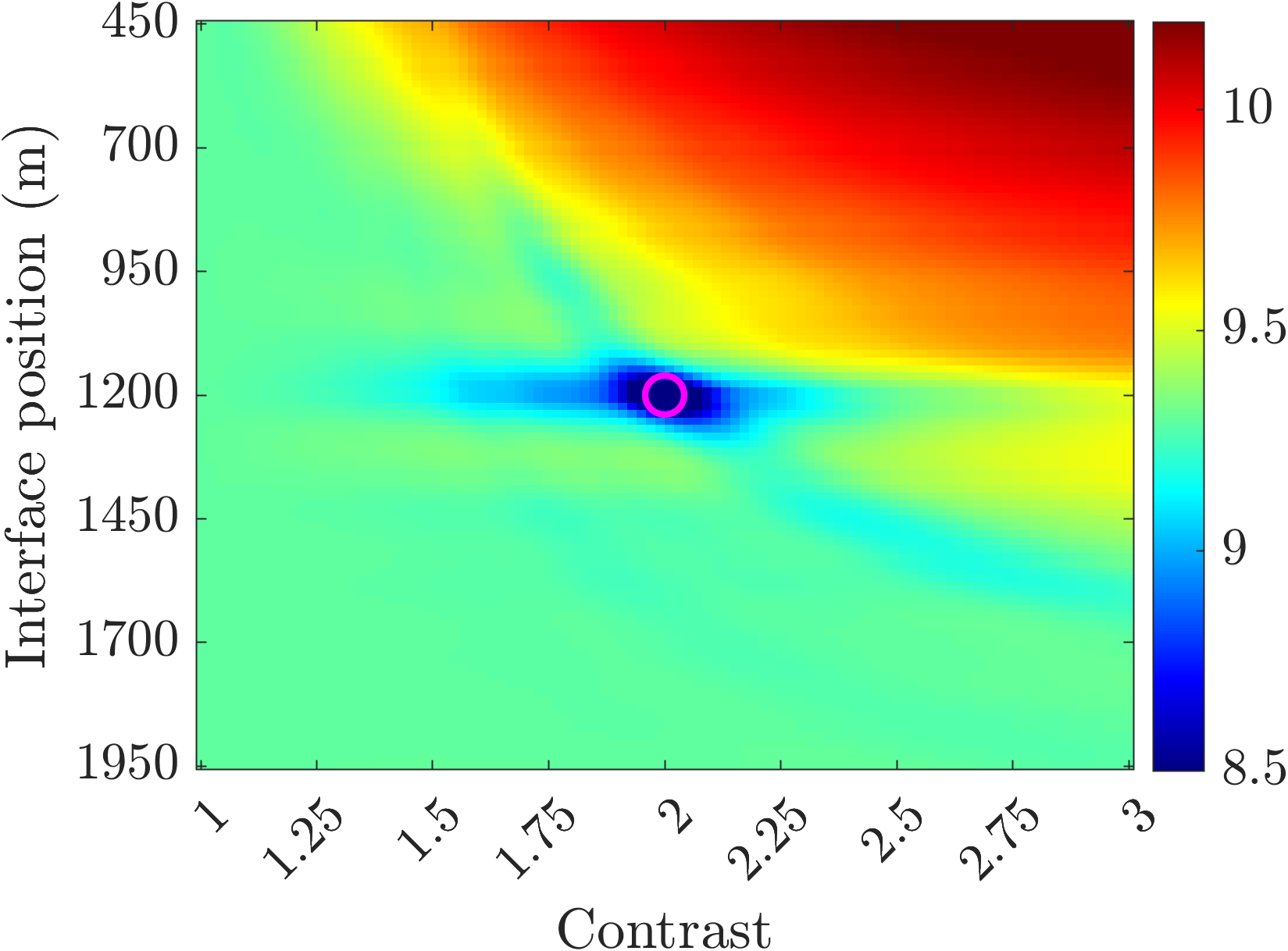} 
\end{tabular}
\end{center}
\vspace{-0.2in}
\caption{Decimal logarithm of objective function (equation \eqref{eq:ROMObjr}) vs. the interface position and velocity contrast.
The true parameters (shown in Figure \ref{fig:topo}) are indicated by $\textcolor{magenta}{\bigcirc}$.}
\label{fig:topog}
\end{figure*}

We observe that due to the block algebra, even if we do not use a spectral truncation, i.e., set $r = \n$, 
the ROM operator in equation \eqref{eq:noisyRV} is not identical to the one computed with Algorithm \ref{alg:arom}.
Nevertheless, they behave the same with respect to the inversion, as illustrated in Figure \ref{fig:topog},
where we plot the logarithm of the objective function 
\begin{equation}
\mathcal{O}^{\RM,r}({v}) = \left\|{\rm Triu} \big(\bA^{\RM,r}({v})-\bA^{\RM,r}\big)\right\|_2^2
\label{eq:ROMObjr}
\end{equation}
for the same experiment as in Figure \ref{fig:topo}, for the cases $r = \n$ and $r = \n-4$. 
There is little difference between Figure \ref{fig:topo}c and Figure \ref{fig:topog}a-\ref{fig:topog}b.

\append[app:rreg]{The regularization threshold}

\begin{figure*}[h!]
\begin{center}
\begin{tabular}{cc}
(a) Singular values of mass matrices &(b) $\left| \sigma_j^\cN / \sigma_j^o - 1 \right|$ and $\varepsilon_\sigma$ \\
\hspace{-0.12\columnwidth}\includegraphics[width=0.33\textwidth] 
{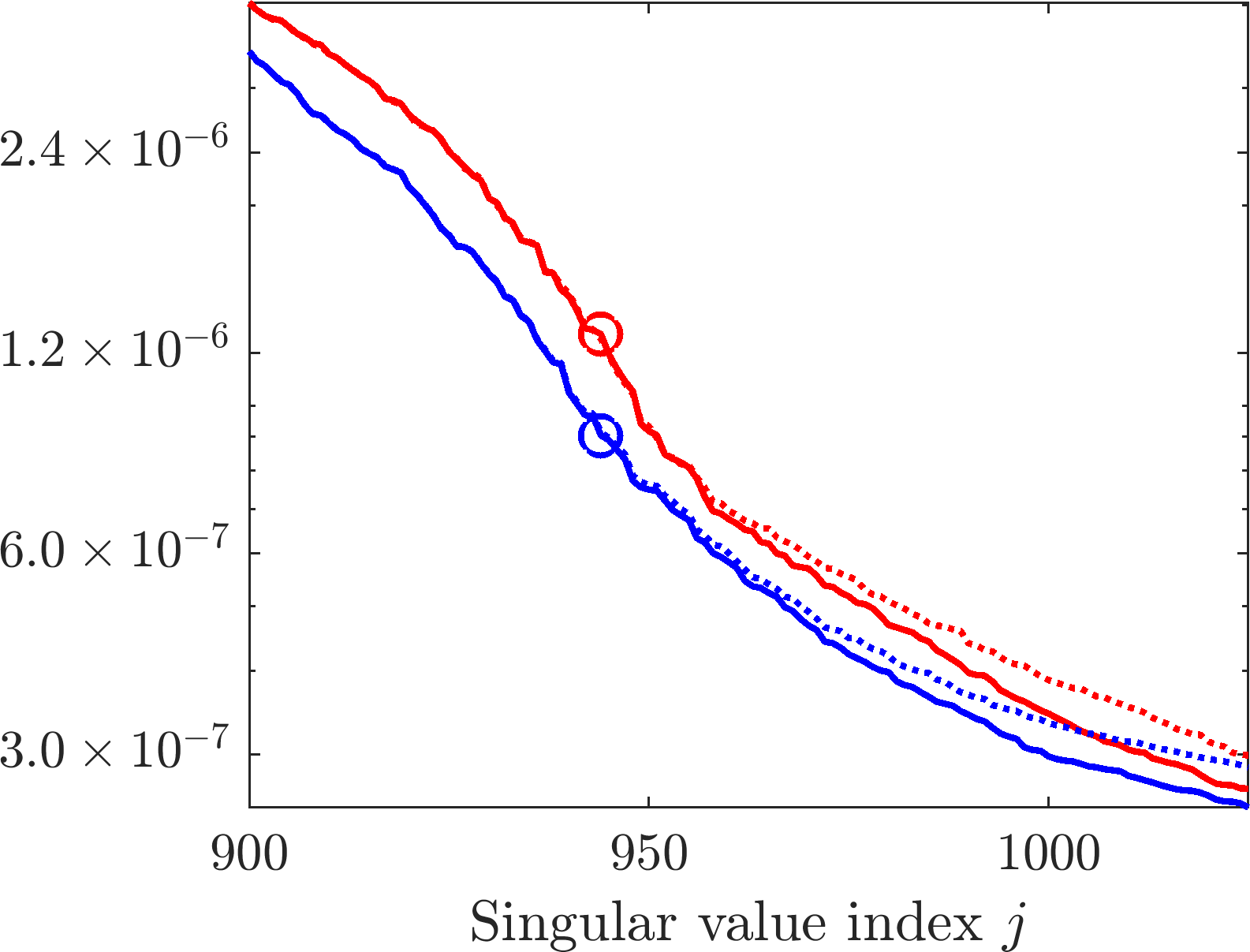} &
\includegraphics[width=0.33\textwidth] 
{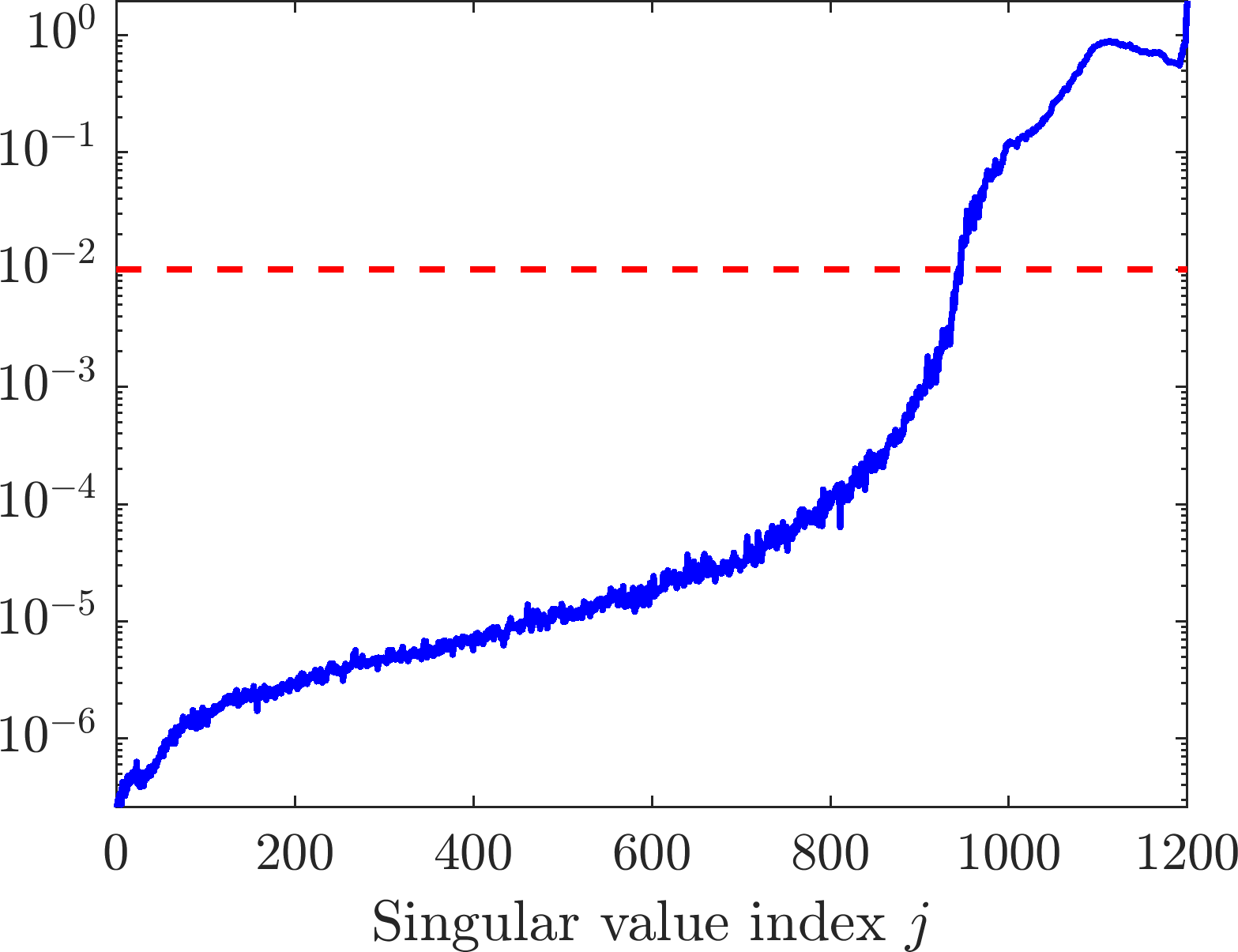}
\end{tabular}
\end{center}
\vspace{-0.2in}
\caption{Regularization threshold illustration:
(a) Singular values of mass matrices $\bM$ (solid red), $\bM^{{\protect\cN}}$ 
(dotted red), $\bM(c_o)$ (solid blue) and $\bM^{{\protect\cN}}(c_o)$ (dotted blue). 
The circles correspond to $j = R^{{\protect\cN}}$;
(b) Left-hand side of equation~\ref{eqn:singrat} (solid blue) and $\varepsilon_\sigma$ (dashed red).}
\label{fig:rreg}
\end{figure*}
Here we explain how we choose the regularization threshold $r$ for the ROM regularization 
procedure in equations \eqref{eq:restrZ} and \eqref{eq:noisyAROMV}. The idea is that $r$ can be determined 
from the part of the spectrum of the mass matrix $\bM^\cN$ that is perturbed by the uncertainty. 
This can be estimated using the mass matrix $\bM^\cN(c_o)$ corresponding to the initial guess 
velocity $c_o(\bx)$, and perturbed in a similar way.
%: For the towed-streamer  type measurements, this means that we use the same approximation procedure. 

{
For noisy sensor array measurements, 
 the matrices
\begin{equation}
% \tag{D-1}
{\itbf E}_j^\cN = \frac{1}{\sqrt{2}} \left( \bD_j^\cN - (\bD_j^\cN)^T \right),
\quad j = 0,\ldots,2 \n
-1,
\end{equation}
can be considered as realizations of the additive noise.  This is because the true wave signals 
are reciprocal ($\bD_j$ are symmetric matrices),  while the additive noise is not. For towed-streamer data other 
noise estimation procedure is needed. For example, measurements at $\m' \gg \m$ locations, separated by a small distance with respect to the wavelength, can be used to estimate the noise at the $\m$ receivers used in the ROM construction. }

Consider the mass matrices $\bM(c_o)$ and $\bM^\cN(c_o)$ computed by Algorithm~\ref{alg:arom} 
from the noiseless background data $\{\bD_j(c_o)\}_{j = 0}^{2\n -1}$ and the artificially generated 
contaminated background data $\{ \bD_j(c_o) + {\itbf E}_j^\cN \}_{j = 0}^{2\n -1}$, respectively.
Let $\{\sigma_j^o\}_{j=1}^{\n\m}$ be the singular values of $\bM(c_o)$, and
$\{\sigma_j^\cN\}_{j=1}^{\n\m}$ the singular values of $\bM^\cN(c_o)$, sorted in decreasing order.
Choose a small $\varepsilon_\sigma$, the largest relative deviation of singular values past 
which we consider them contaminated by noise. Let $R^\cN$ be the smallest among $j$ such that
\begin{equation}
% \tag{D-2}
\Big| \frac{\sigma_j^\cN}{\sigma_j^o} - 1 \Big| \geq \varepsilon_\sigma.
\label{eqn:singrat}
\end{equation}
Then, we can estimate $r = \lfloor R^\cN/ \m \rfloor$.

Note that the estimation can be adaptive. We can choose at iteration $i$ in Algorithm~\ref{alg:prowi} 
the value $r_i$ obtained as above but with $\bM(v(\cdot; \bet^{(i)}))$ instead of $\bM(c_o)$.
However, in our examples this was not necessary, since using $\bM(c_o)$ provided a 
robust if somewhat conservative estimate, as shown in the numerical example described below.
%In our examples, the estimate $r$ was stable along the iterations.

In Figure~\ref{fig:rreg} we illustrate the choice of regularization threshold for the Marmousi model
in the setting outlined in the numerical results section ($\m= 30$, $\n=40$, $1\%$ additive noise). 
Figure~\ref{fig:rreg}a shows the singular values $\sigma_j^o$ and $\sigma_j^\cN$ for a range 
$j = 900, 901, \ldots, 1025$, while also comparing them to the singular values of $\bM$ and $\bM^\cN$. 
Setting $\varepsilon_\sigma = 10^{-2}$, we obtain $R^\cN = 944$ from equation \eqref{eqn:singrat},
as shown in Figure \ref{fig:rreg}b. This gives the value $r = \lfloor 944 / 30 \rfloor = 31$ used in the 
numerical experiments. Note that this process estimates well the point after which the singular values of 
$\bM^\cN$ diverge from those of $\bM$, as observed in Figure \ref{fig:rreg}a.

\bibliographystyle{seg}  % style file is seg.bst
\bibliography{biblio}

\end{document}